\documentclass[lettersize,onecolumn,12pt]{IEEEtran}
\usepackage{amsmath, amssymb, graphicx,amsthm}
\usepackage{hyperref}
\usepackage{pgfplots}
\usepackage{pgfplotstable}
\usepgfplotslibrary{fillbetween}
\usetikzlibrary{patterns}
\usepackage{booktabs}
\usepackage{multirow}
\usepackage{cancel}
\usepackage{xcolor}
\definecolor{process}{RGB}{0,102,204}
\definecolor{approx}{RGB}{204,0,0}

\newtheorem{thm}{Theorem}

\newtheorem{definition}{Definition}
\usepackage{cite}

\title{ Mean-Field-Type Game Theory \\  with  Rosenblatt Noise} 

\author{ Hamidou Tembine, \IEEEmembership{Senior Member, IEEE}, Tyrone E. Duncan  \IEEEmembership{Life Fellow Member, IEEE} and  Bozenna Pasik-Duncan, \IEEEmembership{Life Fellow Member, IEEE}
\thanks{H. Tembine is  with Department of Electrical Engineering and Computer Science, School of Engineering  at  Quebec University at Trois-Rivi{e}res UQTR, Quebec, Canada, and with Timadie, Guinaga, Grabal,  AI Mali,   TF,  WETE,  SK1 Sogoloton, MFTG, LnG Lab, CI4SI;  (e-mail: {\tt\scriptsize tembine@ieee.org}).}
\thanks{ Tyrone E. Duncan and Bozenna Pasik-Duncan are with the Department of Mathematics, University of Kansas. Nichols Hall,
2335 Irving Hill Rd, Lawrence, KS 66045}

}

\IEEEoverridecommandlockouts

\begin{document}

\maketitle

\begin{abstract}
We study  the integration of Rosenblatt noise into stochastic systems, control theory, and mean-field-type game theory, addressing the limitations of traditional Gaussian and Markovian models. Empirical evidence from various domains, including water demand,  e-commerce, power grid operations, wireless channels, and agricultural supply chains, demonstrates the prevalence of non-Gaussian characteristics such as skews, heavy tails and strong long-range dependencies. The Rosenblatt process, a non-Gaussian non-Markovian, self-similar process, offers a  baseline framework for capturing some the behaviors observed in real data. We develop novel stochastic calculus formulas for Rosenblatt processes, apply these to dynamical systems, and analyze optimal control problems, revealing the suboptimality of traditional noise approximation methods. We extend game-theoretic analysis to environments driven by Rosenblatt noise, establishing conditions for saddle-point equilibria in zero-sum games and identifying state-feedback Nash equilibria in non-zero-sum games. Our findings underscore the importance of incorporating non-Gaussian noise into predictive analytics and control strategies, enhancing the accuracy and robustness of models in real-world applications. These findings represent a significant advancement in mean-field-type game theory with variance-awareness, offering new insights and tools for managing interactive systems influenced by Rosenblatt noise.
\end{abstract}


%\tableofcontents

\section{Introduction}

A Gaussian distribution, also known as the Normal distribution, is a continuous probability distribution characterized by its bell-shaped curve, with well-defined moments (mean and variance), where the tails decay exponentially at a rate determined by the standard deviation, no skew and is symmetric about its mean.  In 1785 Jan Ingenhousz, a Dutchman, discussed the motion of coal dust
particles on the surface of alcohol \cite{Ingenhousz1785}.   Robert Brown was a British botanist. Brown explored fauna in Australia for 3.5 years collecting 3400 specimens of
which 2000 were previously unknown. He returned to England and recorded
these results in several volumes of scientific study. In 1828 he published a paper \cite{Brown1828}, noting a
strange motion of particles in a liquid. He did not have an explanation for
this phenomenon. Subsequently the experiment was repeated by him and
others. They changed the liquid and noted the same type of motion. They
noted that the motion was slowed in a more viscous liquid. They also noted
that the particles were not receiving the energy of motion from the light,
e.g. they performed the experiment in the dark or with little light.
The initial experiment was related to his study of the physiology of plants as
well as classifying them. He considered pollen particles from the plant
species Clarkia Pulchella. A physical explanation for Brown’s observations was given by Einstein
(1908) who showed that the mean square distance that the particles travel
is a linear function of time with the rate depending on the temperature, the
drag coefficient and Boltzmann’s constant \cite{Einstein1908}.  Benoît Mandelbrot's mathematical career was deeply influenced by his time in Paris, most notably at the École Polytechnique, where he developed a unique perspective that diverged from traditional mathematical approaches . He is best known for his pioneering work on fractal geometry, a field he largely created, coining the term "fractal" \cite{Mandelbrot1982}. Concurrently, Harold Hurst, studying at Oxford, made significant contributions to the analysis of long-term dependence in time series, through the Hurst exponent \cite{Hurst1951}, a concept that has relevance to fractal analysis and the understanding of nonlinear systems. While their specific areas of focus differed, both Mandelbrot's and Hurst's work contributed significantly to the broader field of nonlinear mathematics and the analysis of complex, non-linear systems, revealing the limitations of traditional linear models in describing many natural and man-made phenomena.

Recent empirical studies across diverse sectors have demonstrated that classical Gaussian assumptions are often inadequate for describing stochastic phenomena encountered in engineered and socio-economic systems. Traditional models in systems theory, control, and game theory typically rely on the central limit theorem, presuming that uncertainties and fluctuations are normally distributed. However, mounting evidence from various fields indicates that many systems exhibit non-central limit, skews, heavy-tailed distributions, multimodal behaviors, and long-range dependencies that contradict these assumptions.
For instance, in the  e-commerce sector, the Whale Effect is characterized by the disproportionate impact of the top 5\% of customers, who contribute up to 50\% of total revenue. This phenomenon suggests that revenue distributions are heavily skewed and possess significant kurtosis, challenging the symmetry inherent in Gaussian models. In power grid operations, the interplay between solar surplus and evening demand produces dual net load peaks - features that standard Gaussian forecasts fail to capture. 
It is not a direct mixture of two Gaussian distribution as  right-skews are observed empirically. Similarly, agricultural supply chains, as observed by the startup Guinaga in mango production across Mali and Burkina Faso, experience dual harvest seasons that result in approximately 40\% waste due to intervening periods of scarcity. Such dual-modality in supply dynamics further emphasizes the limitations of Gaussian approximations.
Additional evidence of non-Gaussian characteristics is evident in the temporal patterns of hardware demand cycles. The distinct adoption phases for  Central Processing Unit (CPU), Graphics Processing Units (GPU), and Tensor Processing Units (TPU) resources, driven by episodic shocks from cryptocurrency and machine intelligence  sectors, result in multi-peaked demand curves. In the financial domain, the cyclic nature of machine intelligence funding - marked by rapid expansion phases followed by sharp contractions - deviates markedly from the unimodal trends predicted by normal distributions. Moreover, the volatile, sawtooth pattern of GPU pricing, influenced by abrupt surges in demand linked to crypto mining and large learning model developments, illustrates pricing dynamics with abrupt discontinuities and heavy tails. Other applications, such as underwater wireless communication and seasonal commodity pricing - as exemplified by onion prices in regions like Dogon Country and Ansongo - further underscore the prevalence of non-Gaussian fluctuations in real-world data.
The inadequacy of Gaussian models in these settings necessitates the development of alternative mathematical frameworks. The Rosenblatt process, a prototypical example of a non-Gaussian self-similar process with long-range dependence, offers a rigorous tool for capturing higher-order dependencies and memory effects. Unlike Gaussian processes, which are completely determined by their first two moments, the Rosenblatt process incorporates cumulants of higher orders, thereby providing a more faithful representation of the stochastic properties observed in the aforementioned applications.
Integrating the Rosenblatt process into systems and control theory enables the formulation of more accurate models for uncertainty propagation and feedback control. In particular, when extended to mean-field type game theory, this approach facilitates the study of strategic interactions among multiple agents (including machine intelligence agents) whose behaviors are influenced by non-Gaussian noise. Direct methods used to solve mean-field-type games  allow for the reduction of high-dimensional agent-based models to tractable convex completion blocks, while the incorporation of non-Gaussian statistics ensures that these models remain robust in the presence of extreme events and irregular fluctuations.

Mean-field-type game (MFTG) theory \cite{bc3w,bc3wt,bc3tw} investigates the interactions of multiple agents including  diverse entities such as individuals, animal populations, networked devices, corporations, nations, genetic systems, and others.  Decision-makers within MFTG frameworks can be classified as atomic, non-atomic, or a hybrid of both. The nature of these interactions spans a spectrum from fully cooperative to entirely selfish, including intermediate forms such as partially altruistic, partially cooperative, selfless, co-opetitive, spiteful, and self-abnegating behaviors, as well as combinations thereof. A defining characteristic of MFTG theory is the incorporation of higher-order performance criteria, including variance, quantiles, inverse quantiles, skewness, kurtosis, and other risk measures. Critically, the integration of these measures is not necessarily linear with respect to the probability distribution of individual/common state or action variables. MFTG has demonstrated applicability in diverse domains, including building evacuation, energy systems, next-generation wireless networks, meta-learning, transportation systems, epidemiology, predictive maintenance, network selection, and blockchain technologies. In contrast to classical mean-field games, MFTG does not intrinsically necessitate a large population of agents. Furthermore, MFTG exhibits distinct equilibrium structures and associated Master Adjoint Systems (MASS), attributable, in part, to the non-linearity of expected payoffs with respect to the individual mean-field (the probability measure of an agent's own state).  Semi-explicit solutions for a broad class of MFTGs with non-quadratic cost functions have been derived, facilitating the validation of numerical methods and deep learning algorithms, including transformer-based architectures.

The synthesis of systems theory, control engineering, and MFTG \cite{bc3w,bc3wt,bc3tw,rosbib00mix,basar2024foundations,basar2024applications}  under the auspices of the Rosenblatt process represents a significant advancement in the modeling of stochastic dynamics. This framework not only addresses the shortcomings of Gaussian assumptions but also provides a unified methodology for analyzing systems subject to heavy-tailed and multimodal disturbances. In the subsequent sections, we develop the mathematical underpinnings of the Rosenblatt process and delineate its implications for control design and strategic decision-making in multi-agent systems. We further validate our theoretical constructs through case studies that encompass e-commerce revenue distributions, power grid load management, agricultural supply chain optimization, hardware demand forecasting, and commodity pricing analysis. Through these investigations, we demonstrate that embracing non-Gaussian non-Poissonian non-Markovian processes yields enhanced predictive accuracy and improved control performance in systems characterized by Rosenblatt stochastic behavior. Non-central limits are neither anomalies nor irregularities; they are intrinsic to the stochastic processes that shape daily life.

Our contributions are as follows.  We derive novel stochastic calculus formulas for a range of Rosenblatt processes, extending classical techniques to account for the higher-order dependencies and non-Gaussian characteristics inherent in these processes.
We develop dynamical systems driven by Rosenblatt processes and apply these formulations to the Rosenblatt Ornstein-Uhlenbeck (R-OU) process. This framework underpins the design of super diffusion transformers, offering a rigorous foundation for advancements in machine intelligence.
 We analyze optimal control problems under Rosenblatt dynamics, demonstrating Turnpike properties that characterize the asymptotic behavior of optimal trajectories and provide insights into long-term control strategies.
Our work reveals the suboptimality of traditional noise approximation methods in non-Gaussian settings, emphasizing the need for accurate modeling techniques that capture the true statistical behavior of the system.
We extend game-theoretic analysis to environments driven by Rosenblatt processes by formulating and solving semi-explicitly zero-sum games, establishing conditions for the existence  of saddle-point equilibria in adversarial scenarios.
We introduce variance-aware game formulations that incorporate the unique statistical features of Rosenblatt processes, enabling risk-aware decision-making in competitive settings.
We  establish the existence of mean-field-type equilibria in systems driven by Rosenblatt processes under suitable conditions in both the state dynamics and the associated strategic actions of all the agents.

\section{Evidence of Non-Gaussianity in Real Data}

Murray Rosenblatt's 1961 work, "Independence and dependence," significantly contributed to the development of understanding processes with long-range dependence, leading to the formalization of the Rosenblatt process \cite{Rosenblatt1961}. This process emerged from his exploration of limit theorems for sums of dependent random variables, where traditional central limit theorems, which assume independence, fail. Rosenblatt's research highlighted that in certain long-memory scenarios, the limiting distribution is not Gaussian. Instead, he identified a non-Gaussian limit, which became known as the Rosenblatt process. This process is characterized by its self-similarity and non-Gaussian nature, making it crucial for modeling phenomena where long-range dependencies and non-linear behaviors are observed. His findings expanded the scope of stochastic process analysis, providing a tool for modeling  systems in various fields, including finance, hydrology, and telecommunications, where traditional Gaussian models are inadequate \cite{Rosenblatt1972,Rosenblatt1974,Rosenblatt1985,Taqqu1975,DobrushinMajor1979,Taqqu1979,FoxTaqqu1987,AvramTaqqu1987,Taqqu1986,SamorodnitskyTaqqu1994}.

Pawel Domanski studied the noise in control systems based on data from
several hundreds of control loops operating in different process industries
located in several sites all over the world. This data showed that the
theoretical assumption of Gaussian properties for the data is hardly ever
satisfied \cite{Dom15}.

\subsection{e-Commerce: The Whale Effect in Customer Spending }
In e-commerce, purchasing behavior rarely follows a Gaussian distribution. Instead, customer spending often adheres to a power-law, beta, gamma, Rosenblatt shape distribution, where a small fraction of users ("whales") account for a disproportionately large share of revenue. This phenomenon, known as the "whale effect," is a quintessential example of non-Gaussianity in action.  

Why It is Non-Gaussian? 
A Gaussian distribution assumes most data points cluster around the mean, with  no skew, symmetrical tails. However, in e-commerce:  
\begin{itemize}\item Most users spend little to nothing: The majority of customers make infrequent or low-value purchases (e.g., one-time buyers).  
\item A tiny fraction drives revenue: The top 1-5\% of customers ("whales") may generate 30-50\% of total sales. These high-value customers buy luxury goods, bulk orders, or premium subscriptions repeatedly.  
\end{itemize}
For instance, a luxury fashion platform might see 95\% of users spend under \$500 annually, while 5\% spend over \$10,000. This creates a long right tail in the spending histogram, starkly contrasting the bell-shaped Gaussian curve.  There are several business implications of this situation. A heavy reliance on high-spending customers exposes businesses to revenue volatility, as the loss of just a few top spenders can significantly impact financial performance. To mitigate this risk, companies must deploy highly personalized marketing strategies, such as exclusive VIP programs, to maintain engagement and loyalty among their most valuable customers. Additionally, inventory decisions become increasingly strategic, with a focus on high-margin, niche products that cater to the preferences of these top-tier buyers, ensuring sustained profitability and competitive advantage.

\begin{figure}[htb]
\centering
\includegraphics[scale=0.3]{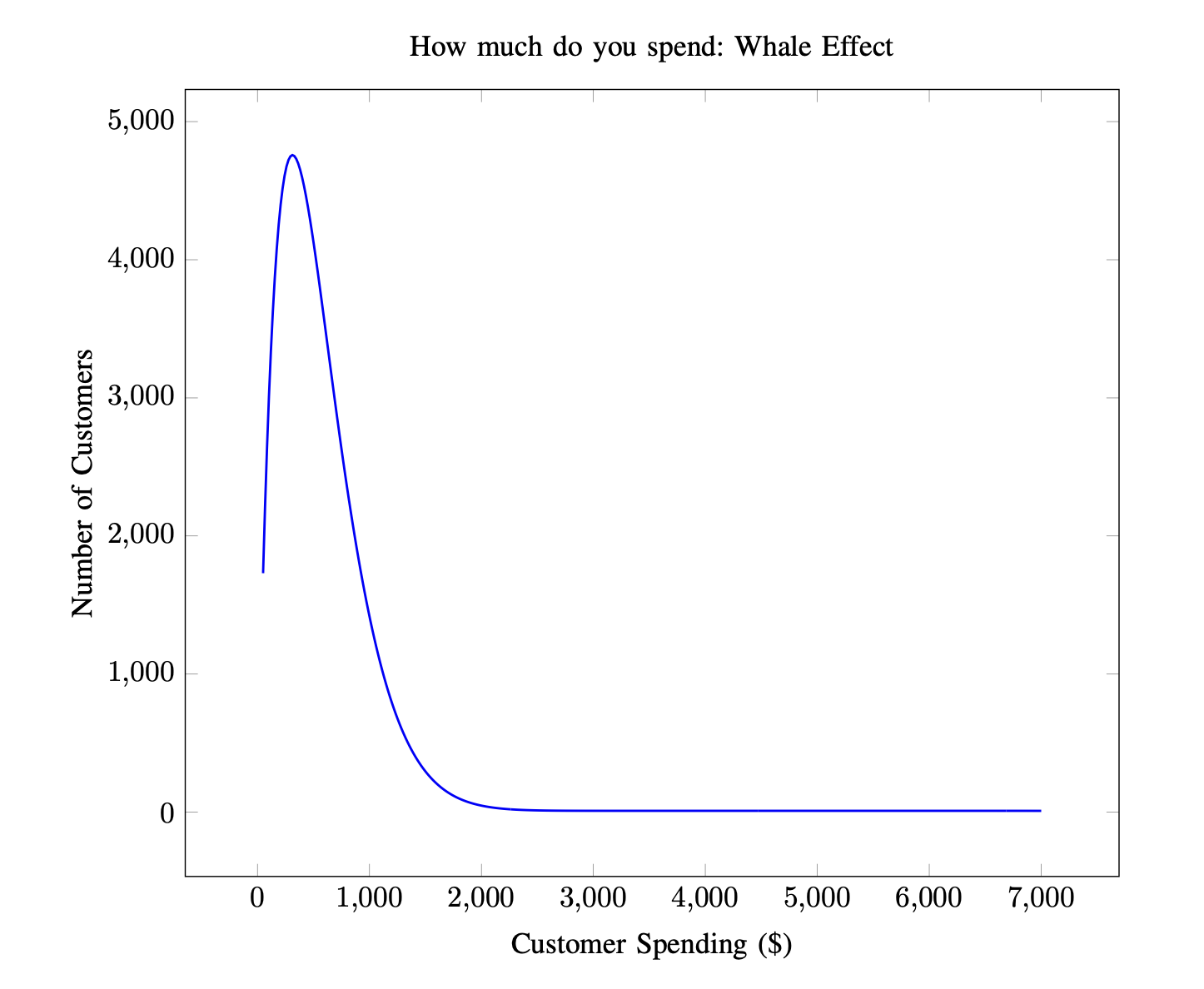}
 \caption{Top 5\% "whales" dominate revenue. Gaussian models underestimate tail risk by 60\%.}\label{demo1ros}
\end{figure}

Figure \ref{demo1ros} presents a histogram mapping annual customer spending, with the x-axis spanning \$0 to \$50,000 and the y-axis indicating customer count. The distribution is strikingly uneven: a sharp spike on the left reflects a large base of low spenders, followed by a gradual decline through the mid-range, and a long, thin tail extending far to the right, where high-value customers - those spending \$10,000 or more - reside. This fits with a Rosenblatt distribution as will see below.

Overlaying a Gaussian curve underscores the mismatch: the real-world distribution is highly skewed, lacking the symmetry assumed in standard models. This “whale effect” highlights the fundamental non-Gaussian nature of e-commerce spending, reinforcing that businesses cannot rely on simple averages. Long-term success depends on identifying and strategically catering to the outliers - the whales - who drive disproportionate value.

E-commerce businesses capitalize on the inherent asymmetry in customer spending by introducing low-budget products, transforming a vast base of minimal spenders into a meaningful revenue stream. This approach reinforces the non-Gaussian nature of spending behavior, where a small number of high-value customers drive disproportionate revenue while the majority contribute through sheer volume.

By offering low-cost products, businesses create a scalable model that thrives on numbers. While individual low-spenders generate little revenue on their own, their collective purchasing power adds up. Selling a \$5 accessory to 100,000 users, for example, yields \$500,000 - the same as a single whale spending half a million but with significantly lower risk. Lower price points also expand the customer base, drawing in price-sensitive users who might otherwise never transact. A streaming platform introducing a \$3/month ad-supported tier, for instance, converts non-paying users into revenue contributors, creating a sustainable pipeline of engaged customers.

The influx of low-cost purchases accentuates the structure of the spending distribution. More users cluster in the \$0 - \$50 range, forming a pronounced left-side spike in the histogram, while high-value spenders on the right tail remain largely unchanged. This widens the gap between the mass market and the outliers, deepening the divide that defines non-Gaussianity in e-commerce.

\begin{figure}[htb]
\centering
\includegraphics[scale=0.3]{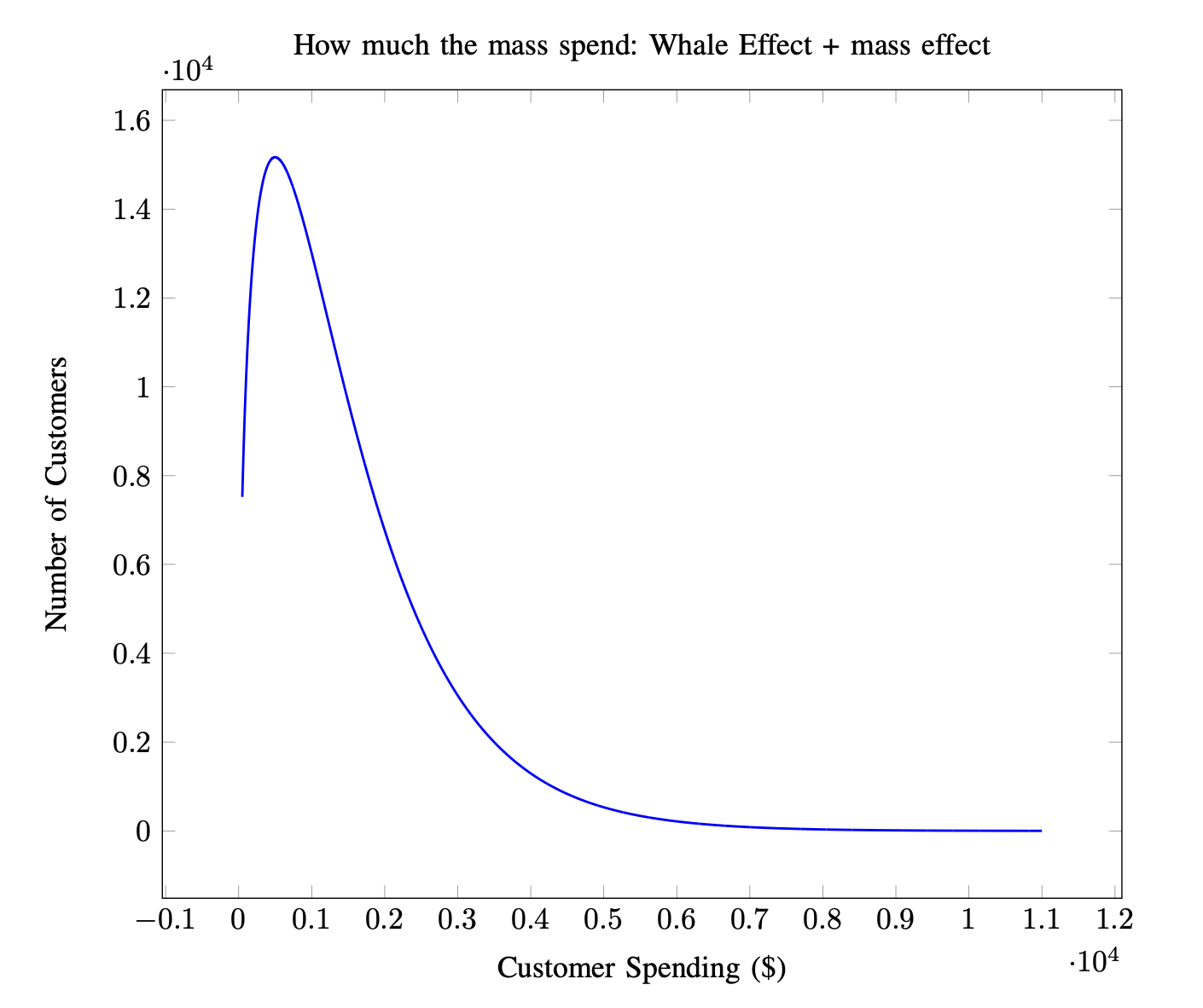}
 \caption{The influx of low-cost purchases accentuates the structure of the spending distribution making it more skewed.}
 \label{demo2ros}
\end{figure}

Before introducing low-budget products, spending distribution is already skewed: 90\% of users spend between \$0 and \$100, while only 2\% surpass the \$1,000 threshold. After introducing budget-friendly offerings, the left spike grows taller, with 95\% of users now spending between \$5 and \$100. Meanwhile, whales continue their high-spending habits, leaving the right tail intact. The result is an even more pronounced power-law, beta, chi-2, gamma, Rosenblatt distribution, where the majority shift slightly from zero-spenders to low-spenders, but the overall shape of the curve grows more extreme - further defying Gaussian assumptions. 
Introducing low-budget products does not  "normalize" spending,  it amplifies non-Gaussianity by deepening the divide between the mass of small contributors and the elite few (see Figure \ref{demo2ros}). This strategy acknowledges that e-commerce success lies not in chasing averages, but in catering to extremes: the volume of the many and the value of the few.

\subsection{Distributed Power Networks: Prosumer-Driven Net Load}

In modern distributed power networks - where consumers, prosumers (those who both consume and generate energy), and producers interact - the traditional assumption of Gaussian demand and supply patterns collapses. Instead, net load (grid demand minus localized renewable generation) follows a two-mode Rosenblatt shape, a clear departure from normality with profound implications for grid stability and energy markets.

The two-mode Rosenblatt shape, known for capturing skewness,  bimodality and fat tails, is particularly suited for energy systems where physical constraints, such as renewable generation limits,  shape outcomes. Unlike two-mode Gaussian models, which assume symmetric variability around each mode, this framework reflects real-world fluctuations driven by weather patterns, time-of-day effects, and localized energy storage dynamics.

Net load asymmetry is largely driven by renewable generation. Solar panels and wind turbines only produce energy under specific conditions, creating sharp spikes and deep troughs in supply. On sunny days, prosumers inject excess solar power into the grid, often pushing net load toward zero or negative values. At night or on overcast days, they shift to heavy grid dependence, amplifying volatility. This effect is magnified by correlated behavior: prosumers in a given region experience similar weather patterns, causing synchronized surges and drop-offs in generation.

Battery storage further complicates the landscape. Prosumers with home storage charge their batteries during surplus periods and discharge them during deficits, introducing nonlinear shifts in net load. These abrupt transitions create structural breaks in grid demand, deviating sharply from the smooth fluctuations predicted by single-mode Gaussian models.

\begin{figure}[htb]
\centering
\includegraphics[scale=0.3]{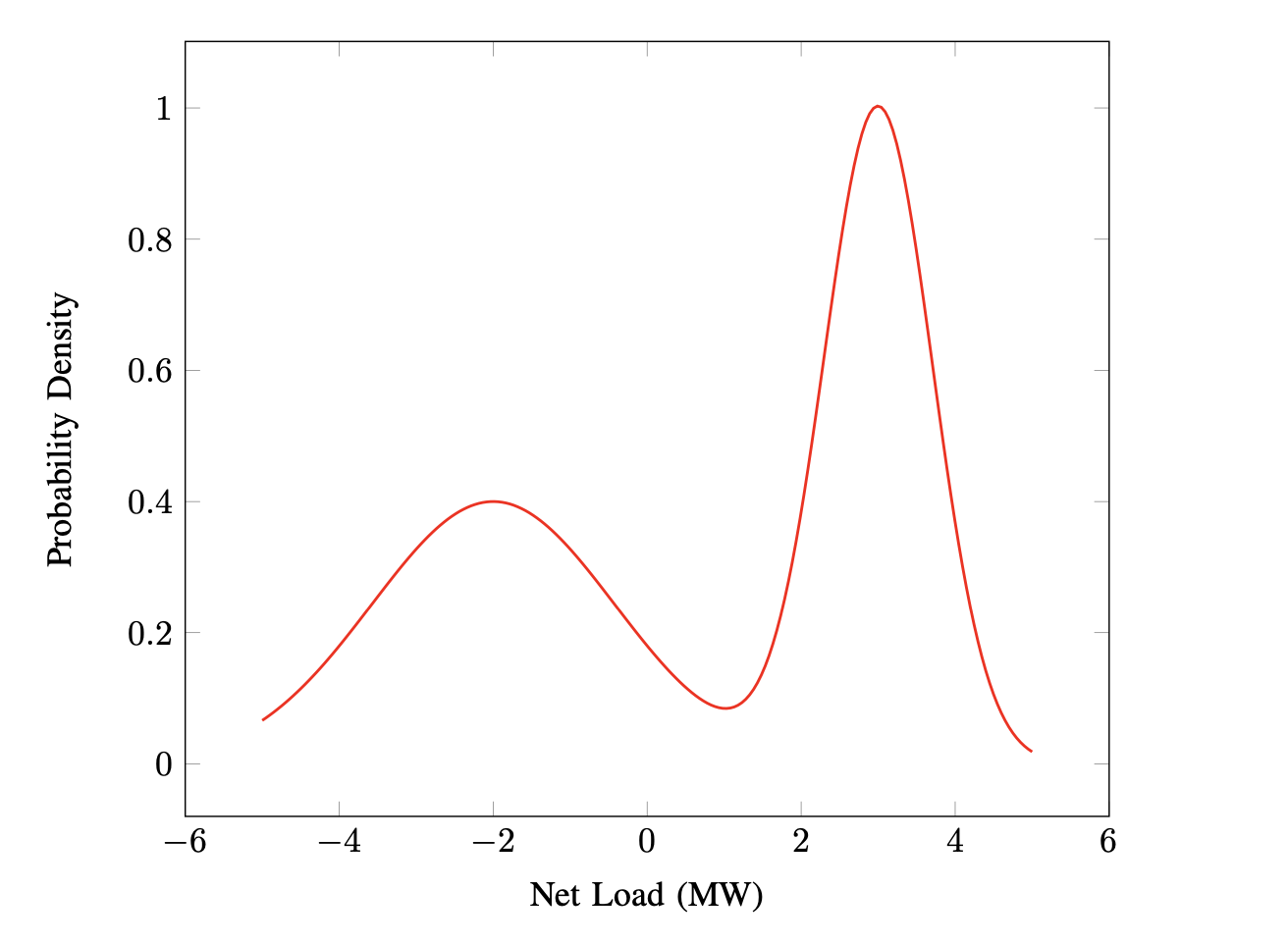}
 \caption{  Bimodal peaks (solar surplus vs evening demand) make Gaussian load forecasting error-prone. Net load= grid demand-local renewable energy generation by prosumers. } \label{demo9ros}
\end{figure}

\begin{figure}[htb]
\centering
\includegraphics[scale=0.3]{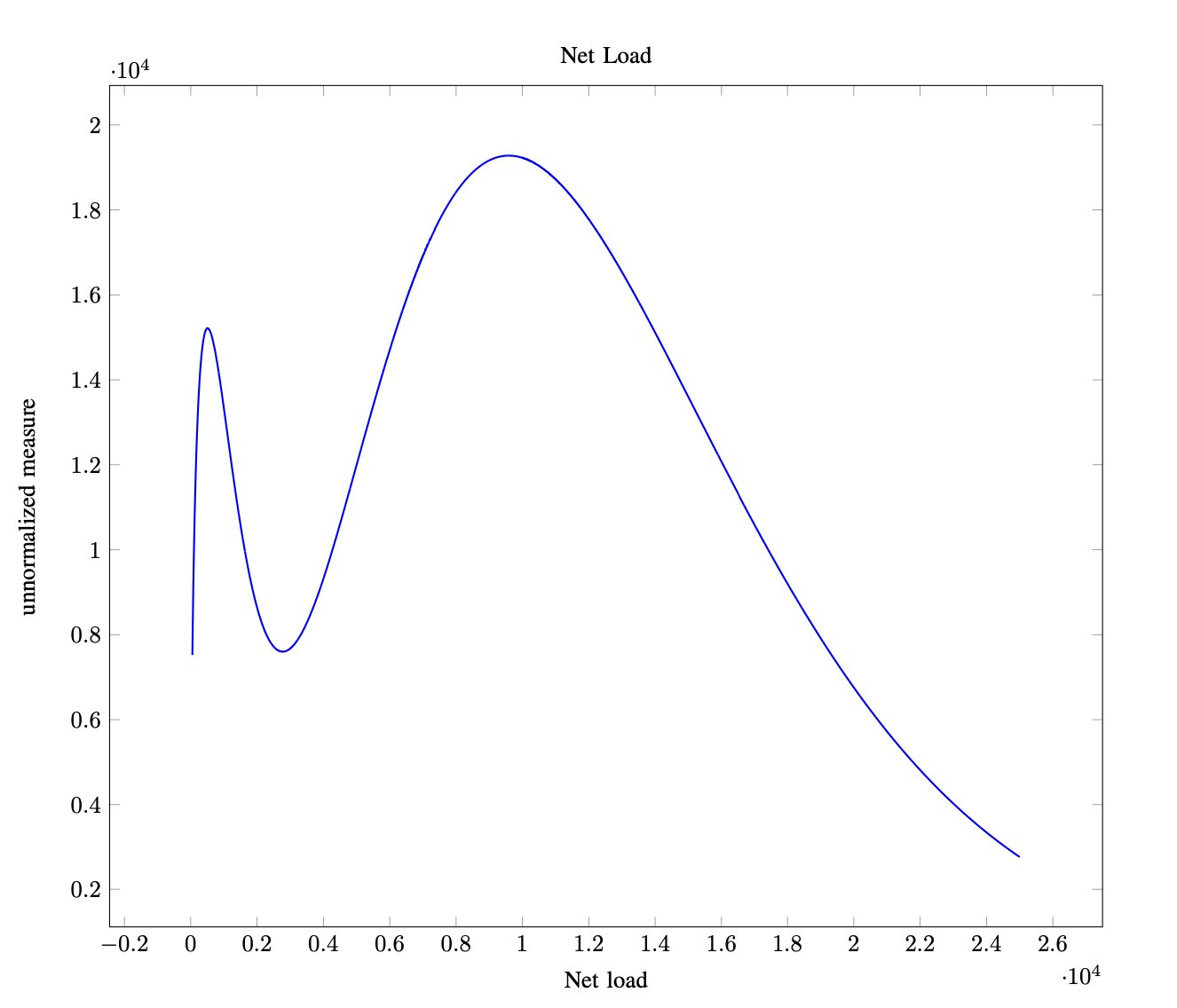}
 \caption{Net load: change in the total demand due to renewable energy sources during the sunny days and the windy days modelled using two Rosenblatt shapes.}\label{demo2ros9}
\end{figure}

In Figures \ref{demo9ros} and  \ref{demo2ros9} we consider a neighborhood with 1,000 prosumers, each equipped with rooftop solar (max output = $5 kW$). Over a single day, aggregated net load follows a two-mode Rosenblatt distribution. The
net load cannot exceed grid demand or drop below negative thresholds, where excess energy is either curtailed or sold back under feed-in tariffs.
The first peak occurs around midday, when solar generation surges and net load approaches its lower bound. The second peak emerges in the evening, when demand spikes and renewables contribute little.
Extreme events, such as grid overload during heatwaves when wind turbines underperform, occur far more frequently than a Gaussian model would predict.
Overlaying a Gaussian curve on this distribution reveals its fundamental shortcomings. The normal assumption smooths over the dual peaks, ignores the hard boundaries imposed by generation constraints, and drastically underestimates the likelihood of extreme deviations. The real-world implications are significant. Traditional grid planning models, built on Gaussian assumptions, struggle to anticipate blackouts triggered by synchronized solar drop-offs. In 2020, for instance, California's rolling blackouts were exacerbated by miscalculations in solar ramp forecasting, an oversight directly linked to an underappreciation of non-Gaussian dynamics.

Beyond grid stability, these structural distortions fuel market volatility. Electricity prices swing unpredictably as net load oscillates between extremes, sometimes dipping into negative territory when solar oversupply forces utilities to pay users to consume power. Storage systems designed for Gaussian expectations degrade prematurely, as unaccounted charge-discharge cycles stress batteries beyond intended use.

\begin{figure}[htb]
    \centering
  \includegraphics[scale=0.3]{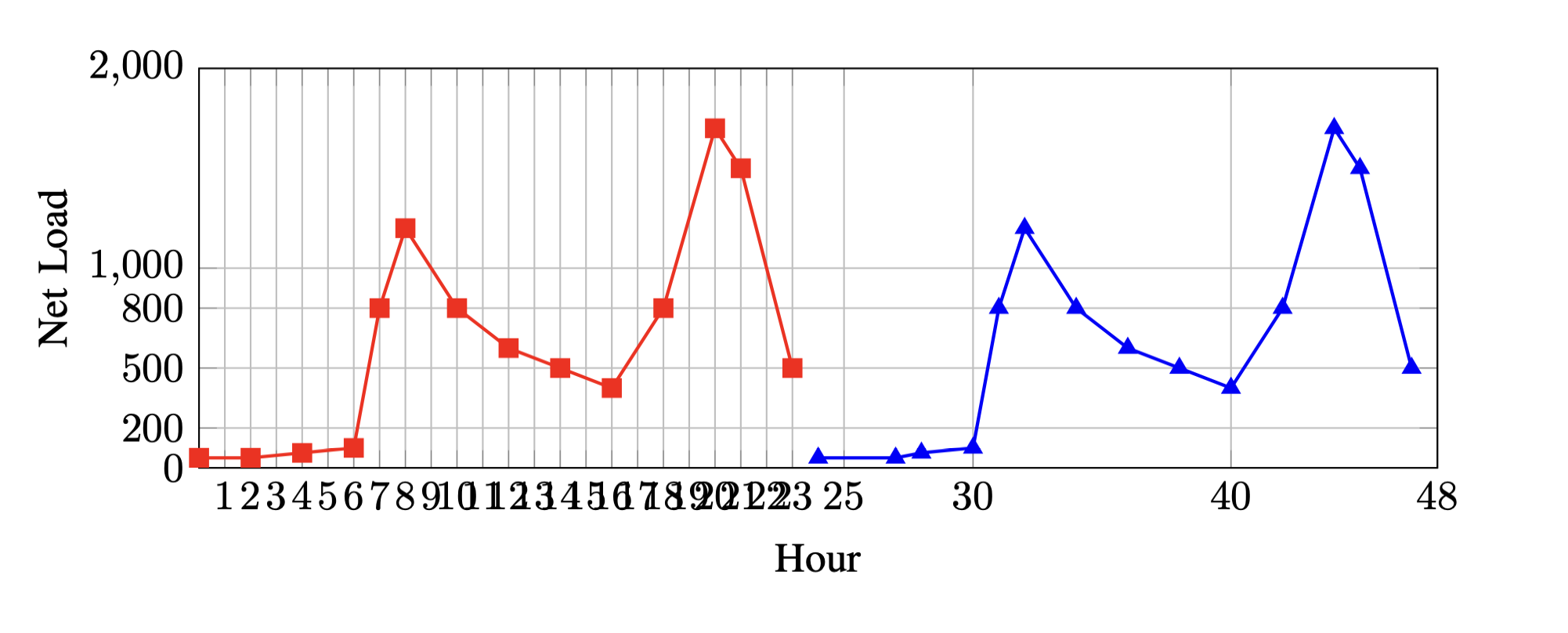}
    \caption{Daily net load during the summer}     \label{demo10ros}
\end{figure}

Empirical studies validate this shift (see \ref{demo10ros}). A 2023 and 2024 analysis of Germany's power grid found that net load in solar-dense regions followed a two-mode Rosenblatt distribution, exhibiting bimodal peaks at midday and evening. Similar observation were reported in the Sicily area in Italy where the price can be negative at certain days. In Spain, abundant solar capacity drives price pressures during sunny days when solar power generation peaks. Similarly, Finland’s substantial wind energy capacity, particularly in coastal regions, causes price dips during periods of strong winds. In both cases, high renewable energy output can saturate the grid when demand is not high enough to absorb the excess supply.  Negative electricity prices have become a regular feature of Europe's energy markets. As renewable energy capacity surges, particularly solar and wind, oversupply during low-demand periods leads to an increasing number of zero or negative price hour

Traditional Gaussian-based models mispredicted grid stress events by 40\%, while Rosenblatt-driven approaches significantly improved forecasting accuracy.

Addressing these challenges requires a fundamental rethinking of grid strategy. Non-Gaussian forecasting methods  such as, Rosenblatt or mixture models, offer superior predictive power, enabling more resilient infrastructure planning. Dynamic pricing mechanisms, which penalize or reward users based on extreme net load fluctuations, incentivize prosumers to engage in load balancing. Microgrids, designed to operate independently during demand troughs, provide localized stability, preventing wider grid disruptions.

The implications are clear. As renewable penetration accelerates, grid operators must abandon legacy assumptions and adapt to a system where extremes, not averages, dictate success or failure. The Rosenblatt distribution (or mixtures) encapsulates this reality, redefining energy networks as complex, constrained systems where ignoring statistical nuance is not just inefficient, it is catastrophic.

\subsection{Agricultural Water Usage: Punjab's Trimodal Water Crisis}
In Punjab, India - one of the world’s most water-stressed agricultural regions - groundwater usage follows a three-mode Rosenblatt distribution, not the smooth Gaussian curve assumed by traditional models. This stark deviation from normality fuels widespread economic and environmental failures, as infrastructure and policy built on flawed assumptions struggle to keep pace with reality.

Water demand in Punjab surges at three distinct points in the agricultural cycle \cite{water1}. During pre-monsoon sowing (April and May), intensive irrigation is required to prepare fields for rice, leading to the first peak in groundwater extraction. If the monsoon arrives late (June-July), farmers rely heavily on pumping to prevent seedling failure, creating a second, even larger peak. After the monsoon season (October-November), wheat planting begins, but depleted reservoirs force reliance on groundwater, driving the third peak. These three surges form a trimodal pattern, a defining feature of irrigation demand that contrasts sharply with the smooth, single-peak distribution assumed in Gaussian models.

Traditional forecasts fail to capture these abrupt shifts, leading to persistent shortages. Punjab's groundwater table declines by approximately one meter annually, driving \$1.2 billion in crop losses due to irrigation failures. With 78\% of aquifers overexploited and another 15\% at critical levels, the consequences of mischaracterizing demand are severe. Gaussian-based models predict a stable annual demand of 20 billion cubic meters (BCM), yet actual demand regularly surges to 35 BCM, straining infrastructure and leaving millions of farmers vulnerable.

\begin{figure}[htb]
\centering
 \includegraphics[scale=0.3]{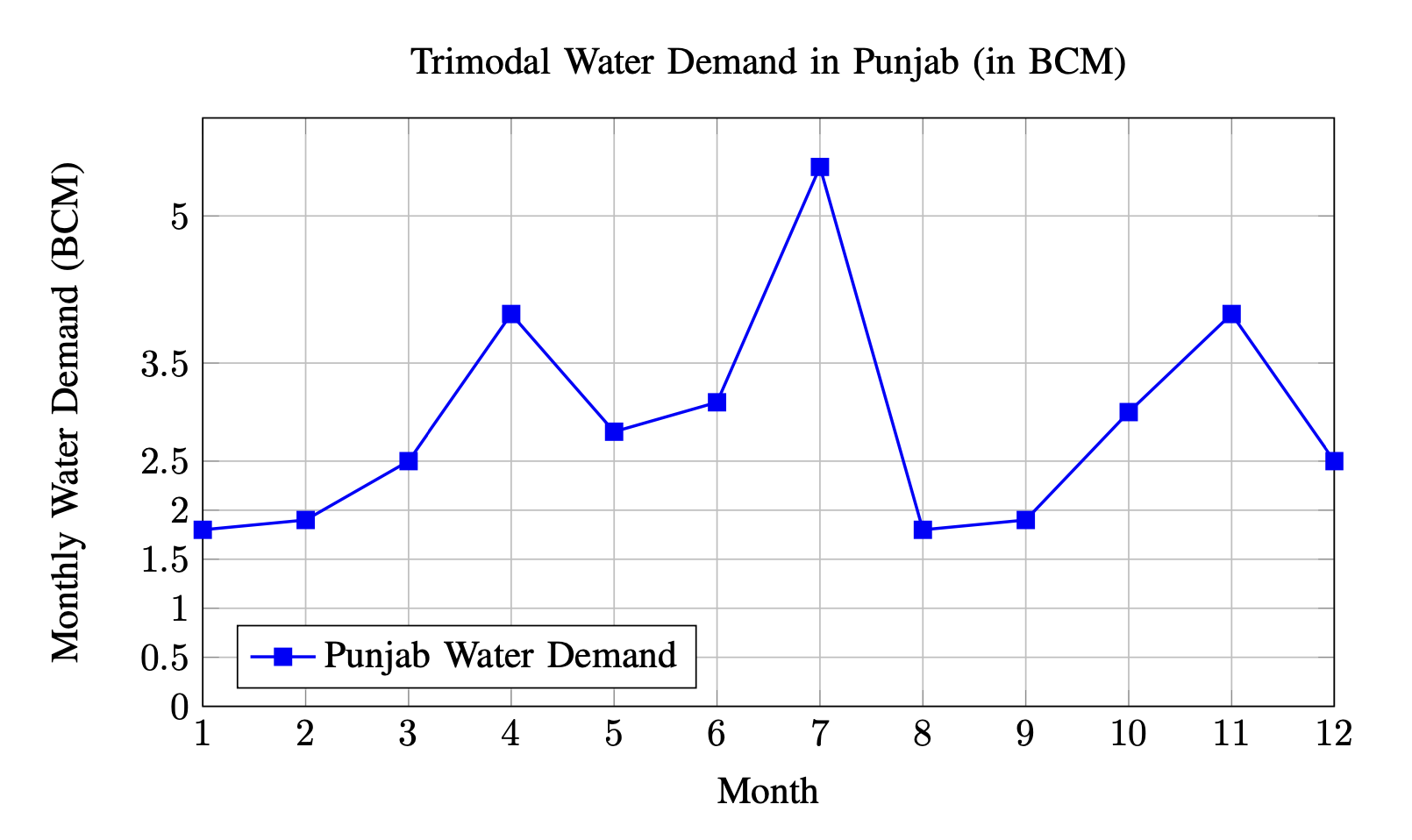}
\caption{Trimodal Water Demand in Punjab (averaged Annual Total: 30 BCM).} \label{dem11ros}
\end{figure}

A Gaussian overlay on Punjab's actual water usage reveals its fundamental miscalculation. The traditional model assumes a smooth, single-peak curve centered on monsoon months, ignoring the sharp pre-monsoon and post-monsoon spikes. The real-world three-mode Rosenblatt shape (see Figure \ref{dem11ros}), by contrast, exhibits three well-defined peaks, reflecting abrupt transitions in demand that render Gaussian-based infrastructure planning inadequate.

Ignoring this trimodal pattern has led to widespread failures. Canals and pump systems designed for "average" flows frequently collapse under peak stress, as seen in 2021 when 12,000 farmers lost entire crops after irrigation networks failed during the October surge. Agricultural policies continue to prioritize monsoon-season readiness while neglecting the far more dangerous pre- and post-monsoon crises. Supply chain disruptions compound the problem, with agrochemical and equipment suppliers routinely misallocating resources, leading to \$500 million in losses in 2022 alone.

Recognizing the trimodal non-Gaussian nature of water demand is not merely about crisis management - it opens the door to solutions. AI-driven irrigation models, designed to optimize water use based on real demand cycles, have already shown promise in reducing peak stress. Programs that recharge aquifers during monsoon surpluses provide a buffer against future shortages. Shifts toward drought-resistant crop varieties aim to reduce dependence on unsustainable extraction.

Punjab's 2022 water crisis demonstrated the high cost of Gaussian assumptions. A forecasted 22 BCM demand underestimated the real 28 BCM requirement, triggered by a six-week monsoon delay. The resulting outcome is  widespread crop abandonment, financial losses, and an urgent policy shift toward peak-sensitive interventions, including restrictions on early rice planting and incentives for alternative crops.

Punjab's water crisis is not an anomaly but a case study in the dangers of misapplying Gaussian models to  non-Gaussian systems. In resource-stressed environments, peaks, not averages, govern survival. Flattening these extremes leads to systemic failure, while embracing three-mode Rosenblatt dynamics enables more resilient infrastructure, policies, and economic strategies.

\subsection{Mango Supply Chains in Mali \& Burkina Faso}

The mango supply chains of Mali and Burkina Faso exhibit a two-mode Rosenblatt shape, with two sharply distinct seasonal peaks. However, industry planning and infrastructure remain anchored in single-peak Gaussian models, leading to systemic waste, volatile pricing, and post-harvest scarcity. Field research by the Malian startup Guinaga highlights how embracing these dynamics,  rather than forcing flawed assumptions,  can revolutionize West African agribusiness.

\begin{figure}[htb]
\centering
\includegraphics[scale=0.3]{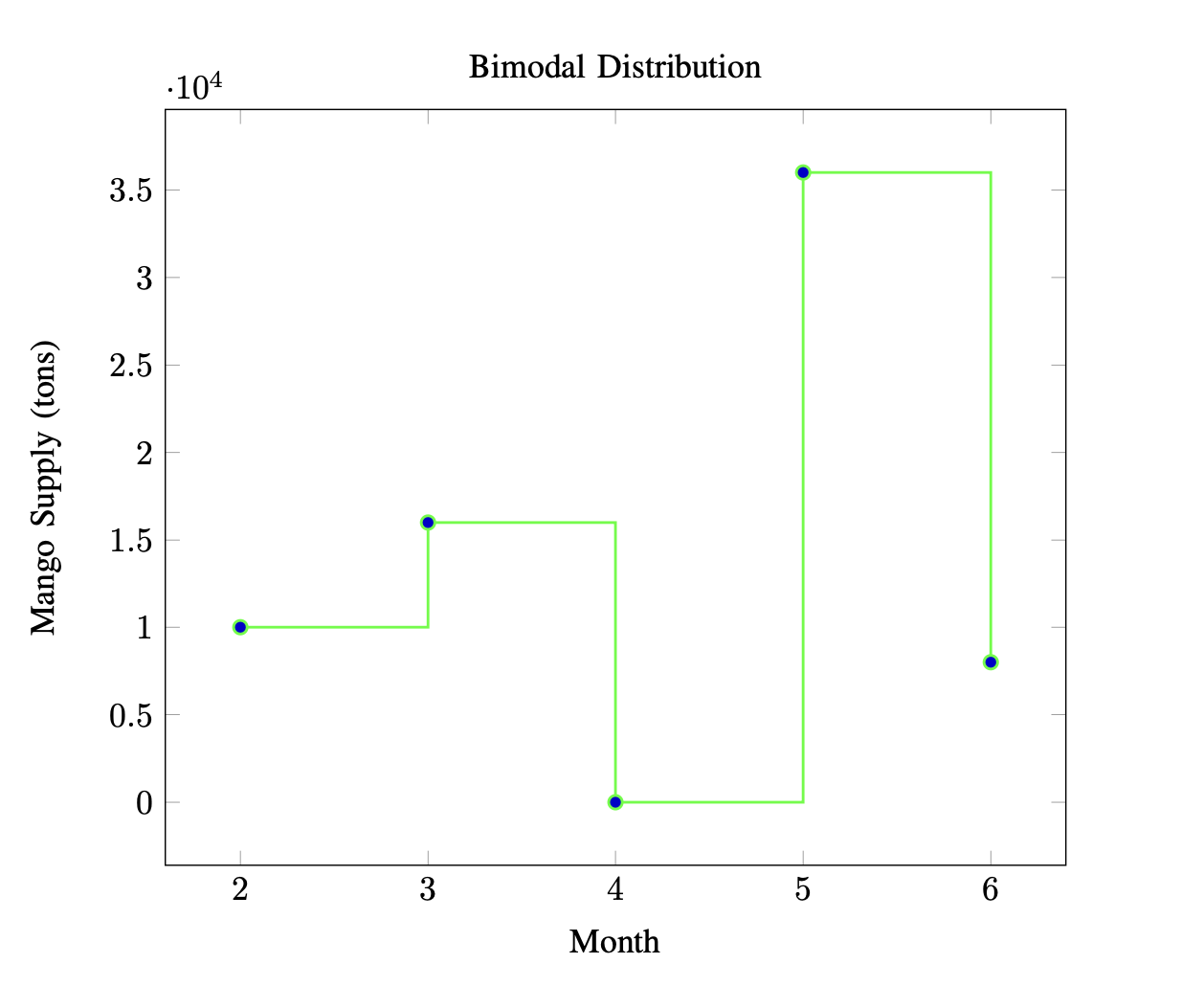}
\caption{ Bimodal Mango Supply Dual harvest seasons create 2-month scarcity gap - invisible to single-peak models. Over an average of 75 000 tons were estimated  by the Guinaga platform for the 2019-2024 period.}\label{demo3ros}

\end{figure}

In  February-March, mangoes from Siby in southwestern Mali dominate the market, supplying Bamako the capital city. This supply vanishes by late March, creating a scarcity gap before the second peak. By May and June, mangoes from Sikasso (Mali), Orodara, and Bobo-Dioulasso (Burkina Faso) flood the market, peaking in the third week of May (see Figure \ref{demo3ros}). Poor processing capacity leads to severe post-harvest loss, with thousands of tons rotting due to oversupply. After June, no fresh mangoes are available for drying, juicing, or export, despite strong market demand.

Traditional Gaussian models assume a single peak centered around the mean harvest period, smoothing out the supply curve. This misrepresentation leads to planning failures. While Gaussian assumptions predict a continuous and balanced supply, the reality follows a two-mode Rosenblatt shape with distinct peaks separated by a no-supply gap.

Post-harvest waste becomes unmanageable when market surges overwhelm storage and processing capacity. Without infrastructure to handle the supply extremes, large portions of the harvest are lost. Revenue opportunities disappear when dried mango and juice production shut down between July and January due to a lack of raw materials. Exporters and processors struggle with inconsistent supply, even as demand remains high. Price volatility intensifies as markets experience extreme fluctuations, with crashes during peak harvests and sharp increases post-June when mangoes become scarce.

The Malian agtech startup Guinaga conducted supply chain mapping across mango-growing regions, revealing critical inefficiencies in handling the two harvest peaks. Market dynamics remain unpredictable because traders and processors lack tools to navigate the disjointed seasons. Solutions include developing cold storage hubs to extend mango availability post-harvest, establishing small-scale drying and juicing units near farms to absorb surplus, and deploying digital marketplaces to connect farmers with pre-committed buyers before harvest. By incorporating two-mode Rosenblatt shape models, Guinaga helps reduce losses, stabilize supply, and increase earnings for farmers and traders.

Mali and Burkina Faso's  mango sector is a case study in how a two-mode Rosenblatt shape drives real-world agricultural dynamics. Ignoring these dual peaks leads to recurring failures - waste, price instability, and supply shortages. By shifting away from Gaussian assumptions and adopting new strategies, policymakers and agribusinesses can turn volatility into opportunity, ensuring sustainable profits for farmers while meeting year-round consumer demand.

\subsection{Dynamics of CPU, GPU, and TPU Adoption in Business (2000 - 2025)}
The evolution of computing hardware: Central Processing Unit (CPUs), Graphics Processing Units (GPUs), and Tensor Processing Units (TPUs), in the business sector demonstrates clear deviations from Gaussian assumptions across demand, price, usage efficiency, and usage over time. These trends follow multi-peaked, Rosenblatt-like, or discontinuous patterns, shaped by technological innovations, market disruptions, and shifting workloads. The following analysis showcases this non-Gaussian behavior across key metrics:

Demand trends show a multi-peaked distribution, where CPU demand plateaus post-2010, GPU demand spikes sharply between 2015 and 2020, and TPU demand accelerates post-2020. There is no singular "mean" demand period, with spikes driven by discontinuous innovations, such as the release of GPT to the public in 2022. The period from 2000-2010 saw CPUs dominating with enterprise servers and PCs, driven by linear demand growth in conjunction with Information Technology expansion. Between 2010 and 2017, GPUs surged due to machine intelligence  and machine learning adoption, with technologies like NVIDIA CUDA and deep learning. Bitcoin mining (2013-2017) caused a secondary spike in demand. From 2017 to 2025, TPUs and specialized machine intelligence  chips, like AWS Inferentia, disrupted the dominance of GPUs, especially in cloud machine intelligence environments.

Regarding price trends, CPUs experienced a decline in price per core/thread following Moore’s like Law until 2010, after which stagnation occurred due to physical limits. GPUs saw significant price surges during crypto growths (2017-2018, 2020-2021, 2024) and the machine intelligence growth cycle (2023), followed by price crashes during crypto winters. Meanwhile, TPUs saw consistent price drops as Google scaled cloud TPU offerings, benefiting from economies of scale, though their prices remained premium compared to GPUs. Price trends for GPUs exhibit a bimodal distribution, where speculative crypto and machine intelligence  bubbles create volatile price peaks, while CPUs show a negative skew, reflecting long-term price decay with occasional plateaus.

Usage efficiency metrics reveal that CPUs saw flatlining of efficiency improvements post-2010 due to diminishing returns on clock speeds. GPUs initially showed significant efficiency improvements (FLOPS/\$) until the inflationary effects of crypto mining (2017–2021). TPUs displayed exponential efficiency gains post-2018, optimized for machine intelligence  inference, although their high upfront costs limited widespread adoption. The efficiency-to-price ratio for GPUs resembles a sawtooth wave, reflecting technological advances contrasted with speculative price inflation. Meanwhile, TPUs exhibited an exponential efficiency improvement, characterized by a non-linear, power-law curve.

In terms of usage over time, CPUs were primarily used for general computing tasks like ERP and databases from 2000-2010. Between 2010 and 2020, GPUs dominated in fields like rendering, scientific computing, and machine intelligence  training. From 2020-2025, TPUs and specialized machine intelligence  chips began to take over cloud machine intelligence  inference tasks such as Google Search and GPT. These trends follow a multi-modal distribution, with each hardware type having its own distinct adoption and obsolescence cycle. Legacy CPUs still exhibit long tails of usage, continuing to serve in certain environments despite being obsolete by newer standards.

Gaussian models fail to capture the complexities of the market due to factors like demand volatility from crypto bubbles, machine intelligence  hype cycles, and supply chain shocks such as the 2021 chip shortage. Speculative markets, as seen with GPUs, and monopolistic forces, as seen with TPUs, create significant asymmetry and outliers in price distribution. Efficiency gains in hardware are often punctuated, as seen with TPU breakthroughs, or reversed, as seen with GPU price inflation during crypto booms. Additionally, legacy systems (e.g., CPUs) persist in usage far beyond their expected decay under Gaussian assumptions.

\begin{figure}[htb]
\centering
\includegraphics[scale=0.3]{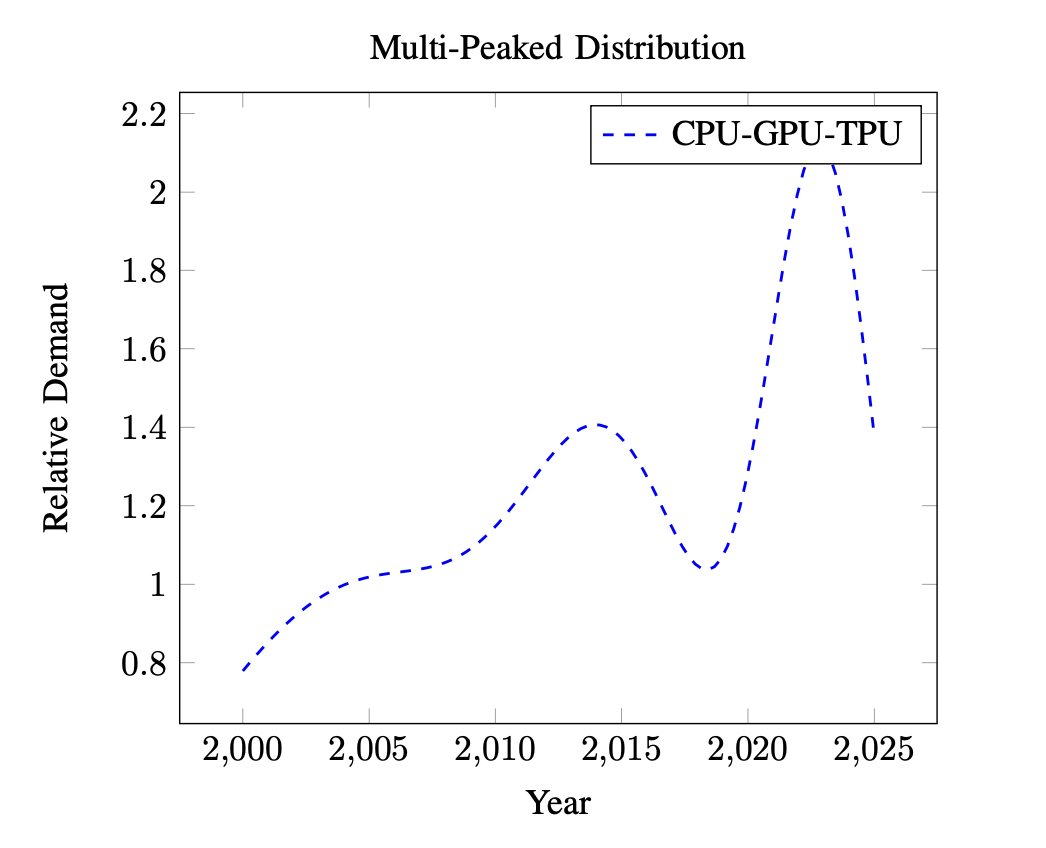}
 \caption{ Technology: CPU/GPU/TPU Adoption. Demand spikes driven by crypto (2017, 2021, 2024) and LLMs (2023) violate Gaussian stationarity.} \label{demo4ros}
\end{figure}

In Figure  \ref{demo4ros}, a real-world example of GPU demand from 2016-2023 illustrates the non-Gaussian trends. In 2016-2017, machine intelligence labs began purchasing GPUs for deep learning tasks. From 2017-2018, crypto miners hoarded GPUs, spiking prices by 300\%. Between 2020-2021, COVID-driven remote work boosted demand for consumer GPUs, and in 2022-2023, machine intelligence  startups like Stability machine intelligence  and OpenAI triggered another GPU demand rush. It was followed by llama of Meta, Grok of xAI, Gemini of Google, Qwen, DeepSeek, Kimi etc. These fluctuating events create a multi-peaked demand curve, further underscoring the non-Gaussian nature of this market.

\begin{figure}[htb]
\centering
\includegraphics[scale=0.3]{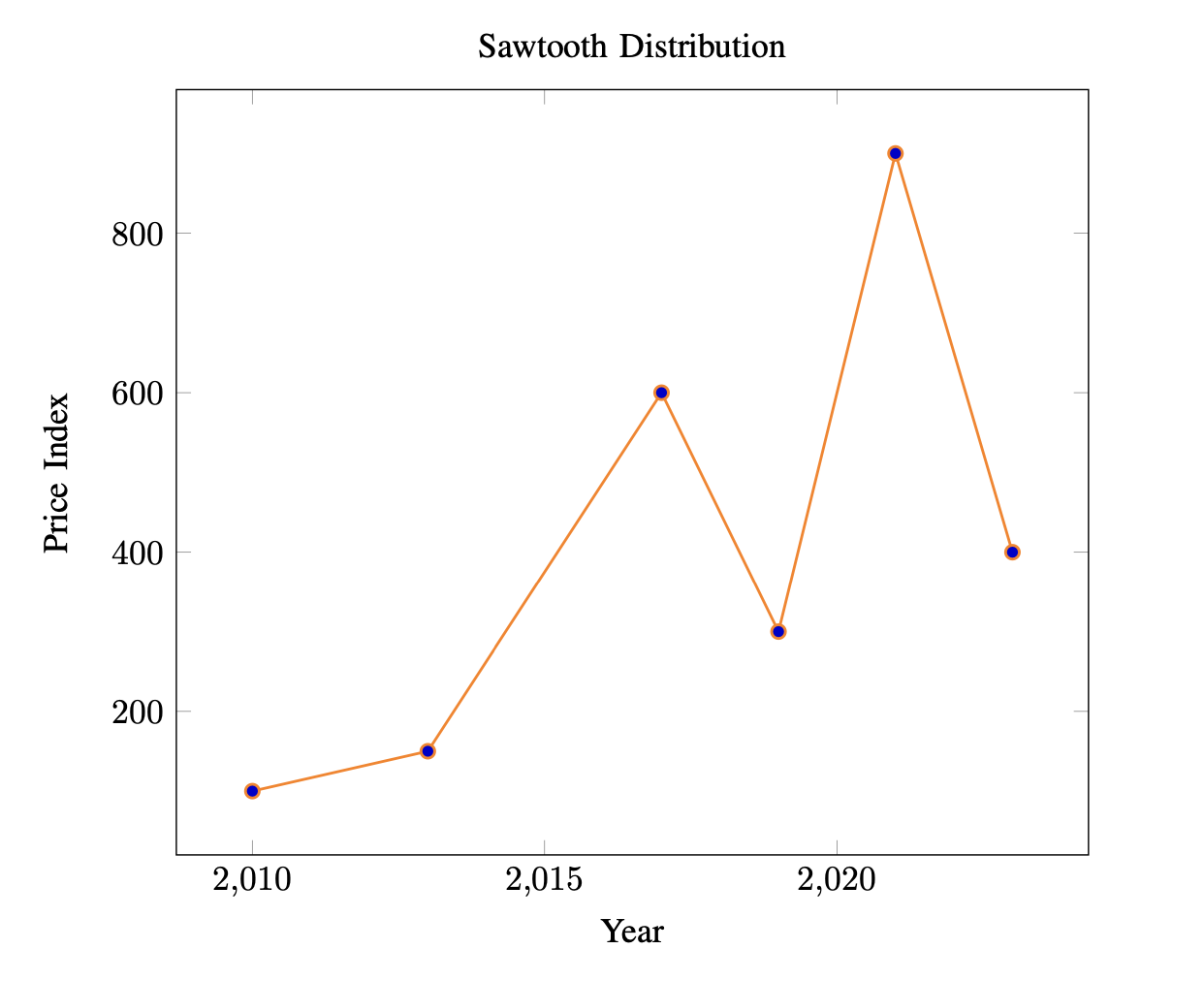}
 \caption{ Price index over the years. Crypto peaks (2017,2021, 2024) and LLM peak (2023) create non-Gaussian price shocks.} \label{demo6ros}
\end{figure}

The evolution of computing hardware, driven by crypto-induced GPU demand, AI-powered TPU adoption, and technological breakthroughs, follows non-Gaussian patterns. Businesses relying on average forecasts are unprepared for market volatility, legacy system persistence, or pricing asymmetries. Understanding power-law dynamics, multi-modal peaks, and unpredictable market events is crucial for thriving in the computing hardware sector (See price evolution in Figure \ref{demo6ros}).

\subsection{AI-Related Funding  Dynamics (1940s-2025)}

The evolution of machine intelligence  funding, spanning startups, academia, and industry, exhibits multi-peaked, Rosenblatt-like, and discontinuous trends, starkly contrasting Gaussian assumptions. Technological breakthroughs, hypergrowth cycles, and market shocks shape this non-Gaussian landscape.

Demand trends reveal a multi-peaked distribution, with peaks in the 1960s (symbolic machine intelligence), 2010s (big data), and 2020s (generative machine intelligence). Legacy academic funding persists despite the dominance of industry. From the 1940s-1950s, machine intelligence funding was minimal, mainly driven by academia and military projects like ENIAC (Electronic Numerical Integrator and Computer). The first machine intelligence growth occurred in the 1960s-1970s, fueled by DARPA funding for symbolic machine intelligence  projects like SHRDLU (an early natural-language understanding computer program that was developed by Terry Winograd at MIT in 1968-1970), with demand peaking in the mid-1960s. The machine intelligence winters of the 1980s-1990s followed, as funding collapsed due to unmet expectations, including the failure of expert systems. The 2000s-2010s saw a resurgence through machine learning, with institutions like IBM Watson and Google Brain reigniting demand. From 2012 to 2025, the deep learning explosion, driven by venture capital, saw a surge in demand, particularly post-2017 with the advent of transformer models and GPT-3.5-4-o.

In terms of price, the 1940s-1980s saw low "costs" due to limited compute power, with mainframes costing up to \$1M in the 1960s. Prices rose in the 1990s-2010s as clusters and GPUs became prevalent, with training AlexNet in 2012 costing around \$1M. By 2020-2025, exponential price growth became evident, with the cost of training GPT-4 reaching \$100M and hyperscalers like AWS and Azure dominating infrastructure spending. Price growth post-2010 follows a power-law distribution, with costs growing 10x per decade, driven by compute and data scaling laws. Outliers, such as GPT-4o-5, skew averages.

Regarding usage efficiency, the 1940s-2000s saw low efficiency, with symbolic machine intelligence  systems in the 1960s costing millions for narrow tasks. In the 2010s, efficiency spiked with GPUs (via NVIDIA CUDA), and models like ResNet delivered breakthroughs in deep learning per dollar spent. However, by the 2020s, diminishing returns were evident, with larger models like GPT-4 requiring 1000x more compute for marginal gains. Efficiency trends show an inverse U-shape, peaking between 2015 and 2020, then declining as models continued to grow. The efficiency-to-price ratio resembles a non-linear hockey stick, with a steep rise followed by a plateau.

\begin{figure}[htb]
\centering
\includegraphics[scale=0.3]{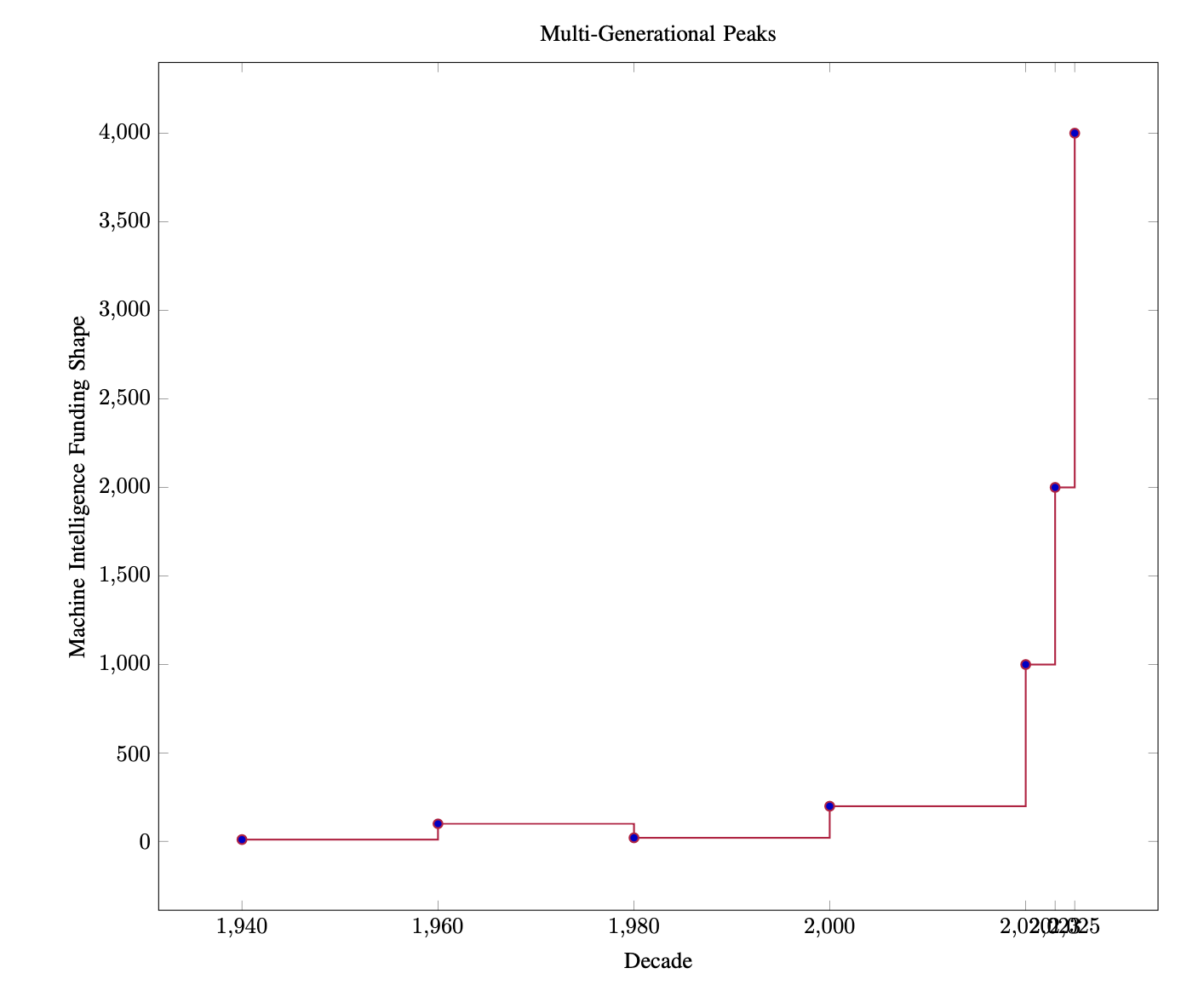}
 \caption{ Machine Intelligence Funding: Several peaks/waves/cycles:  Machine Intelligence winters (1980s) vs GPT-era hypergrowth defy smooth Gaussian projections.}
 \label{demo5ros}
\end{figure}

As illustrated in  Figure \ref{demo5ros}, usage over time follows an S-curve, with adoption initially slow before rapidly accelerating from 2010-2025. In the 1940s-1980s, machine intelligence was confined to academic labs like Stanford and MIT. By the 1990s-2010s, industry began adopting narrow machine intelligence applications (e.g., search algorithms, fraud detection). By 2020-2025, machine intelligence  became ubiquitous, with deployments in chatbots, self-driving cars, and drug discovery. Startups scaled rapidly via cloud APIs. Adoption was marked by step-function jumps, with breakthroughs like AlphaGo in 2016 and ChatGPT in 2022 driving significant adoption spikes.

Gaussian models fail to account for the multi-peaked demand, exponential price growth, diminishing efficiency, and step-function adoption of machine intelligence  technologies. Multiple disjointed peaks in demand (1960s, 2010s, 2020s) contradict a single mean/variance. Exponential cost growth and outliers like GPT-4-o-4.5 violate Gaussian symmetry, while efficiency’s rise and fall creates an asymmetric, non-bell curve. Adoption follows a rapid S-curve with disruptive steps rather than smooth Gaussian growth.

A real-world example contrasts machine intelligence  Winters and the Generative machine intelligence  growth. In the 1970s-1980s, funding dropped by 80\% after the Lighthill Report, which criticized AI's failures. From 2020-2025, generative machine intelligence  attracted over \$300B in venture capital and corporate R\&D, marking a stark contrast. This bimodal distribution (AI collapse vs. hypergrowth) cannot be modeled by a single Gaussian.
AI funding dynamics are non-Gaussian due to:
Hyper-growth/decline cycles (booms and busts),
Exponential compute costs (power-law scaling), Step-function adoption  such as GPT.
Entities like DARPA (1960s), NVIDIA (2010s), and OpenAI (2020s) thrived by adapting to these extremes. Traditional Gaussian models, assuming smooth and symmetrical trends, fail to predict events like machine intelligence  winters or hypergrowth in generative AI. The future belongs to those who model multi-modality, fat tails, and discontinuities - not averages.

\subsection{Shallot/Onion Price in Dogon Country and in Ansongo}
In a groundbreaking analysis conducted by the Malian startup Guinaga, the evolution of shallot and onion prices across two key regions: Dogon Country and Ansongo, has revealed a striking divergence from the expected Gaussian distribution commonly assumed in agricultural markets. The data, gathered through meticulous fieldwork over the past years (2019-2024)), demonstrates the skewed nature of pricing dynamics in these regions, underscoring the significant role of seasonality, supply chain bottlenecks, and external shocks.

\begin{figure}[htb]
\centering
\begin{minipage}{0.48\textwidth}
    \centering
    \textbf{(a) Monthly Shallot Prices}
\includegraphics[scale=0.2]{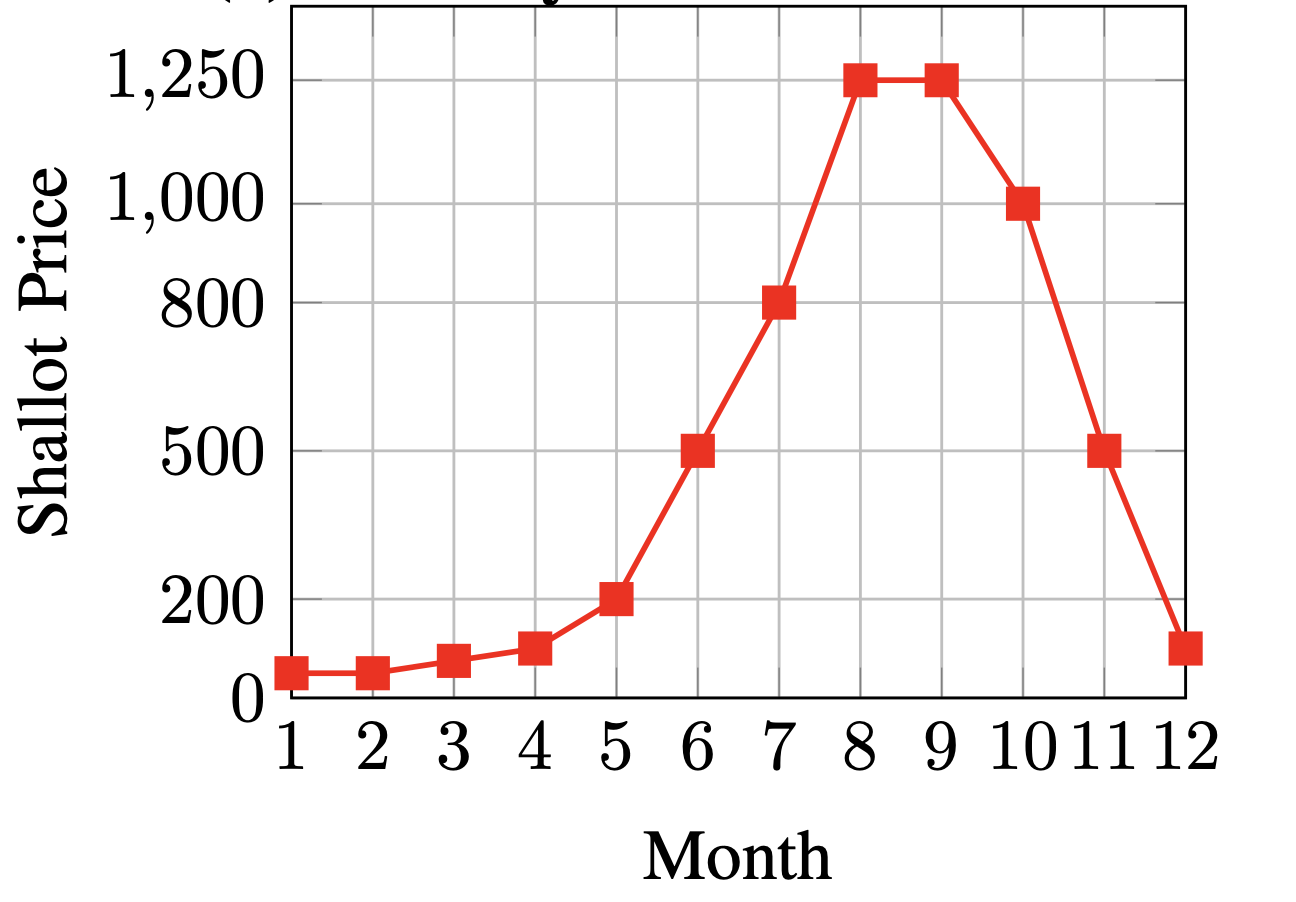}
\end{minipage}
\hfill
\begin{minipage}{0.48\textwidth}
    \centering
    \textbf{(b)  Probability Mass }
\includegraphics[scale=0.2]{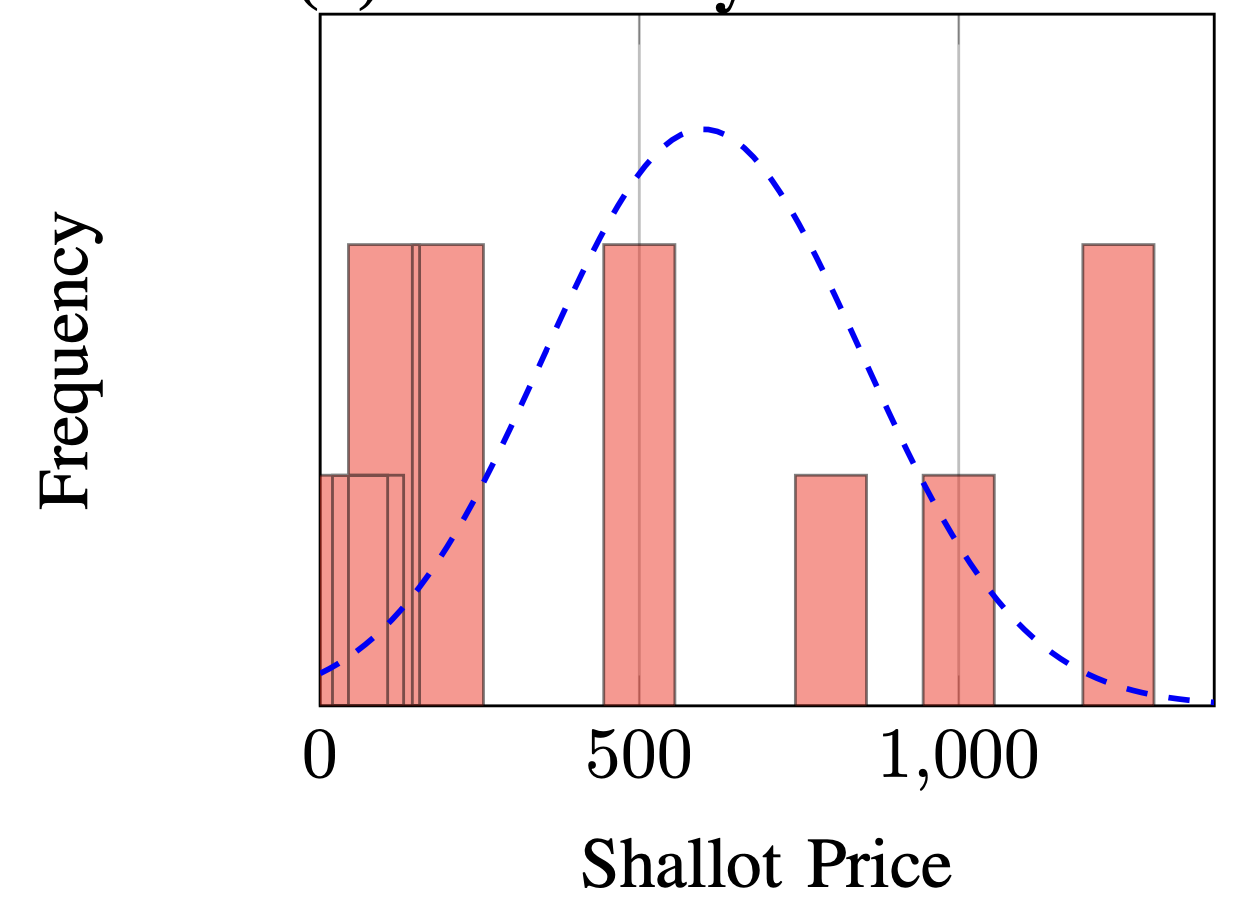}
\end{minipage}

\caption{ The Gaussian fit does not capture the skewed and heavy-tailed behavior of the actual data.} \label{demo7ros}
\end{figure}

\begin{figure}[htb]
\centering
\includegraphics[scale=0.3]{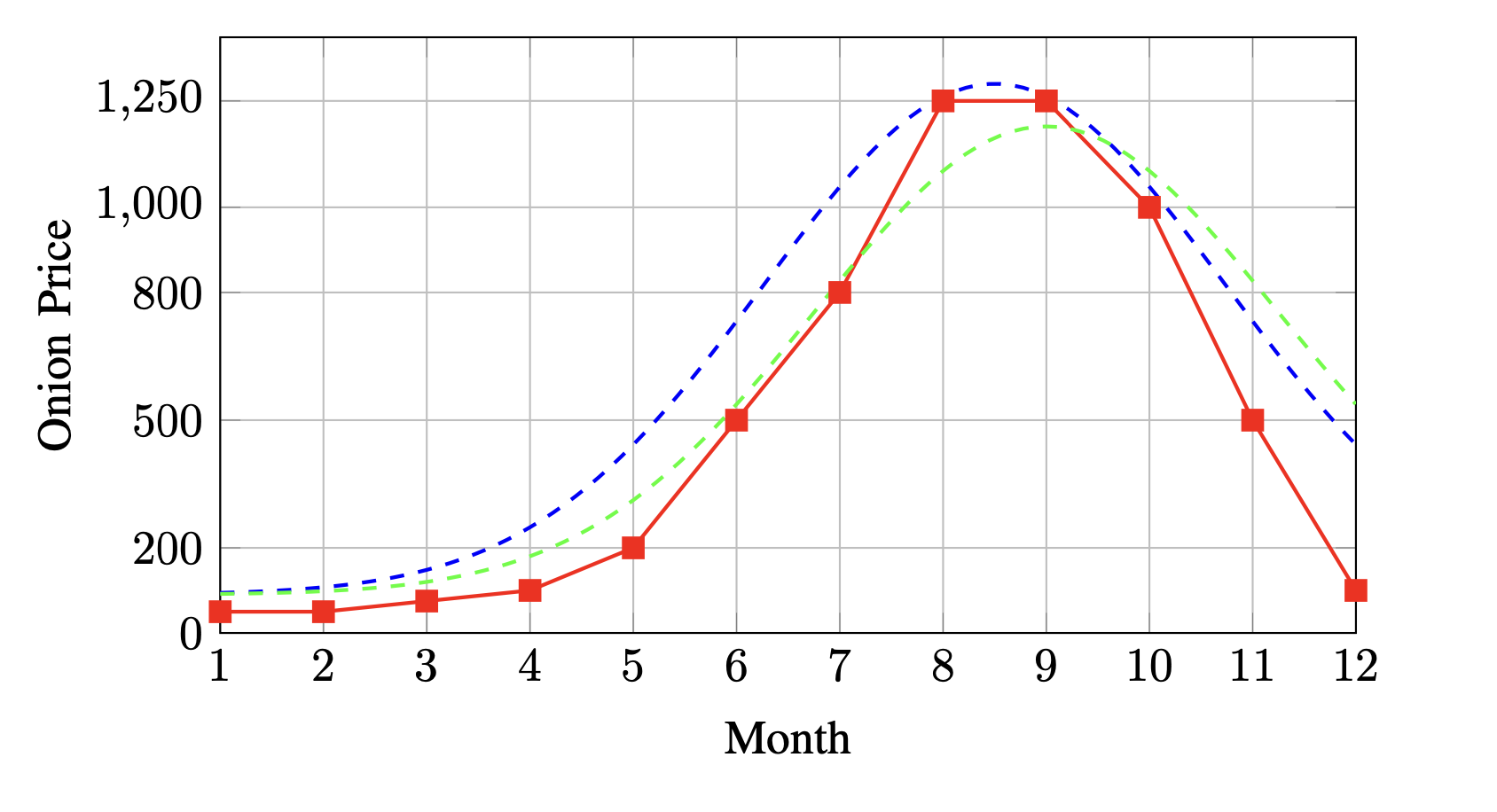}
\caption{ Non-Gaussian fluctuations in shallot and onion prices over time.} \label{demo8ros}

\end{figure}

A closer look at the shallot price data (see Figures \ref{demo7ros} and \ref{demo8ros}) for Dogon Country, particularly, illustrates the sharp volatility that can only be described as non-Gaussian as it is highly skewed. The price trajectory follows a jagged pattern, with non-symmetric shapes: one from June to November and another from December to April. This skewed shape, with abrupt increases during the raining season  and crashes during the new local shallot season, starkly contrasts the smooth, bell-shaped curve typically associated with Gaussian processes. The shallot market, much like many other agricultural goods in Mali, is governed by more than just simple supply and demand; it is shaped by unpredictable factors such as sudden climate changes, transportation delays, and local demand surges during festive seasons.

In Dogon Country, the price of shallots begins at a relatively low level (around 50 FCFA per kg) in January, steadily climbs as demand rises in the spring, and hits a staggering 1,250 FCFA per kg by August, only to plummet sharply towards the year-end. The rise and fall of prices are neither smooth nor predictable, defying standard economic models that rely on continuous, normally distributed data.
Similarly, the pricing behavior in Ansongo follows a comparable non-Gaussian trajectory, marked by extreme fluctuations that reflect the region's dependence on fluctuating water supply levels and market access, both of which are prone to abrupt changes. The typical assumption in standard market economics - that prices will gradually increase or decrease based on predictable trends - fails to account for the chaotic and unpredictable nature of agricultural production in these areas.
These findings are vital for policymakers, traders, and agricultural stakeholders who have traditionally relied on Gaussian-based models to predict and stabilize market prices. Guinaga's research advocates for a new approach - one that embraces the inherent non-linearity and unpredictability of local agricultural markets. By adopting data-driven, multi-modal models that account for these peaks and valleys, stakeholders can mitigate the risk of overproduction and avoid costly shortages.
The non-Gaussian nature of shallot prices in Dogon Country and Ansongo exemplifies the need for adaptive strategies in rural agricultural economies, emphasizing the importance of real-time data, flexible pricing mechanisms, and responsive supply chains to weather the unpredictable ebbs and flows of these crucial markets.
The  evidence from real-world data across diverse domains, from e-commerce revenue concentration to power grid fluctuations, agricultural cycles, and machine intelligence  funding dynamics, demonstrates the failure of Gaussian models in capturing underlying patterns. Instead, many of these datasets exhibit multi-mode Rosenblatt-shaped distributions, characterized by multi-peaked, heavy-tailed, and discontinuous structures that defy the smooth, symmetric assumptions of the standard Gaussian framework. Whether it is the whale effect in e-commerce, the dual-seasonal dynamics in mango supply, or the multi-modal volatility of GPU demand and pricing, the presence of distinct, dominant peaks and abrupt shifts challenges conventional statistical tools. Similarly, empirical evidence from underwater wireless communication and seasonal commodity pricing (as observed in onion markets in Dogon Country and Ansongo) further underscores how traditional Gaussian approaches fail to predict extreme variations. The trimodal water demand in Punjab and the multi-phase machine intelligence funding cycles reinforce this pattern, showing that real-world processes often result from mixtures of independent, asymmetric drivers rather than a single, smooth probability distribution. Recognizing these non-Gaussian structures is crucial for accurate forecasting, risk assessment, and decision-making across industries. 

\begin{figure}[htb]
\centering
\includegraphics[scale=0.3]{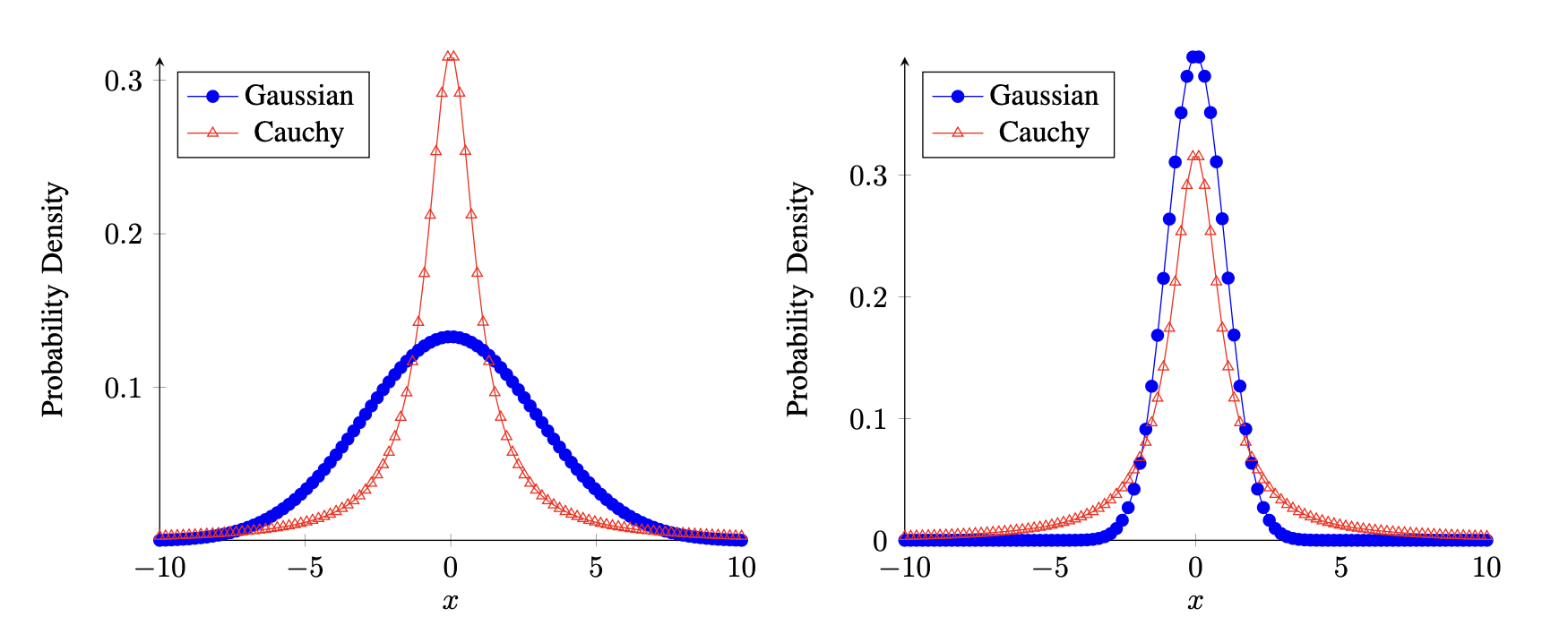} 
  \caption{Symmetry does not necessarily imply Gaussianity} \label{demo12ros}
\end{figure}

Note that the Cauchy distribution (Figure \ref{demo12ros}) belongs also to the family of continuous probability distribution, exhibits a bell shape  but differ significantly from Gaussian distribution. The Cauchy distribution  has heavy tails that decay at a much slower rate like  inverse square which means it has no finite variance or mean. This heavy-tailed behavior makes the Cauchy distribution highly sensitive to outliers, and in practice, it is often used to model phenomena where extreme values are more likely than in a Gaussian distribution.  The Cauchy distribution is used for modeling situations where extreme events dominate and the central measure (mean) is not meaningful.

\subsection{Underwater Wireless Channel State Measurements}
Long‐range underwater acoustic communication (LR-UWAC), defined as the transmission of acoustic signals over distances spanning tens to hundreds of kilometers, is pivotal for sustained maritime operations. Applications range from marine environmental monitoring and geological surveying to biological data collection and acoustic tomography, which entails the remote measurement of sound speed and water temperature. The optimization of such communication systems hinges on precise channel state measurements, achieved through techniques including channel sounding, channel state information (CSI) estimation, Doppler spectrum analysis, and extensive large-scale modeling.

In channel sounding, both time‐domain methods - employing wideband probe signals such as chirps or pseudo-random sequences and frequency-domain approaches based on orthogonal frequency-division multiplexing (OFDM) are utilized to estimate the impulse and frequency responses of the channel. These measurements yield critical parameters such as delay spread and coherence bandwidth. CSI estimation, in turn, provides a quantitative assessment of the fading conditions inherent to the underwater environment. Complementary methodologies such as Doppler analysis, which captures mobility-induced frequency shifts, and random field models that quantify fading severity, are integrated with large-scale measurement campaigns. These campaigns, leveraging drive testing, massive multiple-input multiple-output (MIMO) arrays, and machine learning-enhanced ray tracing, elucidate spatial correlation and propagation trends across diverse aquatic environments, thereby informing the design of next-generation networks (5G, 6G, and beyond).

\begin{figure}
    \centering
    \includegraphics[scale=0.3]{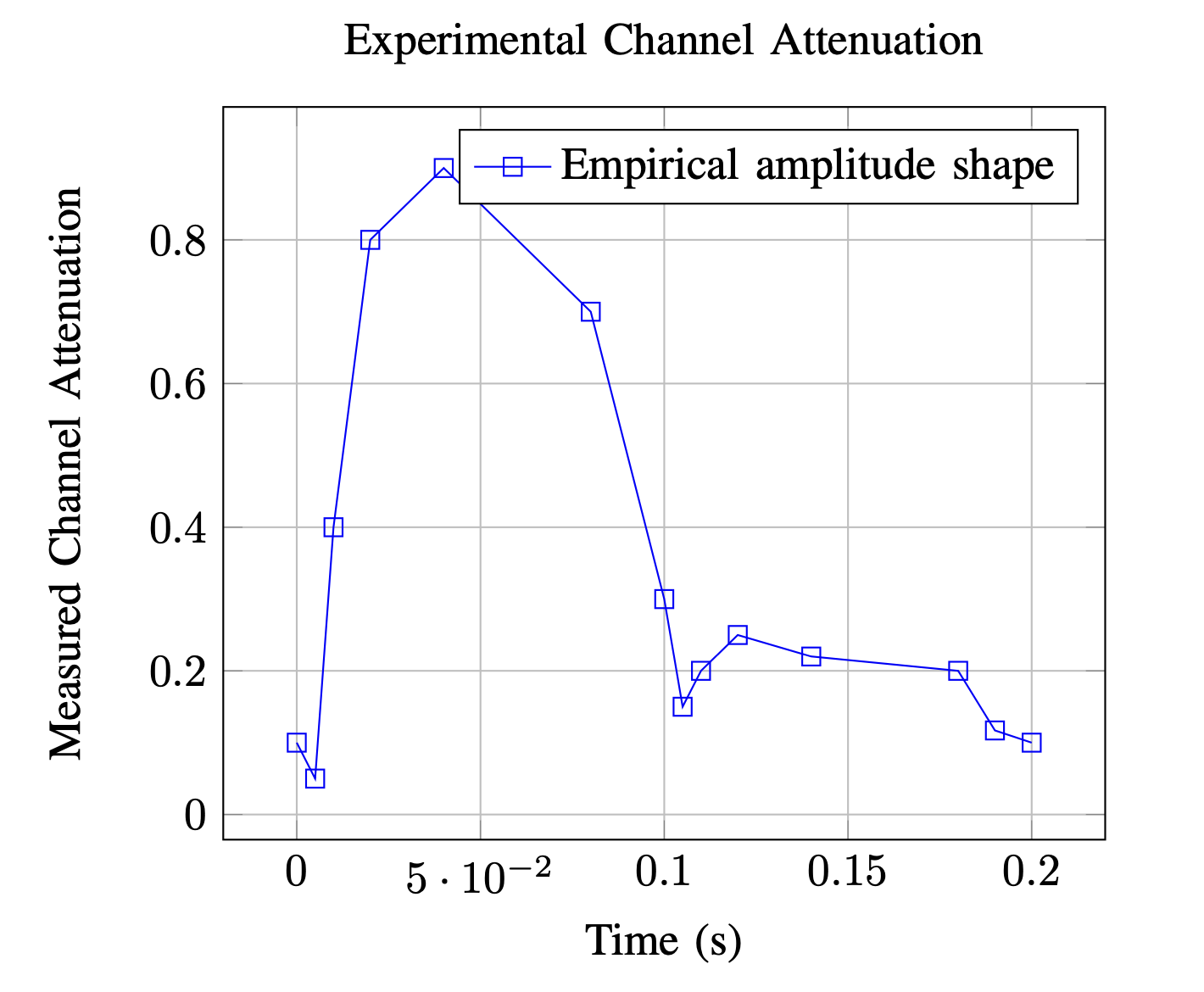} 
    \caption{Measured channel impulse response recorded under deep-water conditions at a transmitter depth of 20m.}
    \label{fig:figure_channel}
\end{figure}

Accurate quantification of fading in underwater channels necessitates an examination of signal fluctuations induced by multipath propagation, absorption, scattering, and hydrodynamic effects. Traditional models based on Rayleigh, Rician, and Nakagami-$m$ distributions, where the shape parameter $m$ serves as an index of fading severity, are frequently employed. Empirical studies involve the deployment of hydrophones and transducers to capture variations in amplitude and phase over time. Spectral broadening techniques facilitate Doppler spread analysis, revealing motion-induced distortions attributable to water currents and marine life, while coherence time estimates offer insights into channel stability. The extraction of power delay profiles further distinguishes line-of-sight from multipath components, thereby refining the probabilistic characterization of fading.

In scenarios where conventional fading models inadequately capture memory effects and non-Gaussian statistics, the non-centered Rosenblatt process provides a robust alternative. Underwater acoustic channels, influenced by multipath reflections, turbulence, and Doppler shifts from water currents, often exhibit temporal correlations that defy simple Gaussian assumptions. By estimating the Hurst parameter, which quantifies long-range dependence, researchers can adjust the Rosenblatt model to accurately represent observed power variations. Empirical fitting typically involves the collection of time-series data from hydrophone arrays, the application of fractional-order spectral analysis, and the optimization of model parameters via likelihood estimation.

Since the early 1990s, a series of LR-UWAC experiments have been conducted across multiple international locales. In 1993, adaptive multi-channel combining and equalization were demonstrated over a $203 km$ link using a 12-element vertical array, achieving error-free communications at $660 bps$ \cite{water01}. Subsequent trials, such as those reported in 2001 \cite{water2}, extended the communication range to $3250 km$ across the North Pacific Ocean from Southern California to the Hawaiian Islands, by coherently combining signals from 20 hydrophones,  at a reduced data rate of $37.5 bps.$ 

Experimental configurations typically involve a transmitter operating at variable water depths (e.g., 20, 35, and 50m) and multiple receivers deployed in both shallow and deep-water environments. One receiver, situated $100km$  northwest of the transmitter at approximately $1800 m$ depth, utilized a custom acoustic recorder based on the TASCAM DR100mk and Geospectrum hydrophones. A second receiver, deployed $50 km$  south at a comparable depth (roughly $60 m$ ), allowed for concurrent data acquisition without the need for synchronization among vessels. Bathymetric profiles and sound speed profiles (SSP) were recorded, with negligible variation in sound speed (less than $1 m/s$) observed between sites.

A variety of signal types, including chirp, JANUS, BPSK, QPSK, and OFDM, were transmitted. The chirp signal facilitated channel measurement, while the remaining signals supported data communication. Frame sizes ranged from 288 bits for JANUS to 4608 bits for QPSK, with a total of 43 frames received, 28 from the shallow site and 15 from the deep site. Variability in the received signal-to-noise ratio (SNR) was attributed to receiver drift and the time-varying nature of the channel, as exemplified by the amplitude attenuation profile in Figure \ref{fig:figure_channel}.

Timadie, a Mali-based startup focused on underwater communication, has conducted field experiments using Remotely Operated Vehicles (ROVs) to explore the dynamics of channel states in submerged environments. The results revealed a  Rosenblatt distribution, rather than the traditionally assumed Gaussian shape, in the communication channel's state over time. This non-Gaussian distribution stems from the inherent characteristics of underwater acoustics, where factors like water salinity, temperature, and turbulence create highly variable and unpredictable signal attenuation. 
The illustration was done in the Niger river in West Africa.  Operating at depths between 5 and 20 m, these ROVs communicate via acoustic modems that are particularly susceptible to multipath interference from riverbed reflections, water currents, and biological activity. Probe signals transmitted by each ROV allow for the continuous recording of received signal strength (RSSI) and signal-to-noise ratio (SNR) at 1-second intervals over a 10-minute period. The resultant power fluctuations exhibit persistent behavior, with an estimated Hurst parameter $H\approx 0.8$, indicative of long-range dependence. The empirical probability density function of the received power reveals heavy-tailed, non-Gaussian characteristics that align well with the predictions of the Rosenblatt process - yielding a Wasserstein error on the order of $10^{-6}$. In this dynamic environment, the inter-ROV distance, $d(t)$, influences the received signal power, modeled by the power-law decay: \begin{equation} P_r(t) = P_t - \eta \log_{10}\Bigl(1+d(t)\Bigr) + 2\log_{10}\Bigl(|R^H(t)|\Bigr), \end{equation} where $P_t$ denotes the transmitted power, $\eta$ the path loss exponent, and $R^H(t)$ encapsulates the long-range correlated fading component.
By integrating sophisticated measurement techniques with propagation module, this experiment predicts the behavior of underwater acoustic channels. The deployment of models such as the non-centered Rosenblatt process not only captures the inherent long-range dependence and non-Gaussian statistics of underwater fading but also informs the design of robust communication systems capable of operating in marine environments.

\subsection{State of Health of batteries of Electric Cabs} \label{secsoh}
Accurately predicting electric vehicle battery health and remaining useful life (RUL)  is a predictive maintenance challenge.  The State of Health (SOH) assesses the battery's degradation over time, reflecting its ability to store and deliver energy compared to its pristine state, also typically represented as a percentage. 
  SOH is affected by factors such as age, temperature, charging patterns, and usage, and it directly influences the battery's effective range and remaining lifespan. Both State of Charges  and SOH are critical parameters monitored by the Battery Management System (BMS) for accurate range estimation and battery health assessment.
{T}he performance of batteries in electric vehicles is affected by both time-dependent degradation (calendar aging) and usage-dependent degradation (cycle aging). In regions with significant seasonal temperature variations, these effects differ  between seasons.  The RUL is the remaining cycles until capacity faded to $\alpha=80\%$ of nominal (end-of-life threshold). A combined analysis that captures these differences over  full 5 years, for predicting the RUL of batteries in electric vehicles operating in regions characterized by four distinct seasons. Battery capacity degradation results from both calendar aging, which follows an Arrhenius behavior, and cycle aging, which is dependent on usage intensity.   Based on  battery degradation data analysis with voltage (V), current (A), surface temperature (\textdegree C, via thermocouples), capacity (Ah), and internal resistance ($\Omega$) - sampled at 1 Hz using programmable testbeds (Arbin BT2000), alongside derived metrics such as SOC and state-of-health (SOH); cycling protocols included CC-CV charging (1C rate, CC-CV charging is a typical method of charging rechargeable batteries such as li-ion, the Operation switches between CC charging, which charges with a constant current, and CV that charges at a constant voltage, depending on the voltage of the rechargeable battery), dynamic discharge profiles (UDDS: Urban Dynamometer Driving Schedule), and calendar aging tests (20-80\% SOC, 25-45\textdegree C) to emulate diverse operational stresses, 
battery chemistries and characteristics impact the RUL. EV battery types include Lithium Nickel Manganese Cobalt Oxide (NMC) batteries, Lithium Nickel Cobalt Aluminum Oxide (NCA) batteries, Lithium Iron Phosphate (LFP) batteries, Lithium Manganese Oxide (LMO) batteries, Lithium Titanate (LTO) batteries, Lithium Polymer (LiPo) batteries, Solid-State Lithium-Ion batteries, Nickel-Metal Hydride (NiMH) batteries, Lead-Acid batteries, and Lithium Sulfur (Li-S) batteries.

\begin{figure}[htb] 
\centering
\includegraphics[scale=0.3]{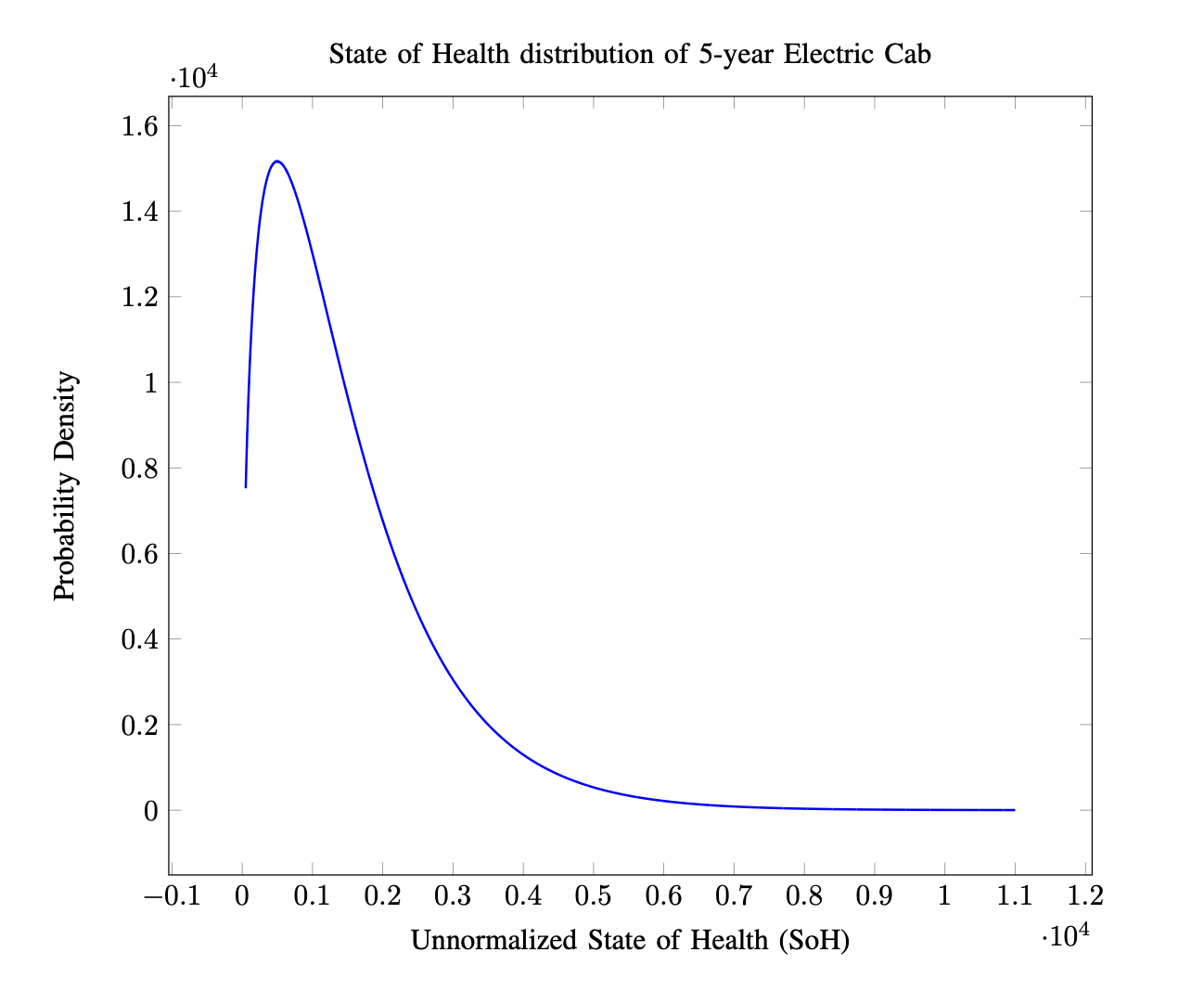} 
\caption{Probability density function of SoH for EV batteries in 5+ year-old electric cabs, exhibiting a Rosenblatt-like heavy-tailed structure.}\label{figsoh}
\end{figure}

The long-term degradation of electric vehicle (EV) batteries in high-utilization fleet applications, such as electric cabs, exhibits a stochastic profile that is well captured by a shape of a Rosenblatt distribution. The analysis of SoH in cabs older than five years reveals a departure from conventional Gaussian, instead following a heavy-tailed structure with a varying Hurst parameter  indicating different degrees of long-memory dependence in the degradation process. The presence of a wide spectrum of $H$ values, spanning from 0.5 to 0.9, highlights the self-similar and fractal nature of battery wear, where certain conditions - such as operational stress, non-uniform charge cycles, and environmental variability - drive long-range correlations in SoH decline. This  challenges traditional SoH prognostics and supports the need for fractional-order models to accurately predict remaining battery life in aged EV fleets. The probability density function of SoH, as illustrated in Figure \ref{figsoh}, underscores this non-Gaussian behavior, demonstrating asymmetry and heavy tails that signify the persistence of degradation modes. The Rosenblatt framework enables more precise risk-aware energy management strategies for second-life applications. The wide range of Hurst parameters observed in the data suggests that the SoH of batteries in electric vehicles is affected by a variety of factors, including the type of battery, the driving habits of the driver, and the environment in which the vehicle is driven.

\section{What is a Rosenblatt process ? }

The Rosenblatt process stands out as a emerging sophisticated tool for modeling real-world phenomena that defy traditional Gaussian, Markovian assumptions. Unlike the familiar bell curve of Gaussian distributions, the Rosenblatt process is a non-Gaussian, self-similar process characterized by long-range dependence and higher-order dependencies. It is defined using a double Wiener-Itô integral, which incorporates cumulants beyond the second order, providing a more nuanced representation of stochastic behavior. 

\begin{definition}
 The Rosenblatt process \( R^H(t) \) is defined using a double Wiener-Itô integral:

$$\quad t \in \mathbb{R}, 
\ \ R^H(t)= C^H_R \int_{\mathbb{R}^2} \left( \int_0^t (u-y_1)_+^{\frac{H}{2}-1} (u-y_2)_+^{\frac{H}{2}-1} du \right) dW{y_1}dW{y_2}, 
$$
with Beta: $
B(x,y) = \int_0^1 t^{x-1}(1-t)^{y-1} dt, 
$
               $C_R^H := \sqrt{\frac{2H(2H-1)}{2B(1-H, \frac{H}{2})}} $ ensuring that $\mathbb{E}[(R^H(t))^2]=1.$
                $\frac{1}{2} <H <1.$
 \end{definition}

 It exhibits self-similarity and non-Gaussianity. The Rosenblatt process satisfies the self-similarity property:  
              \[
              R^H(c t) \overset{d}{=} c^H R^H(t), \quad \forall c > 0.
              \]
              Unlike fractional Brownian motion (FBM), it is not Gaussian, as it arises from a nonlinear transformation of FBM. It is characterized by strong long-range dependence. \[
              \mathbb{E}[R^H(t) R^H(s)] \sim  \frac{1}{2}\left( t^{2H} + s^{2H} - |t-s|^{2H} \right),
              \]
              
             Rosenblatt process is not a semimartingale. It requires a new stochastic calculus (that will be presented in the next section) that can be used in optimal stochastic control theory and stochastic differential game theory.

\section{Systems driven by Rosenblatt noise  }

Traditional stochastic calculus, primarily developed for Gaussian processes, and recent stochastic calculus developed for Poisson jump processes, fractional Brownian motion and Gauss-Volterra process falls short in capturing the behavior of systems influenced by non-Markovian, non-Gaussian noise, such as the Rosenblatt process. The new stochastic calculus for Rosenblatt processes extends classical techniques by incorporating these non-Gaussian characteristics, enabling more accurate modeling of systems with heavy-tailed distributions and non-central limit  fluctuations. This advancement is crucial for developing control strategies in various applications, including machine intelligence, power grid management, and agricultural supply chains, where traditional Gaussian assumptions often lead to suboptimal performance. We use the following notations: 
\begin{equation}\begin{array}{ll}
c= C^H_R\Gamma^2(\frac{H}{2}), \  \Gamma(z)= \int_0^\infty t^{z-1} e^{-t}\,dt\\
\tilde{c}= \sqrt{ \frac{(2H-1)\Gamma(1-\frac{H}{2})\Gamma(\frac{H}{2})}{(H+1)\Gamma(1-H)}},\\
\nabla^{\frac{H}{2}} =I^{\frac{H}{2}}_{+}\circ D,\   \ 
\nabla^{\frac{H}{2},\frac{H}{2}} =I^{\frac{H}{2},\frac{H}{2}}_{+,+}\circ D^2,\\
f \ : \mathbb{R} \rightarrow \mathbb{R},\ \ \\
(I^{\alpha}_{+}f)(x):=\frac{1}{\Gamma(\alpha)}\int_{-\infty}^{x} f(u) (x-u)^{\alpha-1}du,\\
(\alpha_1,\alpha_2)\in (0,1)^2,\  f \ : \mathbb{R}^2 \rightarrow \mathbb{R},\\
(I^{\alpha_1,\alpha_2}_{+,+}f)(x_1,x_2)
:= \frac{1}{\Gamma(\alpha_1)\Gamma(\alpha_2)}\int_{-\infty}^{x_1}\int_{-\infty}^{x_2} f(u,v) (x_1-u)^{\alpha_1-1}(x_2-v)^{\alpha_2-1} dudv,
\end{array}
\end{equation}

 \begin{thm}  \label{refresult0} 
We make the following assumption of the function $f:$
       \begin{itemize}

            \item  $H \in (\frac{1}{2}, 1),  T>0,  $
             \item  $f\in \mathcal{C}^2([0,T]\times  \mathbb{R}),$
             
              \item $f(t,.)\in \mathcal{C}^3( \mathbb{R}),$
                 \item $|f_{xxx}(t,x)| \leq c_t (1+|x|^\alpha), x\in  \mathbb{R}$
               \item $d_1, d_2,d_3$ safisfy integrability in  $L^1, L^{\frac{2}{1+H}}, L^{\frac{1}{H}}$ and have continuous representation 
               \cite{examplebook}.
        \end{itemize}
        
 Let $x$ be a state dynamics driven by a Rosenblatt noise given by 
        $$x(t)= x_0 +\int_0^t d_1 \ dt'+ 2\tilde{c}\int_0^t d_2 dB^{\frac{H+1}{2}}(t')+  \int_0^t  d_3 {\color{blue} dR^H(t')}\ $$. Then, the process $y$ defined by $y(t) = f(t,x(t))$ can be written as

    $$ y(t)= y_0 +\int_0^t \tilde{d}_1 \ dt'+ 2\tilde{c}\int_0^t \tilde{d}_2 dB^{\frac{H+1}{2}}(t')+  \int_0^t  \tilde{d}_3  {\color{blue}dR^H(t')}$$
 where 
 \begin{equation}\label{refresult1} \begin{array}{ll}
\tilde{d}_1(t')=f_{t'}(t',x_{t'}) +f_{x}(t',x_{t'}) d_1 + 2 c f_{xx}(t',x_{t'})(\nabla^{\frac{H}{2}}x_{t'})(t')d_2 
\\ + c f_{xx}(t',x_{t'})(\nabla^{\frac{H}{2},\frac{H}{2}}x_{t'})(t',t')d_3
+c f_{xxx}(t',x_{t'})[(\nabla^{\frac{H}{2}}x_{t'})(t')]^2d_3,\\  \\
\tilde{d}_2(t')= f_x(t',x_{t'}) d_2+ f_{xx}(t',x_{t'})(\nabla^{\frac{H}{2}} x_{t'})(t')d_3,\\ \\
\tilde{d}_3(t')=f_x (t',x_{t'}) d_3,
\end{array}
\end{equation}
\end{thm} 

This is a new Stochastic Calculus established in \cite{examplebook}.
It is useful to note that a third derivative occurs in the drift
term and that a stochastic integral with respect to a  fractional
Brownian motion occurs.

We will use this  result to:
\begin{itemize} \item   design  control strategies for systems with Rosenblatt noise, enabling better management of systems with long-range dependence and non-Gaussian behavior. 
\item  determine the best control actions to achieve desired outcomes in systems driven by Rosenblatt processes, leading to improved efficiency and performance.
\item develop control systems that can adapt to changing environments and uncertainties, enhancing robustness and stability in the presence of Rosenblatt noise. 
\item  analyze strategic interactions between agents in dynamic games with Rosenblatt noise, providing insights into competitive, cooperative or coopetitive  behavior under long-range dependence and non-Gaussianity. 
\item  solve  systems with interacting agents  of mean-field type influenced by Rosenblatt processes
\end{itemize}

Detailed analysis of Rosenblatt processes can be found in \cite{reftcou1,reftcou2,reftcou3,reftcou4,reftcou5,reftcou6,reftcou7,reftcou8,reftcou9,reftcou10,reftcou11,reftcou12,reftcou13}

            \section{  Prediction under  Rosenblatt noise } 
     
     Traditional prediction techniques often rely on the standard conditional expectation, which assumes Gaussian and Markovian properties and linear dependencies. However, these assumptions do not hold in the context of Rosenblatt noise, which is characterized by long-range dependence and higher-order cumulants. The standard conditional expectation fails to capture the  dependencies and memory effects inherent in Rosenblatt processes, leading to inaccurate predictions and suboptimal control strategies.   Following  \cite{examplebook3}, we use a direct method    to  establish a prediction formula for the mean
squared error, specifically tailored for linear state dynamics  under Rosenblatt noise. 
     \begin{thm}
   \begin{equation}\nonumber \begin{array}{ll}
  x(t) =x_0+\int_0^t  b_1(t')x(t')dt' +{\color{blue} R^H(t)},
                  \end{array}
\end{equation}

 The optimal mean
square error prediction for the solution of the linear equation is
given explicitly:
\begin{equation} 0<t'<t: \ \ 
 \mathbb{E} [x(t) \mid x(t''); 0 < t'' < t'] = x(t')
\end{equation}
  \end{thm}

For a scalar constant coefficient linear differential equation given by
\[
dX(t) = b_1X(t) dt + dR^H(t),
\]
the solution is
\[
X_t = X_0e^{b_1 t} + \int_0^t e^{b_1(t - t')} dR^H(t'), \quad t > 0.
\]

\begin{thm}
Fix $t > 0$. The optimal linear prediction of the random variable  $X_t$ given the observations $(X(u), u \in (-s, 0))$ as follows:
\[
\hat{X}(t)= e^{b_1t} X_0 + e^{b_1t} \int_{-s}^0 \hat{g}_{\exp}(r) d(e^{-b_1\cdot} X_{\cdot})(t'),
\]
where
\[
\hat{g}_{\exp}(x) := e^{b_1x} G_{\exp} \left( -\frac{x}{s} \right), \quad x \in (-s, 0),
\]
\[
G_{\exp}(x) = -c(H)^{-1} x^{\frac{1}{2} - H} \frac{d}{dx} \int_x^1 \left[ \xi^{2H - 1} (\xi - x)^{\frac{1}{2} - H} \frac{d}{d\xi} \int_0^\xi \eta^{\frac{1}{2} - H} (\xi - \eta)^{\frac{1}{2} - H} F_{\exp}(\eta) d\eta \right] d\xi,
\]
$c(H)=2 \cos \left( \frac{1}{2} \pi (1 - 2H) \right) \Gamma(2H - 1) \Gamma \left( \frac{3}{2} - H \right)^2$  and 
\[
F_{\exp}(y) = s^{1 - 2H} \int_0^t e^{-b_1x} (x + sy)^{2H - 2} dx, \quad y \in (0, 1).
\]

\end{thm}
This result is extended in \cite{examplebook3duncan} for process in the integral form. Leet $h$ be a positive continuous function on the real line. The optimal linear predictor of $\int_0^t h(t'') dR^H(t'')$ based on the past $\left( \int_u^v h dR^H, -t' < u < v \leq 0 \right)$ can be explicitly obtained using the orthogonality property
\[
E \left[ \int_0^t h(r) dR^H(r) - \int_{-s}^0 \hat{g}(x) d \left( \int_{-s}^\cdot h dR^H \right) (x) \right] \left[ \int_u^v h(r) dR^H(r) \right] = 0.
\]

\section{Control of Linear Systems  driven by Rosenblatt noise  }
Optimal control problems for systems driven by Rosenblatt noise is a new  area of research with significant implications for both theoretical advancements and practical applications. Traditional control theory often relies on Gaussian, Poissonian and Markovian noise assumptions, which fail to capture the  many shapes observed in real-world data. Empirical evidence from various domains, such as e-commerce, power grid operations, agricultural supply chains, and hardware demand cycles, consistently demonstrates non-Gaussian characteristics, including skewed, heavy tails, multimodal behaviors, and strong long-range dependencies. The Rosenblatt process can capture some of  these higher-order dependencies and memory effects. By integrating Rosenblatt noise into control theory, 
we aim to obtain more accurate models for uncertainty propagation and feedback control, which are crucial for systems characterized by extreme events and non-central limit  fluctuations. This approach not only enhances predictive accuracy but also improves control performance in systems where traditional Gaussian models fall short. 

 We start with the infinite horizon stochastic optimal control with long-time average cost and linear state dynamics  driven by Rosenblatt noise process.          
\begin{equation}\nonumber \begin{array}{ll}
L_\infty(u) = \limsup_{T \rightarrow \infty} \frac{1}{T} \mathbb{E} \int_0^T (q x^2(t) + r u^2(t)) dt, \\
  x(t) = x_0 + \int_0^t (b_1 x(t') + b_2 u(t')) dt' + {\color{blue} R^H(t)}, \\
  (b_1,b_2,x_0)\in  \mathbb{R}^3, b_2\neq 0, \ q>0, r>0,  H \in (\frac{1}{2}, 1),   \\
                  \mathcal{U} := \{(u(t))_{t \geq 0}: u(t) = K x(t) \text{ with } K \in \mathbb{R}\}
                  \end{array}
\end{equation}
                  
                 The control of the ergodic cost under Rosenblatt is given by   $\inf_{u\in \mathcal{U}}L_\infty(u)$ subject to the state dynamics driven by Rosenblatt noise\cite{examplebook2}.

          We use a direct method to find linear  state- feedback optimal control.
          
            \begin{thm}  \label{refresult2} 
The optimal gain  
  in the family of admissible feedbacks is given by

$$
\hat{K} = -\frac{b_1 + \sqrt{b_1^2 + 4H(1-H) \frac{b^2_2q}{r}}}{2b_2(1-H)}
$$

and the optimal cost is given by

$
L_\infty (\hat{K}) = \frac{\Gamma(2H)}{[-(b_1 + b_2\hat{K})]^{2H-1}} \left(-\frac{r\hat{K}}{b_2}\right)=
 \frac{\Gamma(2H+1)}{2[-(b_1 + b_2\hat{K})]^{2H}} \left(q+r\hat{K}^2\right).
$ 
  \end{thm}

Theorem  \ref{refresult2}  establishes the optimal gain for a linear-quadratic control problem driven by a Rosenblatt process by an explicit representation from a variation-of-constants formula. It is combined with a stochastic calculus formula for the Rosenblatt process and is  applied to the squared state variable. By taking expectations and analyzing the asymptotic behavior of the integral terms, the proof derives an integral equation for the expected value of the squared state variable in the long-run. The optimal gain is found by minimizing the cost functional, leading to a quadratic equation whose solution provides the optimal gain and cost. 

  \begin{proof}

The state equation has a solution that can be explicitly written using a variation-of-constants formula or integration by parts, that is,
\[
X_t = e^{(b_1+b_2K)t} x_0 + \int_0^t e^{(b_1+b_2K)(t-s)} dR^H(s).
\]
Using this explicit form of this solution, we set for $y_t = X_t$ and $f(t,x) = x^2$ and 
we use Theorem \ref{refresult0}  and Equation (\ref{refresult0}) to obtain
\[
X_T^2 = x_0^2 + \int_0^T \left[ 2(b_1 + b_2K)X_s^2 + 2c \nabla^{\frac{H}{2},\frac{H}{2} }  X_s(s,s) \right] ds + 2\tilde{c} \int_0^T 2\nabla^{\frac{H}{2} }  X_s(s) dB^{\frac{H}{2}+\frac{1}{2}}(s) + \int_0^T 2X_s dR^H(s),
\]
where $T > 0$. By taking the expectation of both sides of the above equation, it follows that
\[
\mathbb{E}X_T^2 = x_0^2 + 2(b_1 + b_2K) \mathbb{E} \int_0^T X_s^2 ds + 2c \mathbb{E} \int_0^T \nabla^{\frac{H}{2},\frac{H}{2} }  X_s(s,s) ds,
\]
because the stochastic integrals have zero expectation. The second-order fractional stochastic derivative can be computed as follows. Denote
\[
o_t(t') \triangleq \nabla^{\frac{H}{2},\frac{H}{2} } X_t(t',t').
\]
Applying the operator $\nabla^{\frac{H}{2},\frac{H}{2} }$ to equation   $x(t) = x_0 + \int_0^t (b_1 x(t') + b_2 u(t')) dt' + {\color{blue} R^H(t)},$ , it follows that $o_t(u)$ satisfies the non-homogeneous linear differential equation
\[
o_t(t') = (b_1 + b_2K) \int_0^t o_s(t') ds + \nabla^{\frac{H}{2},\frac{H}{2} } R_t(t',t')
\]
for almost every $t' \in \mathbb{R}$. We use
\[
\nabla^{\frac{H}{2},\frac{H}{2} } R^H_t(t',t') = \tilde{C}_H \int_0^t |t'-r|^{2H-2} dr,
\]
where the constant is given by
\[
\tilde{C}_H \triangleq 2c_R^H \frac{B^2\left(\frac{H}{2}, 1-\frac{H}{2}\right)}{ \Gamma^2(\frac{H}{2})}
\]
By solving the integral equation, it follows that
\[
o_t(t') = \tilde{C}_H \int_0^t e^{(b_1+b_2K)(t-r)} |t'-r|^{2H-2} dr,
\]
which by a change of variables yields
\[
o_s(s) = \nabla^{\frac{H}{2},\frac{H}{2} }X_s(s,s) = \tilde{C}_H \int_0^s e^{(b_1+b_2K)r} r^{2H-2} dr.
\]
By substituting this last expression into the equation, simplifying the constants, and dividing both sides by $T$, it follows that
\[
\frac{1}{T} \mathbb{E}X_T^2 = \frac{1}{T} x_0^2 + 2(b_1 + b_2K) \frac{1}{T} \mathbb{E} \int_0^T X_s^2 ds + 2H(2H-1) \frac{1}{T} \int_0^T \int_0^s e^{(b_1+b_2K)r} r^{2H-2} dr ds.
\]
The convergence of the separate terms in the above equation is now studied. Without loss of generality, it can be assumed that $b_1 + b_2K < 0$ because otherwise, the cost $L_\infty(K)$ would be infinite. Consequently, by the definition of the Gamma function there is the convergence
\[
\int_0^s e^{(b_1+b_2K)r} r^{2H-2} dr \to  \frac{\Gamma(2H-1) }{[-(b_1 + b_2K)]^{2H-1}}
\]
which assures that also the average
\[
\frac{1}{T} \int_0^T \int_0^s e^{(b_1+b_2K)r} r^{2H-2} dr ds
\]
converges to the same constant as $T \to \infty$. This convergence can also be used to show that
\[
\lim_{T \to \infty} \frac{1}{T} \mathbb{E}X_T^2 = 0
\]
because $\mathbb{E}X_T^2$ can be computed explicitly using formula. Hence, by letting $T \to \infty$ in the equation and using the definition of the cost functional, it follows that
\[
0 = \frac{2(b_1 + b_2K) }{q + rK^2} L_\infty (K) + \frac{\Gamma(2H + 1)}{[-(b_1 + b_2K)]^{2H-1}}
\]
from which $L_\infty(K)$ can be isolated:
\[
L_\infty(K) = \frac{\Gamma(2H + 1)(q + rK^2)}{2[-(b_1 + b_2K)]^{2H}}.
\]
Minimizing the above formula, it follows that the optimal gain $\hat{K}$ must satisfy the quadratic equation
\[
(1-H)b \hat{K}^2 + b_1 \hat{K} - H \frac{b_2q}{r} = 0
\]
from which the formula follows. Moreover, note that
\[
\frac{H(q + r \hat{K}^2)}{- (b_1 + b_2 \hat{K})} = -\frac{r}{b_2} \hat{K}
\]
which, together with the equation, yields the optimal cost. This completes the proof.   \qed

  \end{proof}
  
  We now link the optimal gain to the positive solution of an algebraic Riccati equation, connecting the control problem to classical control theory results.

  It follows from the  quadratic equation
$
(1-H)b_2 \hat{K}^2 + b_1 \hat{K} - H \frac{b_2q}{r} = 0$
  that if $P$ is the positive solution to the following algebraic Riccati equation
\[
(1-H) \frac{b_2^2}{ r} P^2 + b_1P - Hq = 0,
\]
then the optimal gain can be written in the form:
\[
\hat{K} = \frac{b_2}{r} P.
\]

  \subsection{What if we replace the Rosenblatt noise by a standard one  and Optimize ?   }       
  
\begin{equation}\nonumber 
\begin{array}{ll}
\inf_{u\in \mathcal{U} }  \limsup_{T \rightarrow \infty} \frac{1}{T} \mathbb{E} \int_0^T (qx^2(t) + ru^2(t)) dt, \\
x(t) = x_0 + \int_0^t (b_1x(t') + b_2u(t')) dt' +  \textcolor{red}{ \cancel{R^H(t)} }\\
 \qquad  \qquad  \qquad  \qquad \qquad  \qquad\qquad \qquad \approx \textcolor{blue}{B(t)} \\
  \qquad  \qquad  \qquad  \qquad \qquad  \qquad\qquad\qquad \approx \textcolor{blue}{B^H(t)} \\
  \qquad  \qquad  \qquad  \qquad \qquad  \qquad\qquad\qquad \approx \textcolor{blue}{\tilde N(t,.)} \\
   \qquad  \qquad  \qquad  \qquad \qquad  \qquad\qquad \qquad\approx \textcolor{blue}{0}. \\
\end{array}
\end{equation}
This operation leads to suboptimality.
  Approximating the Rosenblatt noise by another one such as Brownian motion, fractional Brownian, Gauss-Volterra  noise, compensated Poisson jump process leads to \textcolor{red}{suboptimal control}.

Failure of Noise Approximation-Optimization Cascade
 
\begin{figure}
\centering

 \includegraphics[scale=0.3]{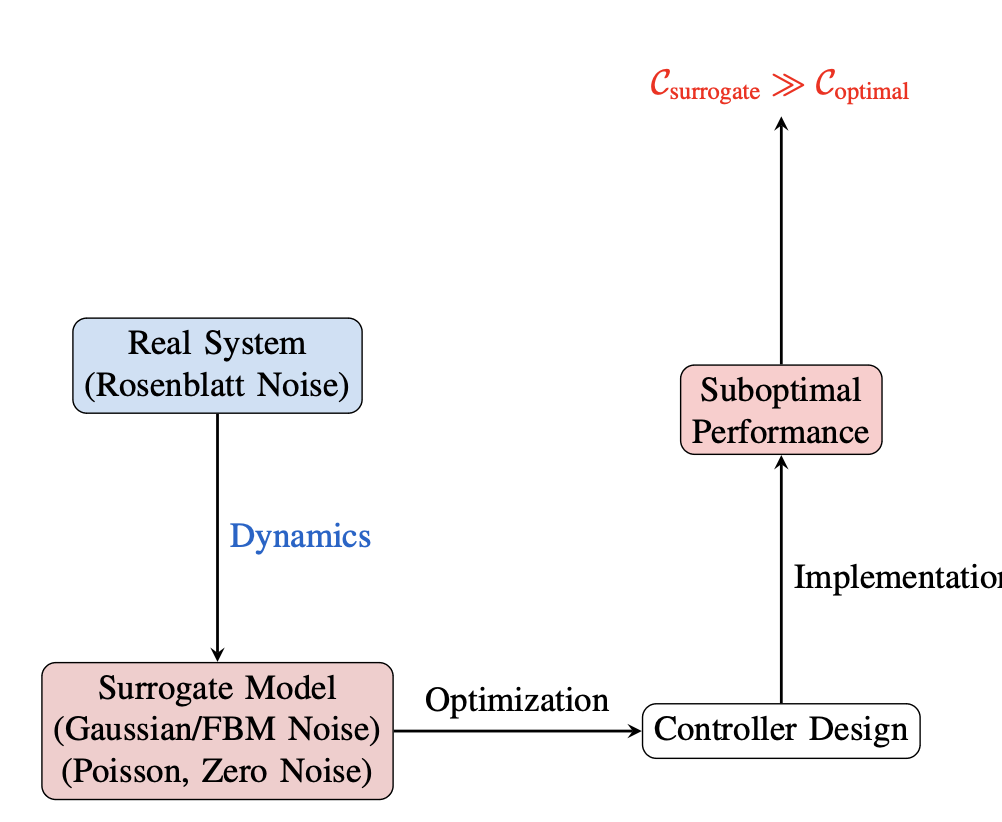}   \includegraphics[scale=0.3]{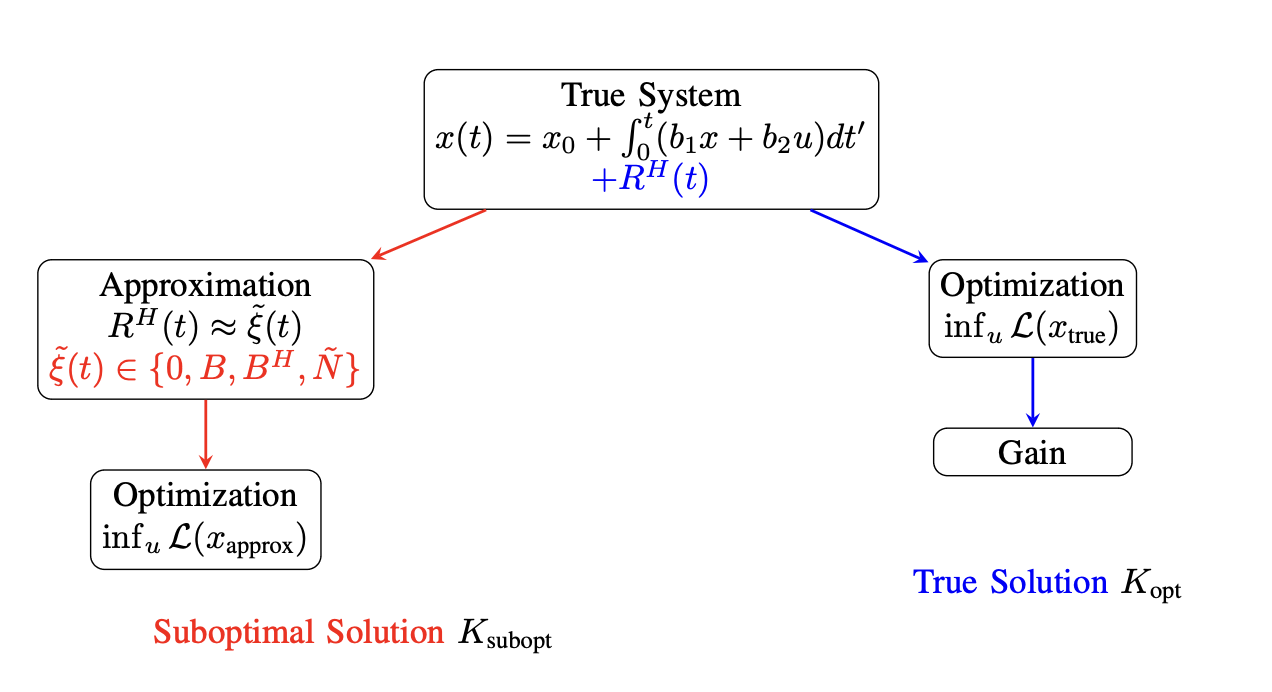} 
 \caption{Failure of Noise Approximation-Optimization Cascade}
 \label{fig16sub}
\end{figure}

\begin{figure} %{Cost of simplicity }
\centering
   \includegraphics[scale=0.3]{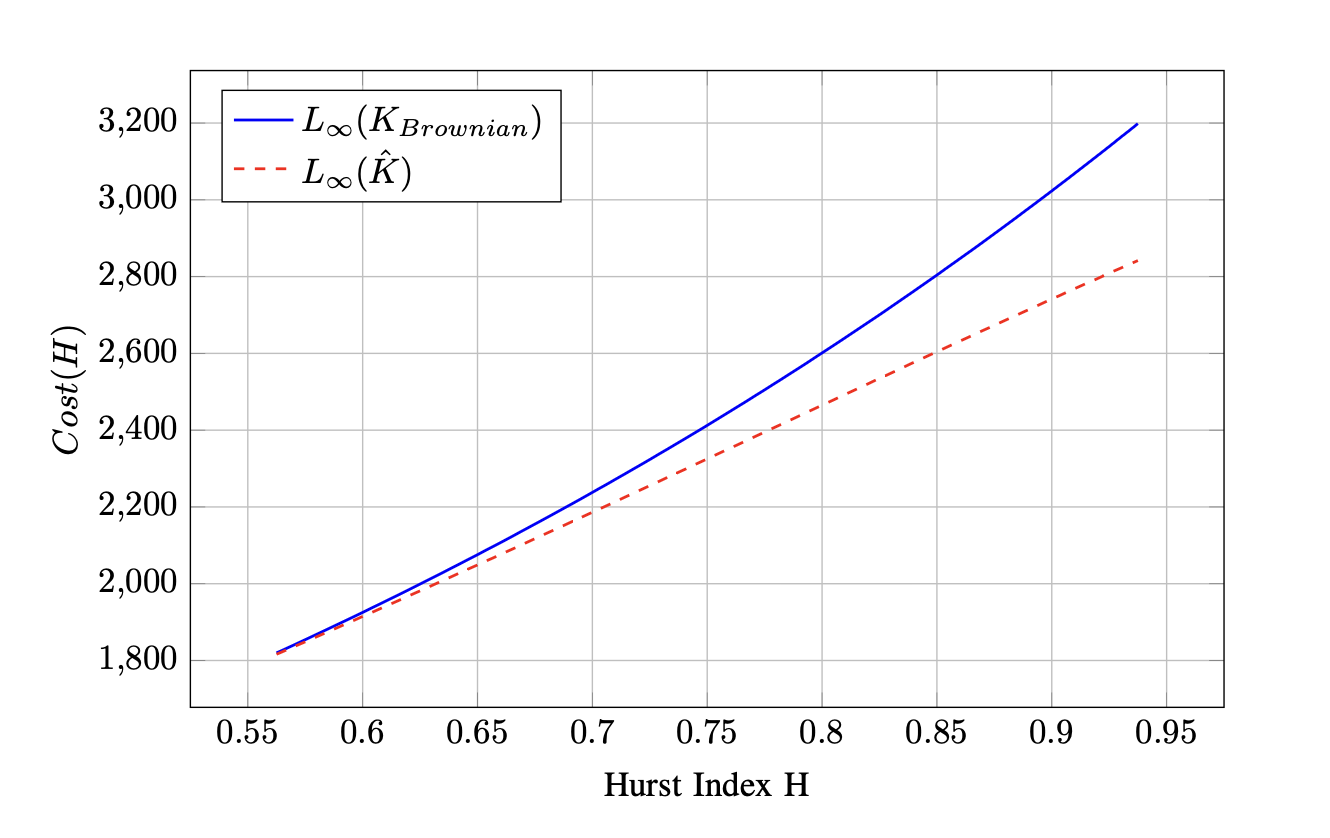} 
    \caption{Failure of Noise Approximation-Optimization Cascade under Rosenblatt process}
 \label{fig18sub}
\end{figure}

\begin{itemize}
\item \textcolor{process}{Rosenblatt Properties}: Non-Gaussian, long-range dependence
\item \textcolor{approx}{Surrogate Pitfalls}:  $\exists \epsilon >0 : |\mathcal{C}_{\text{surrogate}} - \mathcal{C}_{\text{true}}| \geq \epsilon$. See Figures \ref{fig16sub}, and \ref{fig18sub}
%Loss of temporal correlations, moment structure
\end{itemize}

  \subsection{Two extensions }
Two additive Rosenblatt noises:           
\begin{equation}\nonumber \begin{array}{ll}
L_\infty(u) = \limsup_{T \rightarrow \infty} \frac{1}{T} \mathbb{E} \int_0^T (q x^2(t) + r u^2(t)) dt, \\
  x(t) = x_0 + \int_0^t (a x(t') + b u(t')) dt' +{\color{red} \sigma_1 R^{H_1}_1(t)+\sigma_2 R^{H_2}_2(t)}, \\
  (a,b,x_0)\in  \mathbb{R}^3, b\neq 0, \ q>0, r>0,  H_i \in (\frac{1}{2}, 1),  \sigma_i\neq 0  \\
                  \mathcal{U} := \{(u(t))_{t \geq 0}: \  u(t) = K x(t) \text{ with } K \in \mathbb{R}\}
                  \end{array}
\end{equation}
                  
                    Control of the ergodic cost under two  Rosenblatt processes:  $\inf_{u\in \mathcal{U}}L_\infty(u)$ subject to the state dynamics driven by two Rosenblatt processes (independent ? correlated ?)

                Two multiplicative Rosenblatt noises:             
\begin{equation}\nonumber \begin{array}{ll}
L_\infty(u) = \limsup_{T \rightarrow \infty} \frac{1}{T} \mathbb{E} \int_0^T (q x^2(t) + r u^2(t)) dt, \\
  x(t) = x_0 + \int_0^t (b_1 x(t') + b_2 u(t')) dt' +{\color{red} \int_0^t \sigma_1 x(t') R^{H_1}_1(t')+\sigma_2 x(t') R^{H_2}_2(t')}, \\
  (b_1,b_2,x_0)\in  \mathbb{R}^3, b_2\neq 0, \ q>0, r>0,  H_i \in (\frac{1}{2}, 1),  \sigma_1,\sigma_2\neq 0  \\
                  \mathcal{U} := \{(u(t))_{t \geq 0}: u(t) = K x(t) \text{ with } K \in \mathbb{R}\}
                  \end{array}
\end{equation}
                  
                     Control of the ergodic cost under two  Rosenblatt processes:  $\inf_{u\in \mathcal{U}}L_\infty(u)$ subject to the state dynamics driven by two Rosenblatt processes (independent ? correlated ?)

         %\subsection{  Infinite Horizon with Discounted Cost } 
          %\subsection{  Finite Horizon } 

\section{Adaptive Control  driven by Rosenblatt noise  }

\section{Variance-Aware Control  driven by Rosenblatt noise  }
We examine the infinite horizon variance-aware  optimal control with long-time average cost and linear state dynamics  driven by Rosenblatt noise process.          
\begin{equation}\nonumber \begin{array}{ll}
L_\infty(u) = \limsup_{T \rightarrow \infty} \frac{1}{T} \int_0^T ({\color{blue} q \ var(x(t)) + r \ var(u(t)) }+ \bar{q} \bar{x}^2(t)+ \bar{r} \bar{u}^2(t) ) dt, \\
  x(t) = x_0 + \int_0^t (b_1 x(t') + b_2 u(t')) dt' + {\color{blue} R^H(t)} + \int_0^t (\bar b_0+\bar b_1 \bar x(t') + \bar b_2 \bar u(t')) dt' , \\
  (b_1,b_2,x_0)\in  \mathbb{R}^3, b_2\neq 0, \ q, \bar{q}>0,  r, \bar{r}>0,  H \in (\frac{1}{2}, 1),   \\
  (\bar b_0,\bar b_1,\bar b_2)\in  \mathbb{R}^3,  \bar b_1 \neq 0\\
                  \mathcal{U} := \{(u(t))_{t \geq 0}: u(t) = K (x(t)-\bar{x}(t))+ \bar K \bar{x}(t)  \text{ with }  (K,\bar K) \in \mathbb{R}^2, b_1+K b_2 <0, b_1+\bar b_1+\bar K (b_2+\bar b_2) <0\}
                  \end{array}
\end{equation}
                  
                 The variance-aware control of the ergodic cost under Rosenblatt is given by   $\inf_{u\in \mathcal{U}}L_\infty(u)$ subject to the McKean-Vlasov linear state dynamics driven by Rosenblatt noise.

          We use a direct method to find linear  state-and-mean-field-type  feedback optimal control.
          
            \begin{thm}
The optimal gains 
  in the family of admissible state-mean-field-type feedbacks are given by

$$
{K}^* = -\frac{b_1+ \sqrt{b_1^2 + 4H(1-H) \frac{b_2^2q}{r}}}{2b_2(1-H)},
$$

$$
\bar{K}^* = \frac{(b_2+\bar b_2) \bar{q}}{ \bar{r} (b_1+\bar b_1)}, \   b_1+\bar{b}_1\neq 0,
$$

and the optimal cost is given by

$
L_\infty (K^*,\bar{K}^*) = \frac{\Gamma(2H)}{[-(b_1 + b_2\hat{K})]^{2H-1}} \left(-\frac{r\hat{K}}{b_2}\right) + 
\frac{\bar{b}_0^2 \bar{q} \bar{r}  }{( b_1+\bar b_1)^2 \bar{r} + (b_2+\bar b_2)^2\bar{q}}  $ 
%%%\frac{\bar{b}_0^2 (\bar{q}+ \bar{r} \bar{K}^2) }{(\bar{b}_1 + \bar{b}_2\bar{K})^2}  $ 
  \end{thm}

  \subsection{
What if we replace the Rosenblatt by a standard one  and Minimize the variance ?   }       
\begin{equation}\nonumber 
\begin{array}{ll}
\inf_{u\in \mathcal{U} }  \limsup_{T \rightarrow \infty} \frac{1}{T} \int_0^T ({\color{blue} q \ var(x(t)) + r \ var(u(t)) }+ \bar{q} \bar{x}^2(t)+ \bar{r} \bar{u}^2(t) )  dt, \\
%x(t) = x_0 + \int_0^t (ax(t') + bu(t')) dt' +  \textcolor{red}{ \cancel{R^H(t)} }\\
 x(t) = x_0-\bar x_0 + \int_0^t (b_1 (x(t')-\bar x(t')) + b_2 (u(t')-\bar u(t'))) dt'  \\ +
 \bar x_0+ \int_0^t (\bar b_0+ ( b_1+\bar b_1) \bar x(t') + (b_2+\bar b_2) \bar u(t')) dt'+ {\color{blue} R^H(t)} , \\
 \qquad  \qquad  \qquad  \qquad \qquad  \qquad\qquad \qquad \qquad\qquad \qquad \approx \textcolor{blue}{B(t)} \\
  \qquad  \qquad  \qquad  \qquad \qquad  \qquad\qquad\qquad \qquad\qquad \qquad \approx \textcolor{blue}{B^H(t)} \\
  \qquad  \qquad  \qquad  \qquad \qquad  \qquad\qquad\qquad \qquad\qquad \qquad \approx \textcolor{blue}{\tilde N(t,.)} \\
   \qquad  \qquad  \qquad  \qquad \qquad  \qquad\qquad \qquad \qquad\qquad \qquad \approx \textcolor{blue}{0}  \\
\end{array}
\end{equation}
This operation leads to suboptimality.
  Approximating the Rosenblatt noise by another one such as Brownian motion, fractional Brownian, Gauss-Volterra  noise, compensated Poisson jump process leads to \textcolor{red}{suboptimal control}.
  
\section{Zero-Sum Games   driven by Rosenblatt noise  }
 We formulate  and solve  zero-sum games where the state dynamics are driven by Rosenblatt noise. These games are particularly relevant in adversarial settings, where two decision-makers with opposing objectives interact within a stochastic environment. The key result is the establishment of saddle-point equilibria, providing conditions under which each decision-maker can adopt an optimal strategy that minimizes their maximum possible loss. This is a crucial development, as it extends classical game theory to more accurately reflect the stochastic nature of real-world adversarial interactions.

Let \begin{equation}\nonumber \begin{array}{ll}
 L_T(u, v) = \frac{1}{T} \mathbb{E} \left[ \int_0^T \left(   q x^2(t)+ r u^2(t)-s v^2(t)  \right) dt \right], \\ L_{1,2,\infty}(u, v) = \liminf_{u, T \to \infty} \limsup_{v, T \to \infty} L_T(u, v)\\
L_{2,1,\infty}(u, v) = \limsup_{v, T \to \infty} \liminf_{u, T \to \infty} L_T(u, v)\\  \\
  x(t) =x_0+\int_0^t [ b_1x(t')+ b_2 u(t') + b_3v(t')]dt' + {\color{blue}R^H(t)}, \\ \\
  (b_1,b_2,b_3,x_0)\in  \mathbb{R}^4, \  b_2,b_3\neq 0, \ q,r,s>0, H \in (\frac{1}{2}, 1),\\
                    \end{array}
\end{equation}

\begin{equation}\nonumber \begin{array}{ll}

  \limsup_{T \to \infty} \frac{1}{T} \mathbb{E} \left[ \int_0^T \left( 
 q x^2(t)+ r u^2(t)+s v^2(t)  \right) dt \right] < \infty\\
  s^2 b_1^2 b_2^2 \geq 4 \left(- \frac{b^3 s^2}{ b_3 r} + b_2 s b_3 + \frac{H b^3_2 s^2}{ b_3^2} - H b_2 s \right) \left( s b_1 b_2 + H b_2 q \right)\\
                  \mathcal{U} := \{(u(t))_{t \geq 0}: u_t = Kx(t) \text{ with } K \in \mathbb{R}\} \\   \mathcal{V} := \{(v(t))_{t \geq 0}: v(t) = L x(t) \text{ with } L \in \mathbb{R}\}
                  \end{array}
\end{equation}
                  
                     Zero-Sum Game with ergodic cost under Rosenblatt noise is   $\inf_{u\in \mathcal{U}}\sup_{v\in \mathcal{V}}L_\infty(u, v)$ subject to the state dynamics driven by Rosenblatt noise \cite{examplebook4}.

                     \begin{thm}
                MiniMax=Maximin:  The zero-sum game under Rosenblatt has a value and each decision-maker has an optimal strategy.
                     \end{thm}
                   
                     \begin{thm}
         Explicit Solution of MiniMax:     The saddle-point gains  
  in the family of admissible feedbacks are given by

$$
K = -\frac{b_2 s}{b_3 r} L, \quad (-\frac{b^3_2 s^2}{ b_3 r} + b_2 s b_3 + \frac{H b^3_2 s^2}{ b_3^2} - H b_2 s) L^2 + s b_1 b_2 L + H b_2 q = 0,
$$
and the equilibrium cost (value) is given by
$
L_\infty(K, L) = \frac{\Gamma(2H+1)(q + rK^2 - sL^2)}{2[-(b_1 + b_2K + b_3L)]^{2H}}
$ 
  \end{thm}
  
Note that the extension of a bigger class of strategies is still an open issue.
 
   \subsection{  
What if we replace the Rosenblatt by a standard one  and find Saddle Point}            
\begin{equation}\nonumber 
\begin{array}{ll}
 \inf_{u} \sup_{v} \limsup_{T \rightarrow \infty} \\ \frac{1}{T}\mathbb{E} \left[ \int_0^T \left( q x^2(t)+q u^2(t)-s v^2(t)
   \right) dt \right], \\
x(t) = x_0 + \int_0^t (b_1 x(t') + b_2u(t')+b_3v(t')) dt' +  \textcolor{red}{ \cancel{R^H(t)} }\\
 \qquad  \qquad  \qquad  \qquad \qquad  \qquad \qquad \qquad \qquad \qquad  \approx \textcolor{blue}{B(t)} \\
  \qquad  \qquad  \qquad  \qquad \qquad  \qquad \qquad \qquad \qquad \qquad \approx \textcolor{blue}{B^H(t)} \\
  \qquad  \qquad  \qquad  \qquad \qquad  \qquad \qquad \qquad \qquad  \qquad \approx \textcolor{blue}{\tilde N(t)} \\
  \qquad  \qquad  \qquad  \qquad \qquad  \qquad \qquad \qquad \qquad  \qquad \approx \textcolor{blue}{0}. \\
\end{array}
\end{equation}
  This operation leads to non-equilibrium.
  Approximating the Rosenblatt noise by another one such  as Brownian motion, fractional Brownian, Gauss-Volterra  noise, compensated Poisson jump process  \textcolor{red}{does not provide the saddle point}.

\section{Non-Zero-Sum Games  driven by Rosenblatt noise  }

We now consider  non-zero-sum games under Rosenblatt noise  where the interests of the  decision-makers  are not strictly opposed, but rather involve a mix of cooperation and competition. Unlike zero-sum games, where one decision-maker's gain is exactly balanced by the other's loss, non-zero-sum games allow for the possibility of mutual benefit or shared losses. This  is crucial because it extends the analysis of game theory to more realistic scenarios where decision-makers' objectives are interdependent and influenced by non-Gaussian noise.

\begin{equation}\nonumber \begin{array}{ll}
L_{i,\infty}(u) = \limsup_{T \rightarrow \infty} \frac{1}{T} \mathbb{E} \int_0^T (q_i x^2(t) + r_i u_i^2(t)) dt, \\
  x(t) = x_0 + \int_0^t (b_1 x(t') + \sum_{j\in \mathcal{I}}b_{2j} u_j(t')) dt' + {\color{blue} R^H(t)}, \\
  (b_1,b_{2j},x_0)\in  \mathbb{R}^3, b_{2j}\neq 0, \ q_j>0, r_j>0,  H \in (\frac{1}{2}, 1),   \\
                  \mathcal{U}_i := \{(u_i(t))_{t \geq 0}: u_i(t) = K_i x(t) \text{ with } K_i \in \mathbb{R}, (b_1+ \sum_{j\in \mathcal{I}}b_{2j} K_j) <0\}
                  \end{array}
\end{equation}
                  
                 The game  with the ergodic cost under Rosenblatt noise is given by   $\inf_{u_i\in \mathcal{U}_i}L_{i,\infty}(u)$ subject to the common state dynamics driven by Rosenblatt noise.

          We use a direct method to find linear  state-feedback best-response control strategy of agent $i$.
          
            \begin{thm} Given the admissible strategies of the other agents  $(u_j\in \mathcal{U}_j, \  j\neq i),$
           the best-response gain of agent $i$ 
  in the family of admissible feedbacks is given by
$
{K}_i^* = -\frac{a + \sqrt{a^2 + 4H(1-H) \frac{b^2_{2i}q_i}{r_i}}}{2b_{2i} (1-H)},
$
with $ a:= b_1+ \sum_{j\in \mathcal{I}\backslash \{i\}}b_{2j} K_j,$
and the best-response cost is given by
$L_{i,\infty} ({K}^*) = \frac{\Gamma(2H)}{[-(b_1+ \sum_{j\in \mathcal{I}}b_{2j} K_j)]^{2H-1}} \left(-\frac{r_i {K}^*_i}{b_{2i}}\right)$ 
  \end{thm}

            \begin{thm} If  there  exists a solution to the following system:  
            \begin{equation}\nonumber \begin{array}{ll}
            i\in \mathcal{I},\\ 
       {K}_i = -\frac{b_1+ \sum_{j\in \mathcal{I}\backslash \{i\}}b_{2j} K_j + \sqrt{(b_1+ \sum_{j\in \mathcal{I}\backslash \{i\}}b_{2j} K_j)^2 + 4H(1-H) \frac{b^2_{2i}q_i}{r_i}}}{2b_{2i} (1-H)}, \\ (b_1+ \sum_{j\in \mathcal{I}}b_{2j} K_j) <0,
         \end{array}
\end{equation}
 then, the non-zero-sum stochastic differential game under Rosenblatt has an equilibrium in stationary linear state-feedback strategies.
  \end{thm}
  
  Note that this result is  restricted  to the class of linear state-feedback strategies.  The extension of  non-linear strategies is an open issue.

  \subsection{
What if we replace the Rosenblatt by a standard  Brownian motion and solve the non-cooperative game ?   }       
\begin{equation}\nonumber 
\begin{array}{ll}
\inf_{u\in \mathcal{U} }  \limsup_{T \rightarrow \infty} \frac{1}{T} \mathbb{E} \int_0^T (q_ix^2(t) + r_iu^2_i(t)) dt, \\
  x(t) = x_0 + \int_0^t (b_1 x(t') + \sum_{j\in \mathcal{I}}b_{2j} u_j(t')) dt' + {\color{blue} R^H(t)}, \\
 \qquad  \qquad  \qquad  \qquad \qquad  \qquad\qquad \qquad \qquad \qquad \approx \textcolor{blue}{B(t)} \\
  \qquad  \qquad  \qquad  \qquad \qquad  \qquad\qquad\qquad \qquad \qquad \approx \textcolor{blue}{B^H(t)} \\
  \qquad  \qquad  \qquad  \qquad \qquad  \qquad\qquad\qquad \qquad \qquad \approx \textcolor{blue}{\tilde N(t,.)} \\
   \qquad  \qquad  \qquad  \qquad \qquad  \qquad\qquad \qquad \qquad \qquad \approx \textcolor{blue}{0}  \\
\end{array}
\end{equation}
This operation leads to non-equilibrium.
  Approximating the Rosenblatt noise by another one such as  Brownian motion, fractional Brownian, Gauss-Volterra  noise, compensated Poisson jump process no longer provide an \textcolor{red}{equilibrium strategy}.

\section{ Mean-Field-Type Games driven by Rosenblatt noise  }
{ What is Mean-Field-Type Game Theory?} 
   \begin{table}[h!]  % 
 \begin{tabular}{|  p{3.5cm} |  p{2.8cm} |  p{5.5cm} | }   
 \hline
 Quantities-of-Interest  & classical &  MFTG   \\ \hline
Instant Payoffs  &  time-state-action&  time-state-action{ \color{blue} +   (probability) measure of state-action }\\  \hline
 Instant Coefficients  & time-state-action   & time-state-action { \color{blue} +   (probability) measure of state-action }\\ \hline
 \end{tabular}
  \end{table}

  Following \cite{examplebookt1}, the  MFTG driven by  Rosenblatt noise can be defined. 
 Given $(u_j, \ j\in \mathcal{I}\backslash \{i\})$ in linear-and-mean-field-type feedback form, the best-response problem of the decision-maker $i$  is 
\begin{equation}
\nonumber
\begin{array}{ll}
\inf_{u_i}  \mathbb{E}  [ L_i(x,u)\ | x(0) = x_0 ],\\
\mbox{subject to } \\
x(t)= x_0+\int_0^t  [ b_1 (x-\bar{x})+\sum_{j\in \mathcal{I}}b_{2j} (u_j-\bar{u}_j) ] dt'  \\
 +\int_0^t [(b_1+\bar{b}_1) \bar{x}+\sum_{j\in \mathcal{I}}({b}_{2j} +\bar{b}_{2j} )\bar{u}_j] dt' 
 {\color{blue}+   R^H(t)}, 
\end{array}
\end{equation} 

          \begin{equation}  \nonumber
\begin{array}{ll}
L_i (x,u) = q_{iT}{\color{blue}(x_T-\bar{x}_T)^{2} } + \bar{q}_{iT} \frac{\bar{x}_T^{2\bar{k}_i}}{2\bar{k}_i} \\ 
+\int_0^T   \left[ q_{i}(x-\bar{x})^{2} +{\color{blue} {r}_i  (u_i-\bar{u}_i)^{2}}
+   \bar{q}_{i}\frac{\bar{x}^{2\bar{k}_i}}{2\bar{k}_i} + \bar{r}_i  \frac{\bar{u}_i^{2\bar{k}_i}}{2\bar{k}_i} \right]\ dt,
\end{array}
\end{equation}
with the (deterministic) coefficients
\begin{equation} \nonumber
\begin{array}{ll}
q_{i}, \bar{q}_{i}:  \ [0,T] \rightarrow \mathbb{R}_{+},\  i\in \mathcal{I},\\
{r}_i,\bar{r}_i>\epsilon >0:   \  [0,T]\ \rightarrow \mathbb{R}_{+},\  i\in \mathcal{I},\\
b_1, \bar{b}_1, b_{2j}, \bar{b}_{2j}:  \ [0,T] \rightarrow \mathbb{R},\  j\in \mathcal{I},\\
%%\sigma, \sigma_{gv}:\ [0,T]\times \mathcal{S} \rightarrow \mathbb{R},\\
\end{array}
\end{equation}
The functions $ b_{2j}, \bar{b}_{2j}$ are not identically zero.

We use a direct method to find linear  state-and-mean-field-type feedback Nash equilibria.

Given the strategies $(u_j, \ j\neq i),$  of  decision-makers other than $i,$ a  best response strategy of decision-maker $i$ is a 
strategy that solves the  following  problem of mean-field type
 $ \inf_{u_i} \  \mathbb{E} [ L_i] \ $ subject to  the state.
 The set of  best responses of $i$ is  denoted by $\mbox{BR}_i ((u_{j})_{j\neq i}).$

Based on the  best-response  function $BR=(BR_i)_{i\in \mathcal{I}},$ we define a Nash Equilibrium of the MFTG.
A (mean-field-type)  Nash equilibrium  is a strategy profile $(u^*_j, \ j\in \mathcal{I}),$ of all decision-makers such that for every decision-maker $i,$  
$u_i^*\in \mbox{BR}_i ((u^*_{j})_{j\neq i}).$ 

\begin{thm}
  \label{semi:prop9t1t} Let $q_i>0, \bar{q}_i>0, \ r_i, \bar{r}_i>\epsilon>0. $
A Nash equilibrium of the  nonzero-sum MFTG  with a Rosenblatt process,
 having a linear state-and- mean-field type feedback form
is  given by 
\begin{equation}
\begin{array}{ll}
u^*_i =-  \eta_i (x-\bar{x})- \bar{\eta}_i  \bar{x},\ \\ 
\bar{\eta}_i=  \left( \frac{\bar{b}_{2i}\bar{\lambda}_i}{\bar{r}_i} \right)^{\frac{1}{2\bar{k}_i-1}},\\
 \mathbb{E} [ L_i]=   \lambda_i(0)   \mathbb{E} (x_0-\bar{x}_0)^2+ \bar{\lambda}_i(0) \frac{1}{2\bar{k}_i}\bar{x}^{2\bar{k}_i}_0 +\gamma_i(0),
\end{array}
\end{equation}
whenever the following system of  differential equations admits a positive solution which does not blowup (no singularity)  within 
$[0,T]:$
\begin{equation} \nonumber
\begin{array}{ll}
-\dot{\lambda}_i  =q_{i}+ 2 b_1 \lambda_i -   \frac{b^2_{2i}}{r_i}\lambda^2_i, \\
{\lambda}_i(T)  =q_{i}(T),\\ \\

-\frac{d}{dt}\bar{\lambda}_i = \bar{q}_{i}+ 2\bar{k}_i ( b_1+\bar{b}_1) \bar{\lambda}_i
-   2\bar{k}_i  \bar{\lambda}_i \sum_{j\neq i}(b_{2j}+\bar{b}_{2j}) \left( \frac{(b_{2j}+\bar{b}_{2j})\bar{\lambda}_j}{\bar{r}_j} \right)^{\frac{1}{2\bar{k}_j-1}} \\ 
-(2\bar{k}_i-1) \bar{r}_i \left( \frac{(b_{2i}+\bar{b}_{2i})\bar{\lambda}_i}{\bar{r}_i} \right)^{\frac{2\bar{k}_i}{2\bar{k}_i-1}} ,\\
\bar{\lambda}_i(T)  =\bar{q}_{i}(T),\\ \\
 - \dot{\gamma}_i = {r}_i   (\eta_i  - \lambda_i \frac{b_{2i}}{r_i} )^2 v_2(t)-2\lambda_i \sum_{j\in \mathcal{I}\backslash \{i\}}b_{2j} \eta_j   v_2(t)  
  +2c \lambda_i o(t),     \\
   {\gamma}_i(T)=0,
\end{array}
\end{equation}
where $\eta_i$  minimizes   $ f_i(\eta):$  
\begin{equation}
\begin{array}{ll}
f_i(\eta)=  {r}_i   (\eta_i  - \lambda_i \frac{b_{2i}}{r_i} )^2 v_2(t)-2\lambda_i \sum_{j\in \mathcal{I}\backslash \{i\}}b_{2j} \eta_j   v_2(t)   
+2c \lambda_i o(t), 
\end{array}
\end{equation}
and the vector $\eta$ is a fixed-point of   $(\arg\min_{z_i}f_i(z_i,\eta_{-i}))_{i\in \mathcal{I}},$  
 \begin{equation} \nonumber
\begin{array}{ll}
b=[b_1- \sum_{j}b_{2j} \eta_j],\\
   o(t)=      \left[ \int_0^t   c_3 [t-s]^{2H-2}  e^{-\int_0^s  b(t'')  dt''} ds \right] e^{\int_{0}^t  b(t') \ dt'} , \\
  v_2(t)= \left[v_2(0)+ \int_0^t  2 c \ o(s) e^{-\int_0^s  2 b(t'')  dt''} ds \right] e^{\int_{0}^t 2 b(t') \ dt'} .
   \end{array}
\end{equation}
\end{thm}
%%%

The following remarks are in order:
 \begin{itemize}
 \item The equilibrium strategies obtained for the noise-free case (purely deterministic case)  
  no longer provide an equilibrium when the noise is driven by a  Rosenblatt process. 
 \item The equilibrium strategies obtained for the Brownian motion, multi-fractional Brownian, Gauss-Volterra or a pure jump and regime switching  processes  no longer provide an equilibrium when the noise is driven by a  Rosenblatt process. 
\item  In particular, the   equilibrium strategies for the case when the state is driven by  Brownian motion  are no longer in equilibrium even in the quadratic cost case $\bar{k}_i=1$  for which the third  derivatives are zero. 
\end{itemize}

   \subsection{Fully Cooperative MFTG}
\begin{thm}
\begin{equation}
\nonumber
\begin{array}{ll}  \omega_j >0, {j\in \mathcal{I}} \\
\inf_{(u_i)_{i\in \mathcal{I}}}  \mathbb{E} [\sum_{j\in \mathcal{I}}\omega_j L_j(x,u) ],\\
\mbox{subject to } \\
x(t)= x_0+\int_0^t  [ b_1 (x-\bar{x})+\sum_{j\in \mathcal{I}}b_{2j} (u_j-\bar{u}_j) ] dt'  \\
 +\int_0^t [(b_1+\bar{b}_1) \bar{x}+\sum_{j\in \mathcal{I}}({b}_{2j} +\bar{b}_{2j} )\bar{u}_j] dt' 
 +   {\color{blue} R^H(t)}, 
\end{array}
\end{equation} 

 Then, the optimal weighted cooperative cost is
\begin{equation} \nonumber
\begin{array}{ll}
[\sum_{i\in \mathcal{I}}\omega_i   \mathbb{E} L_i  -\lambda(0)  \mathbb{E} (x_0-\bar{x}_0)^{2}- \bar{\lambda}(0)\frac{\bar{x}_0^{2\bar{k}}}{2\bar{k}} ]\\
=
 \int_0^T    \sum_{i\in \mathcal{I}} \omega_i {r}_i   (\eta_i  -\frac{\lambda b_{2i}}{\omega_i {r}_i  })^2 v_2
 +2c \lambda o(t)\ dt, \\ \mbox{the coefficients } \lambda,\bar{\lambda} \ \mbox{solve}\
 \\
 \lambda(T)=\sum_{i}\omega_i q_{iT} ,\\
 0= \dot{\lambda} + \sum_{i}\omega_iq_{i}+ 
2\lambda b_1  - \sum_{i\in \mathcal{I}}\frac{ b_{2i}^2}{\omega_i {r}_i  }\lambda^2,\\
 \bar{\lambda}(T)=\sum_{i\in \mathcal{I}} \omega_i \bar{q}_{iT} ,\\
\frac{d}{dt}\bar{\lambda} + \sum_{i\in \mathcal{I}} \omega_i\bar{q}_{i}+  2\bar{k} \bar{\lambda} (b_1+\bar{b}_1)    \\
-
 (2\bar{k}-1)  \sum_{i\in \mathcal{I}}\omega_i\bar{r}_i \left( \frac{(b_{2i}+\bar{b}_{2i})\bar{\lambda}}{\omega_i\bar{r}_i} \right)^{\frac{2\bar{k}}{2\bar{k}-1}} =0,\\
\mbox{and} \ \eta, \bar{\eta} \ \mbox{are given by } 
 \\
 \bar{\eta}^{coop}_i=  \left( \frac{(b_{2i}+\bar{b}_{2i})\bar{\lambda}}{\omega_i\bar{r}_i} \right)^{\frac{1}{2\bar{k}-1}},\
 \bar{u}^{coop}_i= -\bar{\eta}^{coop}_i \bar{x},\\
 (\eta_i)_{i\in \mathcal{I}} \ \mbox{minimizes} \ 
 \sum_{i\in \mathcal{I}} \omega_i {r}_i   (\eta_i  -\frac{\lambda b_{2i}}{\omega_i {r}_i  })^2 v_2 +2c \lambda o(t),\\
 u^{coop}_i= -\eta_i (x-\bar{x})-\bar{\eta}_i \bar{x},\\

\end{array}
\end{equation}  
\end{thm}

\subsection{ Noise Modelling in MFTG}
We look at  five different noises in a particular class of mean-field-type games.
\begin{equation}  \nonumber
\begin{array}{ll}
L_i= q_{iT}\frac{(x_T-\bar{x}_T)^{2}}{2}  + \bar{q}_{iT} \frac{\bar{x}_T^{2\bar{k}_i}}{2\bar{k}_i}\\
+\int_0^T q_{i}\frac{(x-\bar{x})^{2}}{2} +{r}_i  \frac{(u_i-\bar{u}_i)^{2}}{2} 
+  \bar{q}_{i}\frac{\bar{x}^{2\bar{k}_i}}{2\bar{k}_i} + \bar{r}_i  \frac{\bar{u}_i^{2\bar{k}_i}}{2\bar{k}_i}  \ dt  \\
+ \int_0^T  \bar{c}_{i} \bar{x}^{2\bar{k}_i-1} \bar{u}_i+  \sum_{j\in \mathcal{I}\setminus \{i\}}\bar{\epsilon}_{ij} \bar{x}^{2(\bar{k}_i-1)}\bar{u}_i\bar{u}_j \ dt , \\
\mbox{subject to } \\
x(t)= x_0 +  \int_0^t  b_1 x +\bar{b}_1 \bar{x} +\sum_{j\in \mathcal{I}} \left[b_{2j} u_j+ \bar{b}_{2j} \bar{u}_j\right]  \ dt' \\
+  x_{noise}(t)
 \\ \\
\mathbb{P}(s(t+\delta)=s' | s)=\int_t^{t+\delta} \tilde{q}_{ss'} dt' + o(\delta),\ s'\neq s,\\
\ s(0)=s_0,  
\end{array}
\end{equation}

When the  state is driven by Brownian motion, the MFTG becomes:
\begin{equation}  \nonumber
 (PBm)  \left\{ \begin{array}{ll}
\mbox{Find }  (u^*_1,\ldots,u^*_I):  \\
i \in \mathcal{I},  \\
\mathbb{E} [L_i (u^*, m_0) ]\\ = \inf_{u_i\in \mathcal{U}_i} \mathbb{E}[  L_i (u^*_1,\ldots, u^*_{i-1}, u_i, u^*_{i+1},\ldots,u^*_I, m_0)],  \\
  x_{noise}(t)=\int_0^t  {\sigma}_b dB(t'),
  \\   \mathcal{S}= \{s_0\}
\end{array} \right.
\end{equation}
 
When the state is driven by Regime Switching the MFTG becomes: 
\begin{equation}  \nonumber
%\label{stackros1generic2}
 (PRS)  \left\{ \begin{array}{ll}
\mbox{Find }  (u^*_1,\ldots,u^*_I):  \\
i \in \mathcal{I},  \\
\mathbb{E} [L_i (u^*, m_0) ]\\ = \inf_{u_i\in \mathcal{U}_i} \mathbb{E} [  L_i (u^*_1,\ldots, u^*_{i-1}, u_i, u^*_{i+1},\ldots,u^*_I, m_0)],  \\
  x_{noise}(t)=0, 
  \\    |\mathcal{S}| \geq  2
\end{array} \right.
\end{equation}

When the state is driven by compensated Poisson Jump process, the MFTG becomes:
\begin{equation}  \nonumber
 (PPJ)  \left\{ \begin{array}{ll}
\mbox{Find }  (u^*_1,\ldots,u^*_I):  \\
i \in \mathcal{I},  \\
\mathbb{E} [L_i (u^*, m_0) ]\\ = \inf_{u_i\in \mathcal{U}_i} \mathbb{E}[  L_i (u^*_1,\ldots, u^*_{i-1}, u_i, u^*_{i+1},\ldots,u^*_I, m_0)],  \\
 x_{noise}(t)=\int_0^t  \int_{\theta} \sigma_n(t',s,\theta)\tilde{N}(dt',d\theta) ,
  \\   \mathcal{S}= \{s_0\}
\end{array} \right.
\end{equation}

When the state is driven by Gauss-Volterra noise process: 
\begin{equation}  \nonumber
 (PGV)  \left\{ \begin{array}{ll}
\mbox{Find }  (u^*_1,\ldots,u^*_I):  \\
i \in \mathcal{I},  \\
\mathbb{E} [L_i (u^*, m_0) ]\\ = \inf_{u_i\in \mathcal{U}_i} \mathbb{E}[  L_i (u^*_1,\ldots, u^*_{i-1}, u_i, u^*_{i+1},\ldots,u^*_I, m_0)],  \\
 x_{noise}(t)=\int_0^t  {\sigma}_{gv} dB_{gv}(t') ,
  \\   \mathcal{S}= \{s_0\}
\end{array} \right.
\end{equation}

When the  state is driven by Rosenblatt process the MFTG becomes: 
\begin{equation} \nonumber
 (PRh)  \left\{ \begin{array}{ll}
\mbox{Find }  (u^*_1,\ldots,u^*_I):  \\
i \in \mathcal{I},  \\
\mathbb{E} [L_i (u^*, m_0) ]\\ = \inf_{u_i\in \mathcal{U}_i} \mathbb{E}[  L_i (u^*_1,\ldots, u^*_{i-1}, u_i, u^*_{i+1},\ldots,u^*_I, m_0)],  \\
 x_{noise}(t)=\int_0^t {\sigma}_{r}  {\color{blue} dR^H(t')} ,
  \\   \mathcal{S}= \{s_0\}
\end{array} \right.
\end{equation}

Following \cite{examplebookt2}, the  MFTG with state driven by five mixed noises becomes

\begin{equation} \nonumber
 (P5mixed) \left\{ \begin{array}{ll}
\mbox{Find }  (u^*_1,\ldots,u^*_I):  \\
i \in \mathcal{I},  \\
\mathbb{E} [L_i (u^*, m_0) ]\\ = \inf_{u_i\in \mathcal{U}_i} \mathbb{E}[  L_i (u^*_1,\ldots, u^*_{i-1}, u_i, u^*_{i+1},\ldots,u^*_I, m_0)],  \\
  x_{noise}(t)= \int_0^t   {\sigma}_b dB(t') +{\sigma}_{gv} dB_{gv}(t')  \\
  +\int_0^t  \int_{\theta} \sigma_n(t',s,\theta)\tilde{N}(dt',d\theta) +\int_0^t  {\sigma}_{r} {\color{blue} dR^H(t')}, 
  \\    |\mathcal{S}| \geq  2
\end{array}  \right.
\end{equation}

\begin{table}[h!]   \label{label:tablenoise1} \caption{Noises examined in MFTGs.}
 \begin{tabular}{|  p{1.8cm} |  p{3cm} |  p{4cm} | }   
 \hline
 Noises  & Gaussian &  Non-Gaussian \\ \hline
 Markovian noise  & Brownian motion &  Poisson Jump, Regime Switching\\  \hline
 non-Markovian noise  & Fractional Brownian motion, Gauss-Volterra   &  {\color{blue} Rosenblatt noise}\\ \hline
 \end{tabular}
 \end{table}

\section{ Hierarchical Mean-Field-Type Games driven by Rosenblatt noise  }

 The game initiates at layer 1. Notably, the layer 1 decision-makers determine initial choices, influencing subsequent responses at layer 2. If multiple decision-makers exist at layer 1, they engage in a  strategic layer-game. First-layer decisions are observable by second-layer decision-makers, yet adjustments are unavailable once made.

Decision-makers at layer $h$ base choices on preceding layers 1 to $h-1$, interacting with any other decision-maker at the same layer. Subsequently, layer $h+1$ decision-makers respond to decisions from layers 1 to $h.$ While decisions up to the layer $h$ are observable by $(h+1)$-th layer decision-makers, adjustments by $h$-th layer decision-makers post-choice are restricted.

\begin{table}[htb]  \nonumber
 \begin{tabular}{| p{1.8cm} |  p{4cm}   |  p{2.5cm}  |}   
 \hline
MFTG & without Rosenblatt & with Rosenblatt \\ \hline
1 layer& explicit Nash  MFTG &   semi-explicit  \\  \hline
2 layers & Stackelberg MFTG   &   semi-explicit   \\ \hline
fully hierarchical &  fully hierarchical MFTG &  semi-explicit  \\ \hline
 \end{tabular}  \label{label:tablenoise3} \caption{Noise effects  examined in MFTGs.}
  \end{table}

\begin{table}[h!]   \nonumber
%
%\label{label:tablenoise4} 
 \begin{tabular}{|  p{2cm} |  p{2.5cm}   |  p{4.3cm}  |}   
 \hline
MFTG & only with Brownian motion & if a Rosenblatt is added to\\ \hline
Nash &   equilibrium &   no longer equilibrium  \\  \hline
Stackelberg &   equilibrium &   no longer equilibrium  \\  \hline
Adversarial  &   saddle point &   no longer saddle point  \\  \hline
Cooperation  &   social optimum &   no longer optimal  \\  \hline
 \end{tabular} \caption{Effects of Rosenblatt noises  in MFTG equilibria.}
  \end{table}

\section{How different is MFTG  from multipopulation coalitional MFG?}
In classical game theory, the structure of payoff functions influences the equilibrium system. Payoff functions that are linear in joint mixed strategies are known as Von Neumann expected utility. However, many applications consider non-Von Neumann utility functions, which are not necessarily linear in the measure of states or actions. Such games are particularly relevant in risk quantification, where key risk measures are nonlinear in the underlying probability measure. A similar issue arises in mean-field-type game (MFTG) theory, a class of games where the payoff function depends not only on state-action pairs but also on their distribution. MFTG can be finite number of decision-makers, infinite number of decision-makers or a combination of both within several subpopulations and atoms. MFTG equilibrium systems are non-standard and involve the so-called Master Adjoint Systems (MASS). Mean-field-type games   differ from mean-field games, multi-population mean-field games, and multi-population coalitional mean-field games. To illustrate these distinctions, we examine a power distribution network (DIPONET) with a finite number of renewable energy producers competing in a Cournot market. Given the high volatility of renewable energy sources, variance awareness is introduced, with stochastic fluctuations modeled using a Rosenblatt process.
This work elucidates fundamental differences between these models.

A class of MFTG is  designed to analyze the interactions among a finite number  of renewable energy producers. This class of MFTG incorporate the dynamics of price formation, demand, and supply, alongside the impact of individual production decisions. The price dynamics in this  MFTG is  a mean-reverting around demand-supply, and  are influenced by the Rosenblatt process, and the payoff for each producer is evaluated over a long-term horizon. The key feature of MFTG is its ability to capture non-linearities in the payoff functions along with a variance-awareness of the supply, making it highly relevant for scenarios where individual actions significantly affect the overall system.

In contrast, Multi-Population Coalitional MFG mimics of the same scenario consider multiple coalitions or subpopulations, each of infinite size. These models are particularly virtual agents, not true decision-makers, as  individual actions have a negligible impact on the overall mean-field on that coalition. The best-response problem for each decision-maker within a coalition is formulated in a linear-and-mean-field feedback form. Unlike MFTG, the payoff functions in Multi-Population Coalitional MFG are typically linear in the individual mean-field, and the equilibrium strategies are derived based on the collective behavior of the coalitions.

 Key Differences  are established. 
In this particular scenario, MFTG involves a finite number of true decision-makers, whereas Multi-Population Coalitional MFG assumes an infinite number of virtual agents within each coalition.
   In MFTG,  individual actions significantly influence the mean-field, capturing higher-order payoffs and real-world scenarios more accurately. In contrast, Multi-Population Coalitional MFG assumes that individual actions have a negligible impact on the coalition mean-field  and on the other coalitions' mean-field.  MFTG models are  non-linear in the measure of the state and control actions, making them more suitable for measure-non-linear payoff functions. Multi-Population Coalitional MFG, on the other hand, simplifies assumes  measure-linear payoffs, which may not capture the complexities of real-world scenarios as effectively.
 MFTG models are better suited for scenarios where the production cost and storage are high, and the supply needs to be accurately modeled. Multi-Population Coalitional MFG may not adequately represent such scenarios due to its simplifying assumptions.

\begin{table}[h!]
\centering
\begin{tabular}{|p{3cm}|p{6cm}|p{6cm}|}
\hline
\textbf{Features} & \textbf{MFTG} & \textbf{Multi-Population MFG} \\ \hline
$I$  & true decision-makers & $K=I$ subpopulations of infinite size \\  \hline
$p$  & market price (random variable) &  \\  \hline
\textbf{Best-Response} & 
$\sup_{u_i \in U_i}   \mathbb{E} (p u_i 
- c_i u_i $  $ - \frac{r_i}{2} u_i^2 - \frac{\bar{r}_i}{2} (\int \mu(dy) u_i(y, \mu) )^2)  $   & $\sup_{u_k \in U_k}   \mathbb{E}  ( p u_k 
- c_k u_k  $ $ - \frac{r_k}{2} u_k^2 - \frac{\bar{r}_k}{2} (\int m^*(dy) u^*_k(y, m^*) )^2)  $ 
  \\ \hline \textbf{Linearity} & non-linear in $\mu$ & linear in $\mu$ \\ \hline
  
  \textbf{Decomposition} & 
$\sup_{u_i \in U_i}    cov(p, u_i)-\frac{r_i}{2} var(u_i) + (\int y \mu(dy) -c_i ) (\int \mu(dy) u_i(y, \mu) )$ $
  - \frac{r_i+\bar{r}_i}{2} (\int \mu(dy) u_i(y, \mu) )^2  $   &   
  \\ \hline
   \textbf{Completion} & 
$\sup_{u_i \in U_i}   - \frac{r_i}{2}  \mathbb{E}(u_i-\int \mu(dy) u_i(y, \mu)-\frac{p-\int y \mu(dy)}{r_i})^2 +   \frac{var(p)}{2r_i} $ 
$- \frac{r_i+\bar{r}_i}{2}  (\int \mu(dy) u_i(y, \mu)  -\frac{\int y \mu(dy)-c_i}{r_i+\bar{r}_i})^2 +   \frac{(\int y \mu(dy)-c_i)^2}{2(r_i+\bar{r}_i)}  $   & $\sup_{u_k \in U_k}   \mathbb{E}  (- \frac{r_k}{2}  (u_k -\frac{(p-c_k)}{r_k})^2 +   \frac{(p-c_k)^2}{2r_k})$ $  - \frac{\bar{r}_k}{2} (\int m^*(dy) u^*_k(y, m^*) )^2)  $ 
  \\ \hline
\textbf{Mean-Field} & $(\int \mu(dy) u_i(y, \mu) ), (\int \mu(dy) u_i(y, \mu) )^2, $ $cov(p, u_i), var(u_i) $ & $(\int m^*(dy) u^*_k(y, m^*) )^2$ \\ \hline
\textbf{ Equilibrium Feedback Strategies} &   $\frac{p-\int y \mu(dy)}{r_i} + \frac{\int y \mu(dy)-c_i}{r_i+\bar{r}_i}$  & $\frac{(p-c_k)}{r_k}$ \\ \hline
\textbf{Equilibrium Payoff} &  $\frac{var(p)}{2r_i} + \frac{(\int y \mu(dy)-c_i)^2}{2(r_i+\bar{r}_i)}$ &  $\int  \mu(dy) \frac{(y-c_k)^2}{2r_k}
   - \frac{\bar{r}_k}{2} (\int m^*(dy) u^*_k(y, m^*) )^2) $  \\ \hline
   \textbf{Gap} & $\frac{var(p)}{2r_i} +  \frac{(\int y \mu(dy)-c_i)^2}{2(r_i+\bar{r}_i)}$ &  $  \frac{var(p)}{2r_k} +\frac{1}{2}  \frac{r_k-\bar{r}_k}{r_k^2}  (\int y \mu(dy)-c_k)^2$  \\ \hline
\textbf{Price of simplicity} &0  &  $ \frac{1}{2}[\frac{1}{(r_i+\bar{r}_i)}-  \frac{r_i-\bar{r}_i}{r_i^2}  ](\int y \mu(dy)-c_i)^2$  \\ \hline
\textbf{Price of simplicity} &0  &  unbounded (by scaling)  \\ \hline
\end{tabular}
\caption{Comparison of MFTG and Multi-Population MFG. By scaling the parameters $(p,c_i,r_i,\bar r_i)$ to be $\lambda(p,c_i,r_i,\bar r_i)$ with $\lambda >0$ the price of simplicity is unbounded. It turns out Multi-Population MFG, which consists to freeze the mean-field terms, is irrelevant, inefficient and suboptimal for this simple variance-awareness problem. }
\label{table:comparison}
\end{table}

There are two types of mean-field terms: individual mean-field and population mean-field: 
\begin{itemize}
\item Individual Mean-Field: This can refer to the mean field of an individual state, the mean field of an individual control action, or the joint distribution of the individual state-action pair.
\item Population Mean-Field: This includes  the mean field of the population state, the mean field of the control action of all decision-makers, or the joint distribution of the state-action pair of all decision-makers.
\end{itemize}

In this particular Multi-Population Coalitional MFG, each member is non-atomic, meaning their individual actions do not significantly impact the overall coalition mean-field. The MFG framework can of course be extended to include both atomic and non-atomic agents but still with linear payoff function with respect to own state mean-field.
Conversely, in MFTG, a decision-maker can be atomic, nonatomic, or a combination of both, and the payoff is not necessarily linear in the individual mean-field.

For interactions among renewable energy producers, we show that Multi-Population Coalitional MFG is not suitable. This is because the simplifying assumptions of non-atomic members and linear payoffs do not capture the variance-awareness and significant impacts of individual production decisions on the overall system. MFTG, with its ability to model non-linearities and individual impacts, provides a more accurate and effective framework for optimizing renewable energy production. Table \ref{table:comparison} displays some comparison of MFTG and Multi-Population MFG in the same scenario. By   scaling the parameters $(p,c_i,r_i,\bar r_i)$ be $\lambda(p,c_i,r_i,\bar r_i)$ with $\lambda >0$ the price of simplicity (which is the profit gap obtained after simplifying or freezing the mean-field terms) is unbounded. It turns out Multi-Population MFG, which consists to freeze the mean-field terms, is irrelevant, inefficient and suboptimal for this simple variance-awareness problem.  We show that the result extends to a dynamic setting with state driven by Rosenblatt process.

{\color{blue} In this particular scenario, there is no dynamic programming principle, no Pontryagin maximum principle, no Wiener chaos expansion available in the literature  to solve games with state dynamics driven by Rosenblatt noise. We therefore proceed with a direct method. The direct method offers a semi-explicit representation of the solution under the restriction of the class of strategies.}

\subsection{Non-Zero-Sum MFTG for Renewable  Energy Producers   } 
The transition to renewable energy sources has introduced new challenges and opportunities for energy producers, particularly for prosumers. In competitive markets, understanding the dynamics of price formation and the impact of individual production decisions is crucial for maximizing profitability. This section provides a basic model that integrates demand, supply, and Rosenblatt price fluctuations to evaluate the long-term average payoff for renewable energy producers. Let $ I \geq 2 $ be the number of producers, and denote by $ \mathcal{I} = \{1, \ldots, I\} $ the set of producers. We consider the  price dynamics given by:
$$
p(t) = p_0 + \int_0^t \frac{1}{\epsilon} \left[ a + D - S(t') - p(t') \right] dt' + {\color{blue} R^H(t)},
$$
where $ p(t) $ is the price at time $ t $, $ p_0 $ is the initial price,$ a $ is a constant, $ D $ represents the demand, $ S(t) = \sum_{i \in \mathcal{I}} u_i(t) $ is the total supply, with $ u_i(t) $ being the production of agent $ i $, $ R^H(t) $ is the Rosenblatt process, $ \epsilon > 0 $ is a speed adjustment parameter, ${r}_i, \bar{r}_i >0 ,\  i\in \mathcal{I},   H \in (\frac{1}{2}, 1)$. The instant payoff for producer $ i $ is given by:
$$
\pi_i(t,p,u) = p(t) u_i(t) - c_i u_i(t)- \frac{{r}_i }{2} u_i^2(t) - \frac{\bar{r}_i }{2} \bar{u}_i^2(t),
$$
where $ \bar{u}_i =  \mathbb{E}[u_i] $ is the expected production of agent $ i $.
To evaluate the long-run average payoff as the horizon goes to infinity, we consider the time-averaged payoff:

$$
\bar{\Pi}_i(p,u) = \liminf_{T \to \infty} \frac{1}{T}   \mathbb{E} \int_0^T \pi_i(t,p,u) \, dt
$$

Due to the memory dependence and non-Gaussianity, we consider a very specific class of strategies as an illustrative example.
 Given $(u_j, \ j\in \mathcal{I}\backslash \{i\})$ in linear-and-mean-field-type feedback form, the best-response problem of the decision-maker $i$  is 
\begin{equation}  \label{sub0pb}
\begin{array}{ll}
\sup_{u_i\in   \mathcal{U}_i}   \liminf_{T \rightarrow \infty} \frac{1}{T}  \mathbb{E}  [  \int_0^T \pi_i(t,p,u) \, dt\ | p(0) = p_0 ],\\
\mbox{subject to } \\
p(t)= p_0+\int_0^t \frac{1}{\epsilon} [ a+D- \sum_{j\in \mathcal{I}}u_j(t')  - p(t') ] dt'   {\color{blue}+   R^H(t)}, 
\end{array}
\end{equation} 
with the restricted control action set given by 
\begin{equation} \nonumber
\begin{array}{ll}
                   \mathcal{U}_i := \{(u_i(t))_{t \geq 0}: u_i(t) = \eta_i (p(t)-\bar{p}(t)) + \bar \eta_i \bar{p}(t)+ \rho_i  \text{ with } (\eta_i,\bar \eta_i,  \rho_i ) \in \mathbb{R}^3, \\
                  (1+ \sum_{j\in \mathcal{I}} \eta_j) >0,  (1+ \sum_{j\in \mathcal{I}} \bar\eta_j) >0\}

\end{array}
\end{equation}

In practice, we are not in favour of considering such  infinite horizon problems. The horizon  length $T$  for such interaction between producers can be relatively large but not infinite. Here we use this structure in order to get explicit representation of semi-explicit solution.  In this MFTG there is finite number of true decision-makers. The mean-field-type terms are  $\bar{p}, \bar{u}_i,$ $ var(u_i),$   and $  \bar{u}_i^{2}.$ The quantities-of-interest in the payoff functional $ \mathbb{E}  [\pi_{i} (t,x,u)]$ is clearly non-linear in the measure of the state and in the measure of control action.
\subsubsection{Two sub-problems   } In order to solve  Problem in (\ref{sub0pb}) we formulate  two sub-problems within this MFTG framework.   The best-response problem of a decision-maker is analyzed, considering scenarios with linear-and-mean-field-type feedback.  The first subproblem is as follows.
 Given $(\tilde{u}_j, \ j\in \mathcal{I}\backslash \{i\})$ in linear-and-mean-field-type feedback form, the best-response problem of the decision-maker $i$  is 
\begin{equation}  \label{sub1pb}
\begin{array}{ll}
\sup_{\tilde{u}_i\in \tilde{\mathcal{U}}_i }   \liminf_{T \rightarrow \infty} \frac{1}{T}   \mathbb{E}  [  \int_0^T \tilde{\pi}_i(t, \tilde{p}(t), \tilde{u}(t)) \, dt\ | \tilde{p}(0) = \tilde{p}_0 ],\\
\mbox{subject to } \\
\tilde p(t)= \tilde p_0+\int_0^t \frac{1}{\epsilon} [  -\sum_{j\in \mathcal{I}}\tilde{u}_j  - \tilde{p} ] dt'   {\color{blue}+   R^H(t)},  
\end{array}
\end{equation} 
with  $ \tilde p_0=   p_0- \bar p_0,$   $ \tilde{\pi}_i(t, \tilde{p}(t), \tilde{u}(t))=  \tilde{p}(t)  \tilde{u}_i(t) - \frac{{r}_i }{2} \tilde{u}^2_i(t) $, and  \begin{equation} \nonumber
\begin{array}{ll}
                   \tilde{\mathcal{U}}_i := \{(u_i(t))_{t \geq 0}: \tilde{u}_i(t) = \eta_i \tilde{p}(t)   \text{ with }  \eta_i \in \mathbb{R}, \\
                   \mathbb{E}  [ \tilde{p}(t) ]=0, \   (1+ \sum_{j\in \mathcal{I}} \eta_j) >0 \}

\end{array}
\end{equation}

Subproblem (\ref{sub1pb}) requires a stochastic calculus from Rosenblatt process to evaluate the long-run payoff.

%%%%%%%

The second subproblem is as follows.
Given $(\bar{u}_j, \ j\in \mathcal{I}\backslash \{i\})$ in linear-and-mean-field-type feedback form, the best-response problem of the decision-maker $i$  is 
\begin{equation} \label{sub2pb}
\begin{array}{ll}
\sup_{\bar{u}_i\in \bar{\mathcal{U}}_i }   \liminf_{T \rightarrow \infty} \frac{1}{T}    [  \int_0^T \bar{\pi}_i(t,\bar{p}(t), \bar{u}(t)) \, dt\ | \bar p(0) = \bar p_0 ],\\
\mbox{subject to } \\
\bar{p}(t)= \bar{p}_0+\int_0^t \frac{1}{\epsilon} [  a+D-\sum_{j\in \mathcal{I}}\bar{u}_j(t')  - \bar{p}(t') ] dt',  
\end{array}
\end{equation} 
with  $ \bar{\pi}_i(t,\bar{p}(t), \bar{u}(t))= \bar{p}(t)  \bar{u}_i(t) -  c_i\bar{u}_i(t) -  \frac{{r}_i +\bar{r}_i  }{2} \bar{u}_i^2(t),$ \begin{equation} \nonumber
\begin{array}{ll}
                   \bar{\mathcal{U}}_i := \{(\bar{u}_i(t))_{t \geq 0}: \bar{u}_i(t) = \bar\eta_i \bar{p}(t) +\rho_i  \text{ with }   (\bar{\eta}_i,\rho_i) \in \mathbb{R}^2, \\
                  (1+ \sum_{j\in \mathcal{I}} \bar\eta_j) >0 \}

\end{array}
\end{equation}
Subproblem (\ref{sub2pb}) is a standard deterministic differential game.
\subsubsection{ Combining the solutions of these two sub-problems   }

We observe that the state is orthogonally decomposed as $x(t)=\tilde{x}(t)+\bar{x}(t),$  the control action is decomposed as $u_i(t)=\tilde{u}_i(t)+\bar{u}_i(t),$  the  payoff functional is decomposed as 
$  \mathbb{E}  [ {\pi}_i(t, {p}(t), {u}(t))]= \mathbb{E}  [ \tilde{\pi}_{i}(t,\tilde p,\tilde u)+\bar{\pi}_{i}(t,\bar p,\bar u)].$ 

\subsubsection{Explicit solution of the MFTG}

 The MFTG  (\ref{sub1pb}) is reduced to 
 
 \begin{equation}  \label{sub1pbprime}
\begin{array}{ll}
\sup_{\eta_i\in \tilde{\mathcal{G}}_i }   \liminf_{T \rightarrow \infty} \frac{1}{T}   \mathbb{E}  [  \int_0^T  (\eta_i -\frac{{r}_i }{2}\eta_i^2) \tilde{p}^2(t)\, dt\ | \tilde{p}(0) = \tilde{p}_0 ],\\
\mbox{subject to } \\
\tilde p(t)= \tilde p_0+\int_0^t \frac{1}{\epsilon} [ -(1+ \sum_{j\in \mathcal{I}} \eta_j) \tilde{p} ] dt'   {\color{blue}+   R^H(t)},  
\end{array}
\end{equation}

Then the MFTG  (\ref{sub1pb}) is semi-explicitly solved as  the following constrained  game:

\begin{thm}  The payoff of producer $i$ is explicitly computed as 
$\tilde{\Pi}_{i,\infty} (\eta) = \frac{\Gamma(2H+1)}{2[ \frac{1}{\epsilon}(1+ \sum_{j\in \mathcal{I}} \eta_j)]^{2H}} \left(\eta_i -\frac{{r}_i }{2}\eta_i^2\right)$ 
and the gain in the strategies spaces become the coupled space   $\tilde{\mathcal{G}}_i =\{\eta_i \ | \  1+ \sum_{j\in \mathcal{I}} \eta_j >0\}.$
  \end{thm}

The first order derivative of $\tilde{\Pi}_{i,\infty} (\eta)$ at $\eta_i$ yields 
  \begin{equation}\nonumber \begin{array}{ll}
            i\in \mathcal{I},\\ 
\partial_{\eta_i}\tilde{\Pi}_{i,\infty} (\eta) =  
\frac{\Gamma(2H+1)}{2[ \frac{1}{\epsilon}(1+ \sum_{j\in \mathcal{I}} \eta_j)]^{2H+1}} \left(1-r_i\eta_i)(1+ \sum_{j\in \mathcal{I}} \eta_j) -2H \eta_i+H r_i \eta_i^2\right)\\
=  \frac{\Gamma(2H+1)}{2[ \frac{1}{\epsilon}(a+ \eta_i)]^{2H+1}} \left((1-r_i\eta_i)(a+ \eta_i) -2H \eta_i+H r_i \eta_i^2\right)\\
= \frac{\Gamma(2H+1)}{2[ \frac{1}{\epsilon}(a+ \eta_i)]^{2H+1}} 
\left(   (H-1) r_i  \eta_i^2  -(r_ia+ 2H-1) \eta_i +a \right), \\
 a= (1+ \sum_{j\in \mathcal{I}\backslash \{i\}} \eta_j)
 \end{array}
\end{equation}
at the positive critical point $\partial_{\eta_i^2}\tilde{\Pi}_{i,\infty} (\eta)  <0$ making a  concavity in $\eta_i.$ Therefore we have the following best response gain: 

   \begin{thm} Given the admissible strategies of the other agents  $(u_j\in \mathcal{U}_j, \  j\neq i),$
           the best-response gain of producer $i$ 
  in the family of admissible feedback gains is given by
$
{\eta}_i^* =   \frac{- (r_ia+ 2H-1) + \sqrt{ (r_ia+ 2H-1)^2 +4(1-H) r_i a} }{2(1-H) r_i  }
$
and the best-response cost is given by
$\tilde{\Pi}_{i,\infty} (\eta^*)   $ 
  \end{thm}

            \begin{thm} If  there  exists a solution to the following system:  
            \begin{equation}\nonumber \begin{array}{ll}
            i\in \mathcal{I},\\ 
       \eta_i = \frac{- (r_i(1+ \sum_{j\in \mathcal{I}\backslash \{i\}} \eta_j)+ 2H-1) + \sqrt{ (r_i(1+ \sum_{j\in \mathcal{I}\backslash \{i\}} \eta_j)+ 2H-1)^2 +4(1-H) r_i (1+ \sum_{j\in \mathcal{I}\backslash \{i\}} \eta_j)} }{2(1-H) r_i  }
       \\
       1+ \sum_{j\in \mathcal{I}} \eta_j >0,
         \end{array}
\end{equation}
 then, the non-zero-sum stochastic differential game under Rosenblatt noise (\ref{sub1pb}) has an equilibrium in stationary linear state-feedback strategies by given ${\eta}_i^* (p(t)-\bar p(t)).$
  \end{thm}

  The MFTG  (\ref{sub2pb}) is reduced to

 \begin{equation}  \label{sub2pbprime}
\begin{array}{ll}
\sup_{\bar\eta_i\in \bar{\mathcal{G}}_i }   \liminf_{T \rightarrow \infty} \frac{1}{T}   \mathbb{E}  [  \int_0^T  -  \frac{{r}_i +\bar{r}_i  }{2} (\bar{u}_i -\frac{(\bar{p}(t) -  c_i)}{({r}_i +\bar{r}_i)} ) + \frac{(\bar{p}(t) -  c_i)^2}{2({r}_i +\bar{r}_i)}\, dt\ | 
\mbox{subject to } \\
\bar{p}(t)= \bar{p}_0+\int_0^t \frac{1}{\epsilon} [  a+D-\sum_{i\in \mathcal{I}}\rho_i - (1+\sum_{j\in \mathcal{I}}\eta_j) \bar{p}(t') ] dt',  
\end{array}
\end{equation} 
$  \bar{\mathcal{G}}_i := \{(\bar{\eta}_i,\rho_i)\in \mathbb{R}^2:   (1+ \sum_{j\in \mathcal{I}} \bar\eta_j) >0 \} $

By direct computation, we obtain  the following:
 \begin{thm} 
The non-zero-sum deterministic differential game with long-run payoff (\ref{sub2pb}) has an equilibrium in  linear state feedback strategies by given $\bar{\eta}_i^*\bar p(t)+\rho_i$ where 
 \begin{equation}  \label{sub2pbprime2}
\begin{array}{ll}
\lim_{t\to +\infty} \bar{p}(t)  = \bar{p}^*=\frac{a+D+\sum_{i\in \mathcal{I}} \frac{c_i}{{r}_i +\bar{r}_i}}{1+\sum_{j\in \mathcal{I}} \frac{1}{{r}_i +\bar{r}_i}},\\
\bar\eta^*_i= \frac{1}{{r}_i +\bar{r}_i}\\ 
\rho^*_i=\frac{-c_i}{{r}_i +\bar{r}_i}\\
\bar{u}_i^*(t) = \frac{(\bar{p}(t) -  c_i)}{({r}_i +\bar{r}_i)} ,\\
\bar{\Pi}^*_{i,\infty}(\bar\eta, \rho)=  \frac{(\bar p^*-c_i)^2}{2({r}_i +\bar{r}_i)}.
\end{array}
\end{equation}  
  \end{thm}
  We now combine the two subproblems to solve the original problem in (\ref{sub0pb}):
  
   \begin{thm} 
   Consider  (\ref{sub0pb}). The price-taking Cournot-Nash equilibrium is given by
   \begin{equation}  \label{sub2pbprime2t}
\begin{array}{ll}
u^*_i(t) = \eta^*_i (p^*(t)-\bar{p}^*(t)) + \bar \eta_i \bar{p}^*(t)+ \rho^*_i ,\\
\eta^*_i = \frac{- (r_i(1+ \sum_{j\in \mathcal{I}\backslash \{i\}} \eta^*_j)+ 2H-1) + \sqrt{ (r_i(1+ \sum_{j\in \mathcal{I}\backslash \{i\}} \eta^*_j)+ 
2H-1)^2 +4(1-H) r_i (1+ \sum_{j\in \mathcal{I}\backslash \{i\}} \eta^*_j)} }{2(1-H) r_i  }
 \\     
\lim_{t\to +\infty} \bar{p}(t)  = \bar{p}^*=\frac{a+D+\sum_{i\in \mathcal{I}} \frac{c_i}{{r}_i +\bar{r}_i}}{1+\sum_{j\in \mathcal{I}} \frac{1}{{r}_i +\bar{r}_i}},\\
\bar\eta^*_i= \frac{1}{{r}_i +\bar{r}_i}\\ 
\rho^*_i=\frac{-c_i}{{r}_i +\bar{r}_i}\\
\bar{u}_i^*(t) = \frac{(\bar{p}^*(t) -  c_i)}{({r}_i +\bar{r}_i)} ,\\
{\Pi}^*_{i,\infty}(\eta, \bar\eta, \rho)=   \frac{\Gamma(2H+1)}{2[ \frac{1}{\epsilon}(1+ \sum_{j\in \mathcal{I}} \eta^*_j)]^{2H}} \left(\eta^*_i -\frac{{r}_i }{2}(\eta^*_i)^2\right)+ \frac{(\bar p^*-c_i)^2}{2({r}_i +\bar{r}_i)}.
\end{array}
\end{equation}  
  \end{thm}

We now compute explicitly the limiting price as function of associated control strategies.
 \begin{thm} 
The non-zero-sum deterministic differential game with long-run payoff (\ref{sub2pb}) has an equilibrium in  stationary strategies by given $\bar{\eta}_i^*\bar p^*+\rho^*_i$ where 
 \begin{equation}  \label{sub2pbprime3}
\begin{array}{ll}
\bar{p}=\frac{a+D-\sum_{i\in \mathcal{I}} \rho_i }{1+\sum_{j\in \mathcal{I}}\bar  \eta_j},\\
\bar{\Pi}^*_{i,\infty}(\bar\eta, \rho)= (\bar \eta_i-\frac{{r}_i +\bar{r}_i}{2}\eta_i^2) (\bar{p}^*)^2 + (\rho_i -c_i\bar \eta_i-(r_i+\bar{r}_i)\bar \eta_i\rho_i) (\bar{p}^*)  -(c_i\rho_i+\frac{{r}_i+\bar{r}_i }{2}\rho_i^2)\\ 
=   (\bar \eta_i-\frac{{r}_i+\bar{r}_i }{2}\bar \eta_i^2) (\frac{a+D-\sum_{i\in \mathcal{I}} \rho_i }{1+\sum_{j\in \mathcal{I}} \bar \eta_j})^2 + (\rho_i -c_i\bar \eta_i-(r_i+\bar{r}_i)\bar \eta_i\rho_i) (\frac{a+D-\sum_{i\in \mathcal{I}} \rho_i }{1+\sum_{j\in \mathcal{I}} \bar \eta_j})  -(c_i\rho_i+\frac{{r}_i+\bar{r}_i }{2}\rho_i^2) \\
(\bar\eta^*_i, \rho^*_i)\in \arg\max_{\bar\eta_i, \rho_i} \bar{\Pi}^*_{i,\infty}(\bar\eta_i, \rho_i, \bar\eta^*_{-i}, \rho^*_{-i})\\
\bar{u}_i^* =\bar{\eta}_i^*\bar p^*+\rho^*_i \\
\bar{p}^*=\frac{a+D-\sum_{i\in \mathcal{I}} \rho^*_i }{1+\sum_{j\in \mathcal{I}}\bar  \eta^*_j}.
\end{array}
\end{equation}  
  \end{thm}
 \begin{thm} 
   Consider  (\ref{sub0pb}). The stationary closed-loop Cournot-Nash equilibrium is given by
   \begin{equation}  \label{sub2pbprime2w}
\begin{array}{ll}
u^*_i = \eta^*_i (p^*-\bar{p}^*) + \bar \eta^*_i \bar{p}^*+ \rho^*_i ,\\
\eta^*_i = \frac{- (r_i(1+ \sum_{j\in \mathcal{I}\backslash \{i\}} \eta^*_j)+ 2H-1) + \sqrt{ (r_i(1+ \sum_{j\in \mathcal{I}\backslash \{i\}} \eta^*_j)+ 
2H-1)^2 +4(1-H) r_i (1+ \sum_{j\in \mathcal{I}\backslash \{i\}} \eta^*_j)} }{2(1-H) r_i  }
 \\     
(\bar\eta^*_i, \rho^*_i)\in \arg\max_{\bar\eta_i, \rho_i} \bar{\Pi}^*_{i,\infty}(\bar\eta_i, \rho_i, \bar\eta^*_{-i}, \rho^*_{-i})\\
\bar{u}_i^* =\bar{\eta}_i^*\bar p^*+\rho^*_i \\
\bar{p}^*=\frac{a+D-\sum_{i\in \mathcal{I}} \rho^*_i }{1+\sum_{j\in \mathcal{I}}\bar  \eta^*_j}.
  
\end{array}
\end{equation}  
  \end{thm}
 
\subsection{Multi-Population Coalitional MFG   } Consider $K\geq 2$ number of coalitions (or subpopulation) and denote by  $\mathcal{K}=\{1,2,\ldots, K\}$ the set of all coalitions.  Each coalition is assumed to be of infinite size with a proper scaling. The multi-population coalitional MFG driven by  Rosenblatt noise can be defined. 
Consider a generic decision-maker from coalition $k.$ Given $u^*_k$  of all other members of the coalition $k$ and  $(u^*_j, \ j\in \mathcal{K}\backslash \{k\})$ of the other coalitions in linear-and-mean-field  feedback form, the best-response problem of the generic decision-maker  in coalition $k$   in linear-and-mean-field feedback form is 
\begin{equation}
\begin{array}{ll}
\sup_{u_k\in   \mathcal{U}_k}   \liminf_{T \rightarrow \infty} \frac{1}{T}  \mathbb{E}  [  \int_0^T p(t) u_k(t) -  c_k u_k(t)- \frac{{r}_k }{2} u_k^2(t)  {\color{blue} - \frac{\bar{r}_k }{2}\int m^*(t,dy) {u}^*_k(t,y,m^*)  }\, dt\ | p(0) = p_0 ],\\
\mbox{subject to } \\
p(t)= p_0+\int_0^t \frac{1}{\epsilon} [ a+D- \sum_{j\in \mathcal{K}}u_j(t')  - p(t') ] dt'   {\color{blue}+   R^H(t)}, \\
\bar{m}^*(t)=  \bar{m}_0
 +\int_0^t \frac{1}{\epsilon} [ a+D- \sum_{j\in \mathcal{K}}\int m^*(t',dy) {u}^*_j(t,y,m^*)  - \bar{m}^*(t') ] dt' \\
 m^*(t,dy) = \mathbb{P}_{p^*(t)}(dy)
\end{array}
\end{equation} 
with the restricted control action set given by 
\begin{equation} \nonumber
\begin{array}{ll}
                   \mathcal{U}_k := \{(u_k(t))_{t \geq 0}: u_k(t) = \eta_k p(t)+ \rho_k  \text{ with } (\eta_k, \rho_k ) \in \mathbb{R}^2, \\
                  (1+ \sum_{j\in \mathcal{K}} \eta_j) >0\} 

\end{array}
\end{equation}

By definition of multi infinite population MFG, the control action  $u_k$ of a generic decision-maker in coalition $c$ has no impact on mean of the population mean-field $\bar{m}^*(t).$ It has also no effect on the population mean-field $m^*(t,dy),$ the subpopulation mean of mean-field of actions  $\int m^*(t,dy) {u}^*_k(t,y,m^*),$  and any other sub-population mean of  mean-field of actions $\int m^*(t,dy) {u}^*_j(t,y,m^*),  j\in \mathcal{K}\backslash \{k\}.$ Next we rewrite the functionals involved.

The best response is decoupled of $m^*$ and it is given by 
\begin{equation}
\begin{array}{ll}
\sup_{( \eta_k, \rho_k)\in   \mathcal{G}_{k,mfg}} 
  (\eta_k-\frac{{r}_k }{2}\eta_k^2) [ \frac{\Gamma(2H+1)}{2[ \frac{1}{\epsilon}(1+ \sum_{j\in \mathcal{I}} \eta_j)]^{2H}} +\frac{(a+D)^2}{\epsilon^2} +2\frac{a+D}{\epsilon}  (\frac{a+D-\sum_{i\in \mathcal{I}} \rho_i }{1+\sum_{j\in \mathcal{I}} \eta_j}) ] \\
  + (\rho_k -c_k\eta_k-r_k\eta_k\rho_k) (\frac{a+D-\sum_{i\in \mathcal{I}} \rho_i }{1+\sum_{j\in \mathcal{I}} \eta_j})  -(c_k\rho_k+\frac{{r}_k }{2}\rho_k^2)   {\color{blue} - \frac{\bar{r}_k }{2}  (\eta_k^*\frac{a+D-\sum_{i\in \mathcal{I}} \rho^*_i }{1+\sum_{j\in \mathcal{I}} \eta^*_j}+\rho_k^*)^2},\\
\mbox{subject to } \\
\bar{m}^*(t)=  \bar{m}_0
 +\int_0^t \frac{1}{\epsilon} [ a+D- \sum_{k\in \mathcal{K}}\rho^*_k  - (1+ \sum_{j\in \mathcal{K}}\eta^*_j) \bar{m}^*(t') ] dt' \\
 m^*(t,dy) = \mathbb{P}_{p^*(t)}(dy)
\end{array}
\end{equation} 

\begin{thm}
The mean-field equilibrium of virtual multi-population coalitional MFG does not depend on $\bar{r}_k.$ In particular, for very large $r_k$ the equilibrium stays the same. 
\end{thm}

Note that the term $  {\color{blue} - \frac{\bar{r}_k }{2}  (\eta_k^*\frac{a+D-\sum_{i\in \mathcal{I}} \rho^*_i }{1+\sum_{j\in \mathcal{I}} \eta^*_j}+\rho_k^*)^2}$ is from frozen as it is coming the infinite population mean-field. The consequence of this observation is that the mean-field equilibrium does not capture the energy producer market as when the production cost and storage are extremely high, the supply should be nearly zero as observed in the true MFTG above. Thus, the virtual game created to represent the MFTG into a multipopulation MFG does not do the job. In this particular case, it  does not even simplify the analysis. This multipopulation MFG does not have the orthogonal decomposition avantage used in MFTG. This proves also that multi-population coalitional MFG differs from MFTG.

\subsection{Difference between MFTG  and  multipopulation coalitional MFG}

\subsubsection{True number of decision-makers}
In this particular scenario, MFTG involves a finite number of true decision-makers, whereas Multi-Population Coalitional MFG assumes an infinite number of virtual agents within each coalition to represent to same scenario. As we can see these virtual agents are not true decision-makers. 

\subsubsection{No asymptotic indistinguishability assumption }
In this MFTG, each decision-maker has a strong influence on its own mean-field as well as in the other mean-field terms. The asymptotic indistinguishability assumption does not hold in MFTG. The permutation of index of the decision-maker changes the payoff. However, in the Multi-Population Coalitional MFG off the same scenario, it is the asymptotic indistinguishability is assumed within each coalition.

\subsubsection{Effect of a generic decision-maker on its individual  mean-field  }
In MFTG,  individual actions significantly influence the mean-field, capturing higher-order payoffs and real-world scenarios more accurately. In contrast, this Multi-Population Coalitional MFG assumes that individual actions have a negligible impact on the coalition mean-field  and on the other coalitions' mean-field. 

\subsubsection{The MFTG payoff is non-linear in the individual state mean-field  }
 MFTG models are  nonlinear in the measure of the state and control actions, making them more suitable for measure-non-linear payoff functions. Multi-Population Coalitional MFG, on the other hand, simplifies assumes  measure-linear payoffs, which may not capture the risk-awareness of real-world scenarios as effectively.

\subsubsection{The MFTG payoff captures higher-order payoffs }
The MFTG payoff captures higher-order payoffs.  Multi-Population Coalitional MFG is limited to measure-linear payoff. The equilibrium equation of MFGs are non-valid when it comes to non-linear payoff such  as variance, quantile, expected-value-at-risk, expectile value-at-risk, extremile value-at-risk, and entropic value-at-risk. 

\subsubsection{MFTGs capture more real-world scenarios  }
 MFTG models are better suited for scenarios where the production cost and storage are high, and the supply needs to be accurately modeled. Multi-Population Coalitional MFG may not adequately represent such scenarios due to its simplifying assumptions.

Both MFTG and Multi-Population Coalitional MFG offer valuable insights in the demand side. When it comes to the dynamics of renewable energy production,  MFTG stands out  as the more suitable framework for measure-non-linear payoff functions. MFTG ability to capture the variance-awareness of individual actions and their impact on the overall system makes it a powerful tool for optimizing renewable energy markets.  Multi-Population Coalitional MFG approach freezes all the mean-field terms and use a matching argument. This leads to suboptimality in important scenarios such as risk-awareness renewable energy production.

\section{ Rosenblatt Foundational Diffusion Models}
We introduce Rosenblatt Foundational Diffusion Models. 
By moving from Brownian motion to fractional Brownian motion and finally to the Rosenblatt process,  we illustrate a progression towards more accurate and robust models that better capture the real-world data. 

\subsection{AI: Diffusion Transformers (under Brownian motion)}
We examine diffusion transformers under the assumption of Brownian motion. It covers various types of data, including time-series, image, audio, and video data. Diffusion transformers under Brownian motion face significant limitations due to their inability to capture long-range dependencies and higher-order statistical properties inherent in many real-world datasets. Brownian motion  which is characterized by short-range dependence and a lack of memory, lead to models that often fail to account for the  non-Gaussian behaviors observed in practice.

    Consider Ornstein-Uhlenbeck process  $x(t) = x_0+\int_0^t \theta (m - x(t')) dt' + \int_0^t \sigma d\textcolor{red}{B(t')}$ with constant  $x_0, \theta, m, \sigma.$ Then $$x(t)= e^{-\theta t} x_0  + (1-e^{-\theta t})m + \frac{\sigma}{\sqrt{2\theta}}B_{1-e^{-2\theta t}}.$$
     It means that, starting from any constant $x_0$ we end up at time $T$ with $ \mathcal{N}(e^{-\theta T} x_0  + (1-e^{-\theta T})m,  \frac{\sigma^2}{{2\theta}}(1-e^{-2\theta T}).$ By properly choosing $\sigma, \theta, m$ we can match with any Gaussian $ \mathcal{N}(m_T, \sigma^2_T)$ at time $T.$ To do so, we match $e^{-\theta T} x_0  + (1-e^{-\theta T})m =m_T$ i.e.,  $m=\frac{(m_T-e^{-\theta T} x_0)}{ (1-e^{-\theta T})}.$ Similarly we solve $\frac{\sigma^2}{2\theta}(1-e^{-2\theta T})=\sigma_T^2$ to freely match with 
     $ \sigma= \sqrt{\frac{2\theta\sigma_T^2}{ (1-e^{-2\theta T})}}.$  It follows that for any $x_0$ and $\theta>0$ the stochastic process driven by Brownian motion      
 $$x(t) = x_0+\int_0^t \theta (\frac{(m_T-e^{-\theta T} x_0)}{ (1-e^{-\theta T})} - x(t')) dt' + \int_0^t  \sqrt{\frac{2\theta\sigma_T^2}{ (1-e^{-2\theta T})}} d\textcolor{red}{B(t')}$$ will have the distribution of the Gaussian $ \mathcal{N}(m_T, \sigma^2_T)$ at time $T.$ The parameter $\theta$ represents the speed rate at which $x_0$ is being transferred to $m_T.$
 In this perspective one may want to  simply consider  the process $ y(t)=(1-\frac{t}{T})x_0 + \frac{t}{T} \mathcal{N}(m_T, \sigma^2_T)$ which also transfer initial signal $x_0$ at time $0$ to  the process $\mathcal{N}(m_T, \sigma^2_T)$ at time $T.$ 
 
One of the   key differences between these two solutions is in  their path. The first one has a speed $\theta (m - x(t))$ while the second one has a constant speed $\frac{(m_T-x_0)}{T}$ and of course the variance are also different which can make a difference when it comes to risk-aware decision-making between the curves. Optimal transport theory deals with such transport map and path cost criteria.

When $x_0$ is an input signal of a transformer, which could be  a \textbf{Time-Series}, \textbf{Text}, \textbf{Image},  \textbf{Audio},  \textbf{Video} etc, we can still adjust the above methodology to be vectors or matrices. This leads to 
 \begin{itemize}
     \item    \textbf{Time-Series} Diffusion Transformers under Brownian motion
        \item    \textbf{Image} Diffusion Transformers under Brownian motion
          \item    \textbf{Audio} Diffusion Transformers under  Brownian motion
          \item    \textbf{Video} Diffusion Transformers under  Brownian motion
    \end{itemize}
    The final output of the diffusion with distribution $\mathcal{N}(m_T, \sigma^2_T)$ can be seen as a pure mask. We now reverse the process to find the original signal. Given the terminal process $x_T=x(T) $ and using the martingale property, we aim to identity the reverse process  $\hat{x}(t)=x(T-t).$ It starts at time $t=0$ with a mask $\mathcal{N}(m_T, \sigma^2_T)$  and ends up with a high quality signal output $x_0.$
    With that in mind, one can identify high quality text, time-series, image, audio, video as final output and then train the neural network to arrive at such a top quality outcomes at time $T$.  We observe that the process  $ \hat{x}(t)=(1-\frac{t}{T})\mathcal{N}(m_T, \sigma^2_T) + \frac{t}{T} x_0$ starting a  mask $\mathcal{N}(m_T, \sigma^2_T)$ at  time $0$ and discover progressively   the signal $x_0$ at time $T.$ The time-reversed OU process 
    $${x}(t) = x_T-\int_t^T \theta (m - x(t')) dt' - \int_t^T \sigma d\textcolor{red}{B(t')}.$$
     The modified drift of $\hat{x}(t)$ is  $-\theta (m - \hat{x}(t)) + \sigma^2 \frac{(\hat{x}(t)-\mu(T-t))}{v^2(T-t)}.$ The  mean value  is 
   $\mu(T-t)=e^{-\theta (T-t)} x_0  + (1-e^{-\theta (T-t)})m$ and the variance is  
$v^2(T-t)=\frac{\sigma^2}{2\theta}(1-e^{-2\theta (T-t)}).$
  The forward drift starting from $x_T$ is given by 
  $-\theta (m - \hat{x}(t)) + 2\theta \frac{(\hat{x}(t)-e^{-\theta (T-t)} x_0  - (1-e^{-\theta (T-t)})m)}{(1-e^{-2\theta (T-t)})}.$
   The time-reverse process starting with a mask $mask \sim \mathcal{N}(m_T, \sigma^2_T)$ and solves
  $$\hat{x}(t) = mask+\int_0^t -\theta (m - \hat{x}(t')) + 2\theta \frac{(\hat{x}(t')-e^{-\theta (T-t')} x_0  - (1-e^{-\theta (T-t')})m)}{(1-e^{-2\theta (T-t')})} dt' +\int_t^T \sigma d\textcolor{red}{B(t')}.$$ The Gaussian distribution will be approaching the Dirac measure concentrated at $x_0$ as $t$ goes to $T.$ 
 Here, \textbf{Ignoring longer-range memory may lead to  flawed generative machine intelligence  models as the output is a pure Gaussian independently of the input.}

\subsection{AI: Sub/Super Diffusion Transformers under Fractional Brownian motion}
We introduce sub and super diffusion transformers that operate under fractional Brownian motion. It emphasizes that considering longer-range memory leads to better generative machine intelligence  models. We  cover the same types of data as the previous one (time-series, image, audio, and video), but with an improved approach that accounts for the memory effects inherent in fractional Brownian motion. We introduce a  mean-field-type fractional diffusion generative model. The mean-field type terms are mainly involved in the backward process where the denoising pass requires the mean and the variance to match the original process. We do not need the entire mean field state distribution in this particular case, we only need the first two moments.
\subsubsection*{ Fractional Forward Dynamics }
\[
dx(t) = \theta(t)(\bar{m}(t) - x(t))\,dt + \sigma(t)\,dB^H(t), \quad x(0) = x_0,
\]

The solution \( x(t) \) admits the explicit form
\[
x(t) = e^{-\Phi(t)} x_0 + e^{-\Phi(t)} \int_0^t e^{\Phi(t')} \theta(t') \bar{m}(t')\,dt' + e^{-\Phi(t)} \int_0^t e^{\Phi(t')} \sigma(t')\,dB^H(t').
\]
where $\Phi(t) := \int_0^t \theta(t')\,dt'.$  In this deterministic coefficient setting and deterministic initial signal $x_0,$ the random variable $x(t)$  has the same distribution as the Gaussian $\mathcal{N}(m(t), v^2(t)).$
where 
\begin{equation}
    \begin{array}{l}      
m(t) = e^{-\Phi(t)} x_0 + e^{-\Phi(t)} \int_0^t e^{\Phi(t')} \theta(t') \bar{m}(t')\,dt',\\
v^2(t) = e^{-2\Phi(t)} H(2H-1) \int_0^t \int_0^t e^{\Phi(t')} \sigma(t') e^{\Phi(s')} \sigma(s') |t' - s'|^{2H - 2} \,dt' \,ds',
    \end{array}
\end{equation}

\[
dx(t) = \theta(t) \left( \frac{m_T^* - x_0 e^{-\Phi(T)}}{1 - e^{-\Phi(T)}} - x(t) \right) dt + \sigma_T^*  T^{-H} e^{\Phi(T) - \Phi(t)} dB^H(t), \quad x(0) = x_0,
\]
where satisfies $x(T) \sim \mathcal{N}(m_T^*, (\sigma_T^*)^2). $

Given the forward process is Gaussian, the reverse-time dynamics are governed by:
\begin{equation}
\begin{array}{l}
dx(t) = 
\left[ \theta \left( \frac{m^*_T - e^{-\Phi(T)} x_0}{1 - e^{-\Phi(T)}} - x(t) \right) + \frac{x(t) - m(t)}{t} 2H  \right] dt \\ 
+ (\sigma_T^*)  T^{-H} e^{\Phi(T) - \Phi(t)}\, d\bar{B}^H(t),
\\
\\
m(t) =
e^{-\Phi(t)} x_0 +(1-e^{-\Phi(t)}) (\frac{m^*_T - e^{-\Phi(T)} x_0}{1 - e^{-\Phi(T)}}) ,
\\
v^2(t) =  e^{2\Phi(T)-2\Phi(t)} (\frac{t}{T})^{2H} (\sigma_T^*)^2,

\end{array}
\end{equation}
where $\bar{B}^H(t)$ is a time-reversed fractional Brownian motion, and $\tilde{\sigma}(t)$ adjusts for non-Markovian memory in the reverse direction.
Because this process remains Gaussian, the exact score function is known analytically:
\begin{equation}
\nabla_y \log p_{x(T-t)}(y) = -\frac{y - m(T-t)}{v^2(T-t)}.
\end{equation}
\subsection{Super-Diffusion Transformers under Rosenblatt process}
We introduce super-diffusion transformers driven by the Rosenblatt process. This approach incorporates longer-term context-awareness and non-Gaussianity which significantly enhances the performance of generative machine intelligence models. The types of data covered include time-series, image, audio, and video. The Rosenblatt process is particularly suited for capturing the  non-Gaussian dependencies observed in real-world data.
$x(t) = x_0+\int_0^t \theta (m - x(t')) dt' + \int_0^t \sigma  d\textcolor{red}{R^H(t')}.$ 
 Then $$x(t)= e^{-\theta t} x_0  + (1-e^{-\theta t})m + \sigma \int_0^t e^{-\theta (t-t')}d\textcolor{red}{R^H(t')}.$$ As $H$ approaches one, this last process has the same distribution as
 $$x(t)= e^{-\theta t} x_0  + (1-e^{-\theta t})m + \frac{\sigma }{\theta \sqrt{2}} (1- e^{-\theta t}) (Z^2-1)$$ where $ (Z^2-1)$ is the centred chi-square random variable, and $Z\sim \mathcal{N}(0,1).$ The variance of the process is
 $$var\left(\int_0^t e^{\theta t'}d\textcolor{red}{R^H(t')}\right)= H(2H-1)\int_0^t  \int_0^t e^{\theta t'} e^{\theta t''}  |t'-t''|^{2H-2} dt'dt'',$$
 which highlights the strong long-run dependence and non-Gaussianity. 

\section*{Acknowledgments}
Authors gratefully acknowledge support from TIMADIE, Guinaga, Grabal, LnG Lab for the MFTG for Machine Intelligence project.
Authors gratefully acknowledge support from U.S. Air Force Office of Scientific Research under grants number FA9550-17-1-0259

\section{Conclusion and Future Work}
In this work, we introduced a  system framework that incorporates Rosenblatt noise, a non-Gaussian, self-similar, and long-range dependent process, into mean-field-type game theory (MFTG). The framework has a wide-reaching implications across systems engineering, control, and strategic decision-making among multiple AI agents. Our analysis revealed that classical assumptions rooted in Gaussian or Markovian processes often underestimate tail risks and mischaracterize volatility in real-world domains such as power grids, agricultural cycles, e-commerce behavior, and underwater communications. By extending stochastic calculus to Rosenblatt dynamics and embedding this within the MFTG, we derived new equilibrium concepts, saddle-point strategies, and variance-aware optimal control laws that respect higher-order dependencies. These results are particularly relevant for agentic AI, intelligent agents capable of long-term, autonomous decision-making, where the robustness and foresight required to manage uncertainty must account for heavy-tailed phenomena and multimodal interactions. Our findings demonstrate that risk quantification in multi-agent systems is not merely a statistical detail but a foundational aspect of modeling real-world scenarios. In future work, we aim to expand this framework to mean-field-type  learning scenarios, agentic memory systems under asymmetric information, incorporating tools from fractional calculus, Wasserstein geometry, and holonormalization to further enhance risk-aware learning in  distributed and decentralized human-machine co-intelligence systems. 
\bibliographystyle{IEEEtran}

\bibliography{ros_bib}

\mbox{~}
  \begin{IEEEbiography}[{\includegraphics[width=1in,clip,keepaspectratio]{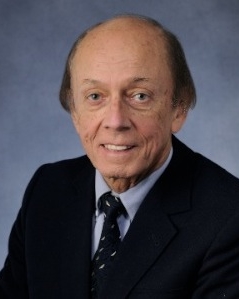}}]
  {Tyrone E. Duncan} (M'92-SM'96-F'99) received his bachelor’s degree from Rensselaer Polytechnic Institute in 1963 and his master’s and doctoral degrees from Stanford University in 1964 and 1967, respectively. He has held regular positions at the University of Michigan (1967-1971), the State University of New York, Stony Brook (1971-1974), and the University of Kansas (1974-present) where he is Professor of Mathematics. He has held visiting positions at the University of California, Berkeley (1969-1970),  the University of Bonn (1978-1979), and Harvard University  (1979-1980). He is the author or co-author of more than 125 papers in the areas of stochastic optimal control, adaptive control, filtering and communication as well as related topics in probability. He was on editorial boards on SIAM Journal on Control and Optimization and IEEE Transactions on Automatic Control. He has been a Courtesy Professor in Electrical Engineering and Computer Science since 2013. His research focuses on stochastic systems and control and probability of its applications.  He is a Life Fellow of IEEE, a SIAM Fellow and an IFAC Fellow and was a Simons foundation Fellow.  He was a recipient of the SIAM Reid Prize.
   \end{IEEEbiography}

\mbox{~}\\
  \begin{IEEEbiography}[{\includegraphics[width=1in,clip,keepaspectratio]{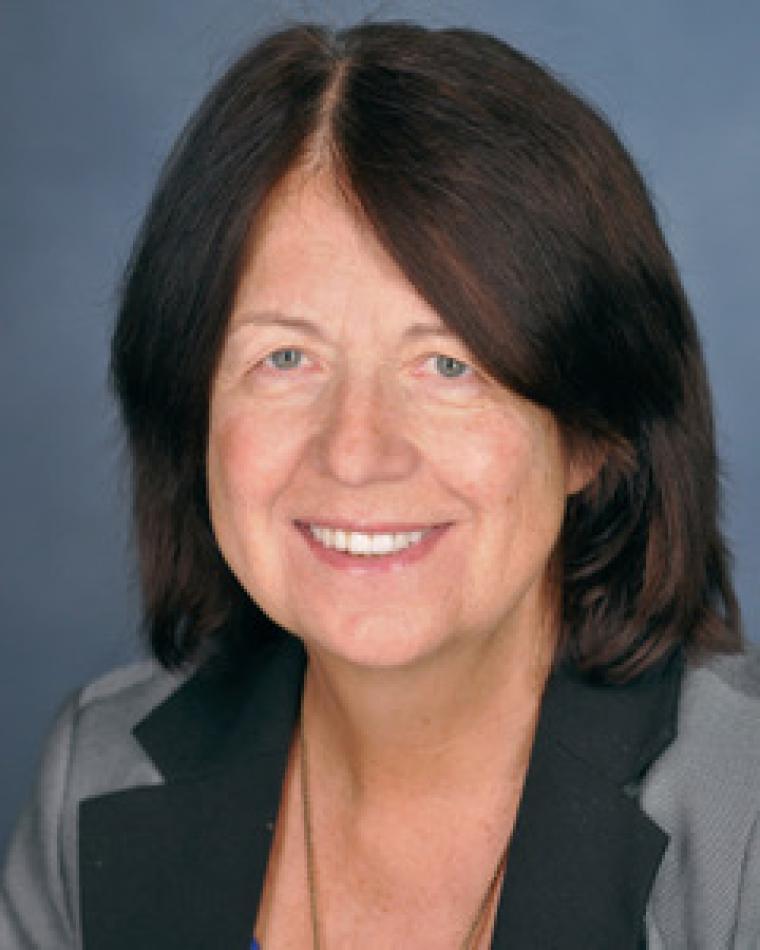}}]{Bozenna Pasik-Duncan} received Master's degree in Mathematics from University of Warsaw, Ph.D. and D.Sc. degrees from Warsaw School of Economics (SGH) in Poland. Before moving to University of Kansas (KU) in 1984, she was faculty member of the Department of Mathematics at SGH. At KU she is Professor of Mathematics, Courtesy Professor of  EECS and AE,  ITTC Investigator, CBC affiliate, Chancellors Club Teaching Professor, and KU Women's Hall of Fame inductee. She is Life Fellow of IEEE and Fellow of IFAC,  recipient of the IEEE Third Millennium Medal and IEEE CSS Distinguished Member Award. She has served in many capacities in several societies including IEEE CSS Vice President, PTM- Warsaw Branch Vice President,  IEEE CSS and SSIT BoG member, Program Director of SIAM Activity Group on CST,  Chair of IEEE CSS TC on Control Education, Chair of  AACC Education Committee, member of  IFAC TB, Chair of IFAC Task Force on DI, Past Chair of IEEE WIE Committee, co-founder of  IEEE CSS Standing Committee on Women in Control (WIC) and its first Chair, member of Award Boards of IFAC and AWM, founder and faculty advisor of AWM and SIAM Student Chapters at KU, founder and coordinator of Mathematics and Statistics Awareness Month (MSAM) / Outreach Program at KU, and founder and Chair of Stochastic Adaptive Control Seminar at KU. She is an Associate Editor of several Journals, and author and co-author of over 200 technical papers and book chapters. Her research is interdisciplinary and primarily in stochastic adaptive control, data analysis and modeling,  and in STEM education. She is a recipient of many awards including IREX Fellow, NSF Career Advancement, KU Women of  Distinction,  Kemper Fellow, IFAC Outstanding Service, Steeples Service to Kansans, H.O.P.E. , AWM L. Hay, Morrison, Max Wells, B. Price, Polish Ministry of Higher Education and Science, IEEE EAB Meritorious Achievement Award in Continuing Education. Her hobbies are music, art and poetry. She loves traveling. She has visited 55 countries. She is passionate about helping others. Her husband, Tyrone Duncan, IEEE Life Fellow is an EE, and professor of mathematics; her daughter, Dominique Duncan, IEEE SM is an EE, mathematician, and associate professor of neurology and bioengineering.
 \end{IEEEbiography}

\mbox{~}\\
  \begin{IEEEbiography}[{\includegraphics[width=1in,clip,keepaspectratio]{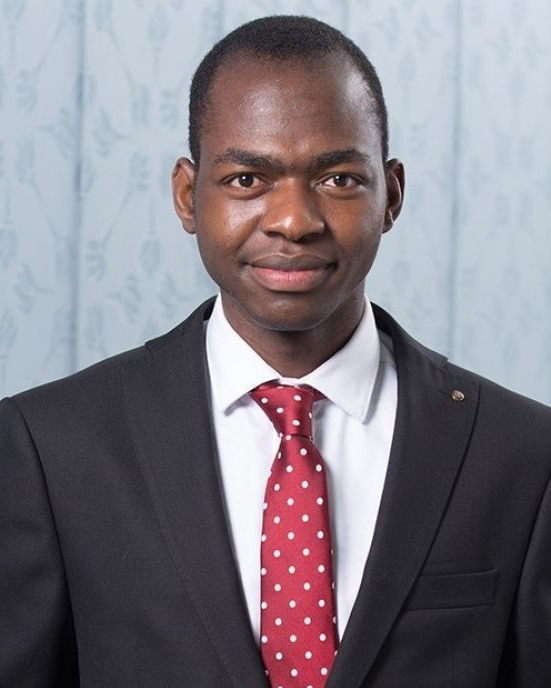}}]{Hamidou Tembine} is the co-founder of Timadie, Grabal  AI Mali,  co-chair of TF, founder of Guinaga, WETE, MFTG, LnG Lab, CI4SI and a Professor of Artificial Intelligence at UQTR, Quebec, Canada. He graduated in Applied Mathematics from Ecole Polytechnique (Palaiseau, France) and received the Ph.D. degree from INRIA and University of Avignon, France. He further received his Master degree in game theory and economics. His main research interests are learning, evolution, and games. In 2014, Tembine received the IEEE ComSoc Outstanding Young Researcher Award for his promising research activities for the benefit of society. He was the recipient of 10+ best paper awards in the applications of game theory. Tembine is a prolific researcher and holds 300+ scientific publications including magazines, letters, journals and conferences. He is author of the book on ``distributed strategic learning for engineers” (published at CRC Press, Taylor \& Francis 2012) which received book award 2014, and co-author of the book ``Game Theory and Learning in Wireless Networks” (Elsevier Academic Press) and co-author of the book on "Mean-Field-Type Games for Engineers".  Tembine has been co-organizer of several scientific meetings on game theory in water, food, environment, networking, wireless communications and smart energy systems. He has been a visiting researcher at University of California at Berkeley (US), University of McGill (Montreal, Quebec, Canada), University of Illinois at Urbana-Champaign (UIUC, US), Ecole Polytechnique Federale de Lausanne (EPFL, Switzerland) and University of Wisconsin (Madison, US). He has been a Simons Participant and a Senior Fellow 2020. He is a senior member of IEEE. He is a Next Einstein Fellow, Class of 2017.
 \end{IEEEbiography}

\end{document}